\documentclass[12pt,notitlepage]{article}
\usepackage{amsmath}
\usepackage{amssymb}
\usepackage{amscd}
\usepackage{amsfonts}
\usepackage{amsthm}
\usepackage{mathrsfs}
\usepackage[matrix,arrow,curve]{xy}

\usepackage{mathtext}         
\usepackage[T2A]{fontenc}     
\usepackage[cp1251]{inputenc} 

 \theoremstyle{definition}
  \newtheorem{df}{Definition}
 \theoremstyle{plain}
  \newtheorem{theorem}{Theorem}
  \newtheorem{prop}[theorem]{Proposition}
  \newtheorem{lemma}[theorem]{Lemma}
  
 \theoremstyle{remark}
  \newtheorem{ex}{Example}
  \newtheorem{rem}{\bf Remark}

\newcommand{\Dj}{\hbox to 8pt{\raisebox{.4\height}{-}\hss D}}

\newcommand{\Hom}{\ensuremath{\mathop{\rm Hom}\nolimits}}

\newcommand{\bc}{\ensuremath{\mathcal B}}
\newcommand{\ac}{\ensuremath{\mathcal A}}
\newcommand{\cc}{\ensuremath{\mathcal C}}

\newcommand{\fc}{\ensuremath{\mathcal F}}
\newcommand{\hc}{\ensuremath{\mathcal H}}
\newcommand{\ic}{\ensuremath{\mathcal I}}
\newcommand{\lc}{\ensuremath{\mathcal L}}
\newcommand{\uc}{\ensuremath{\mathcal U}}
\newcommand{\vc}{\ensuremath{\mathcal V}}
\newcommand{\tc}{\ensuremath{\mathcal T}}
\newcommand{\pc}{\ensuremath{\mathcal P}}

\newcommand{\sca}{\ensuremath{\mathcal S}}
\newcommand{\nc}{\ensuremath{\mathcal N}}

\newcommand{\g}{\ensuremath{\mathfrak{g}}}

\newcommand{\eqdef}{\stackrel{\rm def}{=}}

\newcommand\beq{\begin{equation}}
\newcommand\eeq{\end{equation}}
\newcommand\bqa {\begin{eqnarray}}
\newcommand\eqa {\end{eqnarray}}

\newcommand{\bear}{\begin{array}}
\newcommand{\enar}{\end{array}}

\newcommand{\DR}[1]{\Omega_{DR}({#1})}

\newenvironment{Proof}{\noindent{\it Proof.\ }}{\hfill$\square$}

\newcommand{\inv}[1]{{#1}^{-1}} 
\usepackage{graphicx}
\usepackage{wrapfig}

\begin{document}
\begin{titlepage}
\hfill ITEP-TH-32/08 \vskip 2.5cm

\begin{center}{\LARGE \bf Holonomy, twisting cochains \\
and characteristic classes}
\end{center}

\vskip 1.0cm \centerline{ G.Sharygin \footnote{E-mail:
sharygin@itep.ru} }

\centerline{\sf Institute of Theoretical and Experimental Physics
\footnote{ITEP, 25 ul. B. Cheremushkinskaya, Moscow, 117259 Russia}}

\vskip 2.0cm

\centerline{\large \bf Abstract} The primary interest of this paper is to discuss the r\^ole of twisting cochains in the theory of characteristic classes. We begin with the homological description of holonomy map, associated with a connection on a trivial bundle over a 1-connected manifold. We regard it as a homomorphism from the algebra of differential forms on the structure group to the algebra of differential forms on the based loopspace of the base, represented by the (reduced) bar-complex of differential forms on it. Next we discuss the notion of ``twisting cochains'', or more generally ``twisting maps'', their equivalence relation and give various examples. We show that every twisting map gives rise to a map from the coalgebra to the bar-resolution of the algebra. Further we show that in the case of genuine twisting cochains one can obtain a map from the differential forms on the gauge bundle, associated with the given principal one, to the reduced Hochschild complex of the algebra, which models the base. Then we discuss a concrete example of a twisting cochain that is defined on the polynomial de Rham forms on an algebraic group and takes values in \v Cech complex of the base. We show how it can be used to obtain explicit formulas for the Chern classes. We also discuss few modifications of this construction. In the last section we discuss the construction, similar to the one, used by Getzler, Jones and Petrack in their paper \cite{GJP}. We show that the map we call ``Getzler-Jones-Petrack's map'' is homotopy-equivalent to the map that one obtains from a twisting cochain. This enables us to find a generalization of the Bismut's class, which we regard as an image of a suitable element in the differential forms on the group under the Getzler-Jones-Petrack's map.
\vskip 1.0cm

\end{titlepage}
\tableofcontents
\newpage
\section*{Introduction}
This paper is a result of an attempt to give an algebraic description of the well-known results of Witten, Bismut, Getzler, Jones, Petrack and others, which describe the index of a Dirac operator on a vector bundle in terms of pairing of two characteristic classes in equivariant homology and cohomology of the free loop space of the manifold. In effect, the following formula is true:
$$
ind_ED=\langle ch(D),\,ch_E(\nabla)\rangle,
$$
where $D$ is a Dirac operator on the spin manifold $X$ and $E$ is a vector bundle on $X$. The class $ch(D)$ has degree $0$, it is the so-called Witten's class in equivariant (with respect to the natural action of the circle) homology of the free loop space of $X$, $\mathcal LX$. It is given by a formal expression, involving Feynman integral over the space $\mathcal LX$. The nature of this object is somewhat mysterious and far from being completely understood. On the right hand side there stands another characteristic class, Bismut's class, determined by a connection on the bundle $E$. It belongs to the degree $0$ equivariant \emph{co}homology of $\mathcal LX$. In the paper \cite{GJP} this class was described as a cocycle in the (reduced) cyclic complex of the de Rham algebra of $X$. In the present paper we develop further the ideas of \cite{GJP} and try to express $ch_E(\nabla)$ in the terms of such a fundamental algebraic topological object as \emph{twisting cochain}. We hope that it will be possible to apply similar ideas to the left side of this formula and obtain an algebraic description of Witten's class too. However, we postpone this discussion to a future paper.

Let us recall that the notion of twisting cochains was first introduced in algebraic topology by E.~Brown (see \cite{Br59} and definition \ref{twcochdf} below). In the cited paper of Brown it is shown that for every principal bundle $P$ over a space $X$ with structure group $G$ (we shall always assume that all three spaces are compact closed manifolds and that $X$ is $1$\/-connected) there exists a twisting cochain on the coalgebra of singular chains on the base with values in the Pontriagin algebra of sigular chains on the structure group, such that the corresponding twisted tensor product (see equation \eqref{eqtwdiff}) models the total space of the fibre bundle. Moreover, it can be shown that this correspondence is 1-1 modulo an equivalence relation, similar to the relation between flat connections. Thus one can regard twisting cochains as algebraic counterparts of bundles and connections and ask the question: ``Is it possible to express other important invariants of fibre bundles like characteristic classes, corresponding K-theory elements and index of elliptic operators on bundles (in particular, Bismut's and Witten classes) in terms of their twisting cochains?'' In the present paper we answer part of this questions and present explicit constructions for twisting cochains -- something that is rarely found in papers on algebraic topology. One should keep in mind, that in the cited paper of Brown, the author works in the homological setting, i.e. differentials in all algebras and coalgebras that he uses decrease degrees by $1$, while in our paper we prefer to deal with cohomological situation. It results in the slight change of definition of the twisting cochain and in few other minor corrections. In particular, we cannot use Brown's theorem directly, instead of this we give an explicit construction of twisting cochain in the situation we consider.

Thus we regard the notion of twisting cochain as a convenient algebraic model for principal bundles and connections. This approach enables us to extend the ideas and methods of differential geometry to the wider domain of differential graded algebras and coalgebras, that are not necessarilly de Rham algebras of a manifold. One can work with twisting cochains pretty much like one works with connections (flat or non-flat). In particular, one can ask, if it is possible to use a twisting cochain to find characteristic classes of the principal bundle, especially when it is difficult or impossible to find a connection. The answer to this question is positive: it is well-known that twisting cochains induce characteristic maps from the cobar resolution of the coalgebra to the algebra in which they take values (see \cite{Smirn76, Smirn77, SmSOp}). One can show that the image of this map consists of characteristic classes of the bundle $P$. We use this simple idea and the explicit formula for the twisting cochain (with values in the \v Cech complex of the base) to find explicit formulae for the characteristic classes of $P$ expressed in the terms of the glueing functions (i.e. in the terms of the corresponding noncommutative \v Chech cocycle with values in the structure group.) For the sake of simplicity we shall always assume that our structure group is embedded into a general matrix group over $\mathbb C$ or $\mathbb R$. It doesn't cause any loss in the generality of our consideration, since as it is well-known all smooth compact Lie groups are subgroups of $GL(n)$. Details of this construction one can find in the paper \cite{mypapLF}, where we compare our approach to the subject with that of Bott (see, e.g. \cite{BottTu}, section 12.)

Another natural question that one can ask about the twisting cochain is concerned with the holonomy map that is usually related to a connection on a bundle. One can ask, if it is possible to express this map in the terms of twisting cochain. This question is very natural in the framework of algebraic topology in the view of the well-known Kan's theorem (see \cite{Kan}): there's a 1-1 correspondence between principal bundles over a base and representations of Kan's group of the base (i.e. of the group, homotopy equivalent to the $H$\/-space of based loop on the base.) There are many models of this group, and many algebraic constructions that can replace its singular cochains. One of the most convenient of them is the Chen's iterated integral construction that gives an insight into the structure of this space. This construction gives a homomorphism of algebras from bar-resolution of the de Rham algebra of a manifold $X$ (one can introduce the multiplication in $B(\DR{X})$ with the help of shuffle products) into the algebra of differential forms on the loop space of $X$. One can show (see \cite{Chen}) that this map induces an isomorphism in cohomology if $X$ is 1-connected.

This approach to the loop spaces can be generalized to free loop space of a manifold: it turns out that the iterated integral map can be extended to the map from (reduced) Hochschild complex of $\DR{X}$ into the de Rham forms on the free loop space of $X$ (see \cite{GJP} and \cite{Jones},) which induces an isomorphism on cohomologies under the same condition of 1-connectedness of $X$. Moreover, in the cited papers it is also proved that the equivariant cohomology of free loop space with respect to the cirle action on it (by $S^1$ translations of the argument) is isomorphic to the cyclic homology of the algebra $\DR{X}$ (the map is given by a generalization of the iterated integral). Thus it is natural for us to use the bar-complex, Hochschild and cyclic complexes as models for the loop spaces. It is in one of these complexes, where our map should take values.

Here we should make a remark, concerning the terminology used in this paper: we regard the loop space rather as a topological group. Thus the holonomy of a connection becomes for us a homomorphism of groups. This justifies the name of \textit{homological monodromy map} that we use to describe the inverse image homomorphism, induced by the holonomy of a connection and its purely algebraic analogs that we construct in this paper.

In what follows we shall use twisting cocains to construct a map from a complex, modelling the gauge bundle associated with the given principal one, to the Hochschild complex of the base (recall, that one calls the gauge bundle of a principal bundle $P$ with structure group $G$ the following associated space $\hat P=P\times_{Ad}G$, where $Ad$ denotes the adjoint action of group on itself). We show that it is in fact homotopy equivalent to the homological monodromy map. In effect, the complex that shall play the r\^ole of algebraic model of the gauge bundle is what one can call a ''bitwisted`` tensor product $K\hat{\otimes}_\phi A$ of algebra and coalgebra (see definition \ref{bitwtenpr}). We show that there's a map $\tilde{\hat\phi}$ from $K\hat{\otimes}_\phi A$ into the (reduced) Hochschild complex of $A$. Moreover, if the algebra $A$ modelling $X$ is commutative, then this map intertwines the comultiplication structures on both sides (see proposition \ref{propleftright} and discussion that follows it,) otherwise one should use the $A_\infty$ structure on the left.

On the other hand, if we consider frame bundle $P$ of a vector bundle $E$ of rank $n$ (i.e. $P$ is a principal $GL(n)$\/-bundle,) one can embed $E$ into a trivial rank $N$ bundle ($N>n$) and consider the corresding projector $p$. In this way one obtains a globally defined connection (''Grassmanian connection``) on the bundle $E$. Observe that the connection form in this case is a global $N\times N$ matrix-valued $1$\/-form on $X$. One can use this connection in a way, similar to the trivial case, considered in section \ref{sectbegin} and obtain a map from the differential forms on $\hat P$ to differential forms on the free loop space of $X$, or rather to the reduced Hochschild complex of $\DR{X}$. We shall call this latter map ''Getzler-Jones-Petrack's`` map: it is a homomorphism of the Hopf algebras over $\DR{X}$. We prove (see theorem \ref{omagn}) that this map is homotopy equivalent to the map $\tilde{\hat\phi}$ that we constructed before in a purely algebraic way.

This result is interesting because Bismut's class as described by Getzler, Jones and Petrack can be obtained as the result of application of Getzler-Jones-Petrack's map (or rather of its suitable modification) to an explicit element in the de Rham algebra of $\hat P$. The modification of the Getzler-Jones-Petrack's map is necessary because Bismut's class takes value in equivariant cohomology of $\mathcal LX$. It turns out that it is possible to modify the map $\tilde{\hat\phi}$ in a similar way, so that its values on cocentral elements in $\DR{GL(n)}$ (or rather on similar elements in $\DR{\hat P}$) are closed in the corresponding cyclic complex. Unfortunately one cannot prove directly that these elements are equivalent to Bismut's class, because they all depend on the choice of twisting cochain, just like Bismut's class (and Chern-Symon's classes) depend on the choice of flat connection. However, if the twisting cochain is determined by a connection $\nabla$ (see section \ref{sectconnect},) then it is evident that our construction gives the same class as that of Bismut. In a general case one can use the reasoning similar to the one used in the end of section \ref{sect3} to prove its independence on the choice of the twisting cochain inside the given equivalence class of twisting cochains.

\subsection*{Contents and notations}
The rest of the paper is divided into four big sections, each containing several subsections. In the first section we discuss the most simple case: the case of the trivial bundle. We describe the monodromy map induced by a connection on such a bundle in terms of the Chen's iterated integral map as a homomorphism from the de Rham algebra of the structue group $G$ to the (reduced) bar complex of $\DR{X}$ (here $X$ is the base of the princiapl bundle.) So in subsection \ref{sectbegin} we recall the construction of iterated integral map and the description of geometric monodromy in terms of time-ordered exponent. The next subsection \ref{sect2} is devoted to a detailed construction of homological monodromy map in the terms of the given global connection (this construction will be later used in the chapter, devoted to Getzler-Jones-Petrack's map). In the next subsection, \ref{sect3}, we show that the map we construct doesn't depend on the gauge transformation class of the connection, at least if the range of this map is the \emph{normailized} bar-resolution of the de Rham algebra of the base. After this, in the last subsection of this section, paragraph \ref{sect17}, we discuss a variant of the same monodromy map, this time defined on the algebra of differential forms on the gauge bundle and taking values in the (normalized) Hochschild complex of the base.

In the next section we review the general notion of twisting cochains, their gauge transformations characteristic maps, induced by twisting cochains and other constructions, related to them. Thus, in subsection \ref{sectnext} we recall the definition of the twisting cochain and give a slightly more general definition of a twisting map. Then we give several important examples of such maps that appear in various geometric situations (first of all in the case of trivial bundle over a manifold). Then in subsection \ref{sectnext2} we give the definition of characteristic map, associated with a twisting map $\phi$ from a coaugmented d.g. coalgebra $K$ to a d.g. $A$ (in particular, with twisting cochain.) By definition it is a map from $K$ into the bar-resolution of the algebra. Dually one can construct the map from the cobar resolution of $K$ into $A$. We show that these maps commute with differentials and discuss the conditions under which they become homomorphisms of coalgebras (respectively, of algebras.)  The subsection \ref{sectnnext} is devoted to establishing equivalence relations on twisting maps: we define the equivalence transformations of twisting maps, which we call ''gauge transformations``, and which in the case of twisting cochains reduce to the standard definition. Then we discuss few examples of such equivalences, in particular we show that the cochains on trivial bundles from subsection \ref{sectnext} are all trivial up to a gauge transformation. Besides this we show that the characteristic maps of the previous paragraph don't depend (up to a chain homotopy) on the choice of equivalent twisting maps. Finally, in the last paragraph, of this section, subsection \ref{sectnext4}, we discuss the generalization of the characteristic map and of twisted tensor product. Namely, we define the ''bitwisted`` product $K\hat\otimes_\phi A$ (def. \ref{bitwtenpr}) and construct the map $\tilde{\hat\phi}$ from this product into the (reduced) Hochschild complex of $A$ (at this stage we do not need reduction, so its accurate definition is postponed to the last section of this paper).

Next section, section \ref{sect30} is devoted to the study of twisting cochains, associated to a nontrivial principal bundle. So, in paragraph \ref{sect31} we construct an example of such a twisting cochain, it is defined on the coalgebra (in fact, Hopf algebra) $K=\DR{GL(n)}$ of (polynoial) de Rham forms on the group and takes values in the algebra $A$ of \v Cech-de Rham forms on the base. The main idea that we use in this construction is founded on the results of paper \cite{Shih}. We find the explicit expression of the twisting cochain in the terms of the transition functions $g_{\alpha\beta}$ of the given principal bundle and check, that this formula actually verifies the definition of twisting cochains and that the associated twisted tensor product is isomorphic to an appropriate \v Cech-de Rham complex of the total space. In section \ref{sect32} we plug in this formula into the construction of the characteristic map of previous section (here we use the map from cobar resolution of $K$ into $A$) to obtain explicit formulas for the characteristic classes of the bundle. In fact, we show, that all the characteristic classes of $P$ can be obtained in this way. Details of this construction, as well as an explicit algorithm for calculating Chern classes and its comparison with earlier results of other authors can be found in our paper \cite{mypapLF}. The last subsection of the section \ref{sect30}, paragraph \ref{sectconnect}, is devoted to another example of the twisting cochain associated to the bundle $P$. Once again the range of our cochain is in the \v Cech-de Rham complex of the base. This construction is based on a choice of global connection on $P$ and uses such ingredients as the connection $A$, its curvature $F$ and the transition functions $g_{\alpha\beta}$. We show, that this twisting cochain is in effect equiuvalent to the one we used before. One of the advanteges of this new formula is that it is explicitly related to the differential geometry of the principal bundle. In particular, one can see that the characteristic map of this twisting cochain is related to the monodromy map of the nontrivial bundle $P$.

After this in section \ref{seccc} we discuss possible approaches to the ''globalization`` of our twisting cochains. So, in subsection \ref{sectinfty} we describe a way to ''glue`` the formulas, given in the sections \ref{sect31} and \ref{sectconnect} to obtain a cochain on $K$ with values in the de Rham algebra of the base. To this end we recall the perturbation lemma and use it to define the homotopy inverse of the evident homotopy equivalence from the de Rham algebra to the \v Cech-de Rham algebra of a manifold. Next we describe the $A_\infty$\/-map that intertwines the multiplication on these two complexes. Finally we use the formulas from \cite{SmSOp} that relate $A_\infty$ maps and  twisting cochains to define the final result. Unfortunately, this construction is not very useful in a general case, since it uses extensively such an implicit object as the partition of unity, associated with a given open cover. After this, in sections \ref{sect41} and \ref{sect42} we work in another possible global algebra, associated to the given space: the Dupont realization of the algebra of de Rham forms on the open subsets of the base. First of all, we define this algebra and prove that it is homotopy equivalent to the de Rham algebra of a manifold. After this, in section \ref{sect42} we show how one can obtain a twisting cochain with values in this algebra.

The last section of this paper, section \ref{sect40} is devoted to the study of Bismut classes of the bundle. First, in paragraph \ref{sect43} we use Grassmanian connection on a vector bundle to define a map $\tilde\phi_E$ from the de Rham algebra of the gauge bundle $\hat P_E$, associated with the bundle to the Hochschild complex of the de Rham forms of the base. The construction of this map is similar to the constructions, we use to define the monodromy map in section \ref{sectbegin}, \ref{sect17} and resembles the construction of Bismut class due to Getzler, Jones and Petrack, see \cite{GJP}. Next we show, section \ref{sectcompth}, theorem \ref{omagn}, that this map is homotopy equivalent (as a map of coalgebras over the de Rham algebra of the base) to the characterisc map $\tilde{\hat\phi}$, earlier associted to a twisting cochain (with values in $\DR{X}$, see section \ref{sectinfty} for the definition of the twisting cochain.) Finally, in section \ref{sect44} we describe the equivariant complex of Jones Getzler and Petrack and show that their Bismut class is an image of the ''equivariantized`` characteristic map $\tilde\phi_E$ from previous subsection. After this we show, what modifications should be made in the definition of the characteristic maps $\tilde{\hat\phi}$ to obtain similar classes from twisting cochains.

The following conventions are used throughout the text. All the algebras and coalgebras that we consider are taken over a fixed characteristic $0$ field $\Bbbk$ ($\mathbb C$ or $\mathbb R$ will do,) unless otherwise stated. Throughout the paper we work with a fixed principal bundle $P$ over a base $X$. We assume that both $P$ and $X$ are closed manifolds, that $X$ is compact and that the projection is smooth submersion. The structure group $G$ of this principal bundle is supposed to be an algebraic subgroup of a general linear group (in particular, it can be equal to $GL(n)$ itself.) If $G$ is compact, then we assume that $P$ is also compact. In the most part of the paper we assume that $G$ acts on $P$ \emph{on the left}, which is different from the usual agreements, although in many occasions, where it doesn't cause confusion, we do not distinguish between the left and the right actions. We use the symbol $\mathfrak g$ for the Lie algebra of the structure group $G$ (we regard it as the space of all left-invariant vector spaces on $G$, or equivalently as the tangent space in the unit element of $G$.) Throughout the text $\DR{\cdot}$ will denote the deRham algebra of a smooth manifold, while the symbol $\Omega X$ will be reserved for the based loops space of a based space $(X,\,\ast)$ (we shall usually omit the base point $\ast\in X$ from our notation,) i.e. for the space of all maps $f:S^1\to X,\ f(1)=\ast$. If we want to underline the particular choice of the base point $x\in X$, we shall use the symbol $\Omega_x X$. By $\mathcal LX$ we shall denote the free loops space of a space $X$, i.e. the space of all maps $S^1\to X$ without any restriction on them. If $X$ is a smooth manifold, then we assume that all the loops in $\Omega X$ and $\mathcal LX$ are piecewise-smooth. The same smoothnes assumption will be imposed on all the mapping spaces that appear in the text. The symbol $\Omega^*_{poly}(\cdot)$ will refer to the differential graded algebra of polynomial differential forms on a (smooth) affine variety, e.g. on the group $G$. Moreover we shall often abbreviate $\Omega_{poly}^*(G)$ simply to $\Omega_G^*$ or $\Omega_G$. In particular, $\Omega_{poly}^0(\cdot)$ will denote the algebra of polynomial functions. We shall often abbreviate $\Omega_{poly}^0(G)=\Omega_G^0$ simply to $\ac(G)$. All differntials we consider here will increase the degree of an element by $1$, i.e. we shall work in cohomological context. For a graded algebra $\Omega$, symbol $|w|$ will denote the degree of a homogeneous element $w\in\Omega$. All the other notations will be explained in appropriate places of the text.

\subsection*{Acknowledgements}
The author would like to express his deepest gratitude to Anton Gerasimov, whose attention and numerous discussions at the early stages of this project helpded the author a great deal. I would also like to express my gratitutde to all the institutions where I payed visits while working on this project, especially to IMPAN, Warsaw, where I took part in the programm  ''Noncommutative Geometry and Quantum groups``, contract No MKTD-CT-2004-509794. While working on this project I was supported by the grants NSh-3035.2008.2 and RFFI 09-01-00239a.

\section{The homological monodromy map}
\label{sect11}
In the first two subsections of this section we define the \textit{homological monodromy map} in the case of a trivial principal $G$\/-bundle over a pointed space $(X,\,x_0)$. In this case it is a homomorphism of differential graded Hopf algebras $\Omega_{poly}^*(G)\to B(\Omega_{DR}(X))$, which arises from the time-ordered exponent. Further we discuss the way it commutes with the action of the gauge transformations and in the last paragraph (section \ref{sect17}) we briefly discuss its possible extension to the gauge bundles, associated with arbitrary principal bundles. The latter question will be examined in more detail in the last chapter.

\subsection{Time-ordered exponent and iterated integrals}
\label{sectbegin}
Let $G$ be a compact smooth semisimple algebraic group (over $\mathbb C$ or $\mathbb R$). We can assume, that it is imbedded in the standard unitary (orthogonal) group $U(n)$ (resp. $O(n)$). Let $\mathfrak g$ denote its Lie algebra, $N={\rm dim\/}\mathfrak g$. In the complex case the algebra of regular functions on $G$, $\ac(G)=\Omega_{poly}^0(G)$ can be described as a factor algebra of the polynomial ring of $2n^2$ variables $u_{ij}$ and $\bar u_{ij},\ i,j=1\dots n$ modulo the шideal, generated by the equations that determine the group. This algebra is equipped with antilinear involution $u_{ij}\mapsto\bar u_{ij}$. One should regard $u_{ij}$ as the function, which assigns to an element $g\in G$ the value of the $(i,j)^{{\rm th}}$ entry of the matrix, representing it in $U(n)$. (Evidently, in real case one should use only the first half of these generators and involution is given simply by complex conjugation of coefficients.) The $\mathbb C$-linear homomorphisms
\begin{align}
\label{eq1}
\Delta(u_{ij})&=\sum_k u_{ik}\otimes u_{kj}\\
\label{eq2}
\epsilon(u_{ij})&=\delta_{ij},\\
\intertext{where $\delta_{ij}$ is the Kronecker symbol, and $\mathbb C$-antilinear (anti-)homomorphism}
\label{eq3}
S(u_{ij})&=\bar u_{ji}
\end{align}
determine the Hopf algebra structure on $\ac(G)$. Besides this it is known (by Peter-Weyl theorem), that linearly $\ac(G)$ is equal to the direct sum
\begin{equation}
\label{eq4}
\ac(G)=\bigoplus_{\rho\in{\mathcal T}}\ac^\rho.
\end{equation}
Here $\mathcal T$ denotes the set of all irreducible unitary representations $\rho:G\to U(n_\rho)$ of the group $G$, and $\ac^\rho$ is the subspace, generated by functions $u^\rho_{ij}(g)=\rho(g)_{ij},\ i,j=1,\dots n_\rho$. Observe, that the spaces $\ac^\rho$ are closed under comultiplication, since
\begin{equation}
\label{eq5}
\delta(u^\rho_{ij})=\sum_k u^\rho_{ik}\otimes u^\rho_{kj}.
\end{equation}

Now let $P\xrightarrow[G]{\pi} X$ be a (not necessarily trivial) principal bundle with the structure group $G$ over a smooth manifold $X$. Recall that a connection form $A$ on $P$ is a $\mathfrak g$-valued $1$-form on $P$, satisfying certain covariance and normalization conditions. Let $h_1,\dots h_N$ be the orthonormal (with respect to the Killing form) basis of $\g$. Then we can write $A$ down as $A=\sum_{k=1}^N A_k\,h_k,\ A_k\in\Omega^1(P)$.

For any local section $s:U\to P$, on an open subset $U\subseteq X$, one can consider the inverse image $s^*(A)=\sum_k s^*A_k\,h_k$ of $A$. One can define connection $A$ as a collection of such locally defined $\g$-valued 1-forms on $X$, also known as \textit{gauge potentials}. If $\rho$ is a representation of $G$, then we can consider the collection of local matrix-valued gauge potentials $A^\rho=\sum_k A_k\,\rho(h_k)$. We will omit the section and use the same notation for a representation of group and for the induced representation of Lie algebras. In this way one defines a connection on the associated complex vector bundle $E^\rho=P\times_\rho\mathbb C^{n_\rho}$.

Let $\gamma:[0;1]\to X$ be a path on $X$. One can define the holonomy transformation along this path, associated to the connection $A^\rho$ as follows. A local section $\tilde s$ of $E^\rho$ is contravariant constant along $\gamma$ with respect to the connection $A$, if its covariant derivative along $\gamma$ with respect to $A$ vanishes. This means that $s$ is flat in the direction of $\gamma$, if the following linear differential equation holds:
\begin{equation}\label{eq5 1/2}
  \nabla_{\frac{d}{dt}}(\tilde s)=\frac{d}{dt}\tilde s+\sum_k
A_k(\frac{d}{dt})\rho(h_k)(\tilde s)=0.
\end{equation}
Here we regard $\tilde s$ as a $\mathbb C^{n_\rho}$-valued function on $[0;1]$ (in fact, $\tilde s(t)\in E^\rho_{\gamma(t)}$, and any vector bundle over a segment is trivial.) For a vector $v\in E^\rho_{\gamma(0)}$ one can consider the solution $\tilde s_v$ of equation \eqref{eq5 1/2} with the initial value $\tilde s_v(0)=v$. Then the holonomy is defined as the map, sending $v$ to $\tilde s_v(1)$. If we chose bases in the fibres of $E$ at $\gamma(0)$ and $\gamma(1)$, we shall obtain a matrix $M_A\in \rho(G)$, the subscript $A$ stands for the connection we choose, which changes to $g^{-1}M_A h$, if the bases are changed by matrices $g$ and $h$ in $E^\rho_{\gamma(0)}$ and $E^\rho_{\gamma(1)}$. The map $\gamma\mapsto M_A(\gamma)$ satisfies the relation $M_A(\gamma_1 *\gamma_2)=M_A(\gamma_1)M_A(\gamma_2)$ when the paths $\gamma_i$ are composable, and hence it defines a ``group homomorphism'' from $\Omega_{x_0}X=\Omega X$ to $\rho(G)$ (recall, that $\Omega X$ is not a group or even a monoid itself, but only homotopy equivalent to a group, hence it is not quite correct to speak about homomorphism here).

From now on and through the end of this section, we will suppose that the principal bundle $P$ is trivial. (The general case will be briefly discussed in paragraph \ref{sect17}.) Choose a global section $s:X\to P$ (and hence the global trivialization of $P$, $P=X\times G$ and of all the associated vector bundles.) Then the connection $A$ pulls down to a \g\/-valued differential 1-form on all $X$. Now, it is well known, that one can write down the holonomy matrix $\widetilde M_A$ in the terms of the so-called \textit{time-ordered exponent} as follows:
\begin{equation}
\label{eq6}
M_A(\gamma)=P\exp\int_\gamma A^\rho dt.
\end{equation}
The right hand side of \eqref{eq6} is equal to the infinite sum of Chen's iterated integrals:
\begin{equation}
\label{eq7}
P\exp\int_\gamma A^\rho dt=\sum_{n=0}^\infty \int_{\Delta^n}
A^\rho(t_1)\dots A^\rho(t_n) dt_1\dots dt_n,
\end{equation}
where $\Delta^n=\{(t_1,\dots,t_n)\in\mathbb R^n|0\leq t_1\leq\dots\leq t_n\leq1\}$ is the standard $n-$dimensional simplex. Consider the map $\Psi_\gamma:\Delta^n\to X;\ (t_1,\dots,t_n)\mapsto\gamma(\frac{1}{n}(t_1+\dots+t_n))$. It defines the inverse image map $\Psi_\gamma^*:\Omega^*_{DR}(X)\to\Omega^*_{DR}(\Delta^n)$. If $\omega\in\Omega^1_{DR}(X)$, we define $\omega(t_i)\in C^\infty(\Delta^n)$ as the coefifficients of the decomposition $\Psi_\gamma^*(\omega)=\sum_{i=1}^n\omega(t_i)dt_i$. One defines iterated integral $\int_{\Delta^n} A^\rho(t_1)\dots A^\rho(t_n) dt_1\dots dt_n$ as
\begin{equation}
\label{eq8}
\int_{\Delta^n}
A^\rho(t_1)\dots A^\rho(t_n) dt_1\dots dt_n=\int_0^1 A^\rho(t_1)\Bigl[\int_{t_1}^1 A^\rho(t_2)\Bigl[\int_{t_2}^1\dots\Bigr] dt_2 \Bigr]dt_1
\end{equation}
The holonomy map defines the inverse image homomorphism $M_A^*:\Omega_{DR}(G)\to\Omega_{DR}(\Omega X)$, where the right hand side is defined in terms of the local plots, see \cite{Chen}. By the virtue of the decomposition \eqref{eq4}, it is enough to define the images of $u^\rho_{ij}$ for all $\rho$. Since $\Psi_\gamma^*(A^\rho)=\sum_k\Psi_\gamma^*(A_k)\,\rho(h_k)$, we see that $A^\rho(t_i)=\sum_k
A_k(t_i)\,\rho(h_k)$. So
\begin{equation}
\begin{split}
\label{eq9}
M_A(\gamma)&=\sum_n \int_{\Delta^n}\prod_{i=1}^n\Bigl(\sum_k
A_k(t_i)\,\rho(h_k)dt_i\Bigr)\\
&=\sum_{n, (k_1,\dots,k_n)}\Bigl(\int_{\Delta^n}
A_{k_1}(t_1)\dots A_{k_n}(t_n)dt_1\dots
dt_n\Bigr)\rho(h_{k_1})\dots\rho(h_{k_n})\\
&=\sum_{n, (k_1,\dots,k_n)}\Bigl(\int_{\Delta^n}
A_{k_1}(t_1)\dots A_{k_n}(t_n)dt_1\dots
dt_n\Bigr)\rho(h_{k_1}\dots h_{k_n}).
\end{split}
\end{equation}
Here the summation is taken over all collections $(k_1,\dots, k_n),\ k_i=1,\dots,N$, and we denote the induced representation of the universal enveloping algebra of $\g$ by the same symbol $\rho$. Clearly, the $(i,j)^{\rm th}$ entry of the resulting matrix is equal to the sum of the $(i,j)^{\rm th}$ entries of the terms, so
\begin{equation}
\label{eq10}
M_A^*(u^\rho_{ij})=\sum_{n, (k_1,\dots,k_n)}\Bigl(\int_{\Delta^n}
A_{k_1}(t_1)\dots A_{k_n}(t_n)dt_1\dots
dt_n\Bigr)\rho(h_{k_1}\dots h_{k_n})_{ij}.
\end{equation}
Now recall that the representation of Lie algebra is defined by the formula
\begin{equation}
\label{eq11}
  \rho(X)=\frac{d}{dt}|_{t=0}\rho(\exp(tX)),
\end{equation}
for $X\in\g$. Hence, if we represent the basis $\{h_i\}_{i=1}^N$ by vector fields $\{X_i\}$ on $G$, we obtain for an arbitrary function $f\in\ac(G)$:
\begin{equation}\label{eq12}
  M_A^*(f)=\sum_{n, (k_1,\dots,k_n)}\frac{1}{n!}\Bigl(\int_{\Delta^n}
A_{k_1}(t_1)\dots A_{k_n}(t_n)dt_1\dots
dt_n\Bigr) X_{k_1}\dots X_{k_n}(f)|_e,
\end{equation}
where $e$ denotes the unit element of $G$.

Recall that we have assumed the bundle to be trivial, and fixed the trivialization. In particular, this assumption holds for the pullback bundle $\pi^*(P)$. In this case we can regard $A$ as a $\g$-valued global form on $P$, satisfying certain equivariance conditions. Then the map $M^*_A$ determined by formula \eqref{eq6} takes values in the algebra $\DR{\Omega_{p_0}P}$ where $p_0\in P$ is a point, projecting onto $x_0$. In general we would like to show that in this particular case the image of $M_A^*$ belongs to the subalgebra of the forms, lifted from the loop space of the base. To this end we shall need a better understanding of the map $M^*_A$, so we shall postpone this discussion to the last paragraph of this section.

To go on with our investigation of the holonomy map, we will need a suitable model for the De Rham algebras on loop spaces. One of the most convenient models of such algebras was suggested by Chen, see \cite{Chen}. It is based on the notion of differentiable space and smooth plots rather than the open charts. It is shown in \cite{Chen} that for any manifold $X$ there is a well-defined quasi-isomorphism $\sigma$, called the ``iterated integral map'', $\sigma:B(\Omega_{DR}(X))\to\Omega(\Omega_{x_0}X)$. Here $B(\Omega_{DR}(X))=B(\mathbb C,\Omega_{DR}(X),\mathbb C)$ is the bar-resolution of the algebra of De~Rham fomrs on $X$, defined as
\begin{equation}\label{eq12 2/5}
\sigma([\omega_1|\dots|\omega_n])=\int_{\Delta^n}\omega_1(t_1)\dots\omega_n(t_n)dt_1\dots
dt_n.
\end{equation}
In brief, one can define $\sigma$ as the composite map $p_*\Phi^*_n$ in the diagramm:
\begin{equation}\label{eq12 3/7}
\begin{CD}
    {\DR{\Omega X\times\Delta^n}} @<{\Phi_n^*}<< {\DR{X^{\times n}}} @<\supset<< {\DR{X}^{\otimes n}}\\
    @V{p_*}VV                                           @.                       @.\\
    {\DR{\Omega X},}               @.             {}                  @.   {}
\end{CD}.
\end{equation}
Here $\Phi_n(\gamma,t_1,\dots,t_n)=(\gamma(t_1),\dots,\gamma(t_1))$, $p:\Omega X\times\Delta^n\to\Omega X$ is the projection, and $p_*$ denotes the direct image, i.e. integration along the fibres of the projection.

Below we shall work with the \textit{bar-resolution} $B(\Omega_{DR}(X))$ of $\DR{X}$. It is defined as follows
\begin{align}
\label{eq12 1/2}
B(\Omega_{DR}(X))&=\bigoplus_{n\ge0}\mathbb C\otimes\bigl(\Omega_{DR}(X)[1]\bigr)^{\otimes n}\otimes\mathbb C=\bigl(\Omega_{DR}(X)[1]\bigr)^{\otimes n},\\
\intertext{where $[1]$ denotes the suspension, $\Omega^{n-1}_{DR}(X)[1]\eqdef\Omega_{DR}^n(X)$. Elements from $B(\Omega_{DR}(X))$ are denoted by $[w_1|w_2|\dots|w_n]$. The algebra $\Omega_{DR}(X)$ acts on $\mathbb C$ via $w\cdot1=w(x_0)$, if $w\in\Omega^0_{DR}(X)$ and $0$ otherwise. Now the differential in bar-resolution is given by}
\begin{split}
\label{eq12 2/3}
d_B[w_1|w_2|\dots|w_n]&=\sum_{i=1}^n(-1)^{\varepsilon_{i-1}}[w_1|\dots|dw_i|\dots|w_n]+\sum_{i=1}^{n-1}(-1)^{\varepsilon_i+1}[w_1|\dots|w_iw_{i+1}|\dots|w_n]\\
                      &\quad+w_1\cdot1[w_2|\dots|w_n]+(-1)^{\varepsilon_n+1}[w_1|\dots|w_{n-1}]w_n\cdot1,
\end{split}
\end{align}
where $\varepsilon_k=\sum_{i=1}^k|w_i|$. The first term on the right of \eqref{eq12 2/3} corresponds to the De Rham differential in $\DR{X}$ and is usually denoted by $d_I$. The remaining terms of the formula \eqref{eq12 2/3} are denoted by $d_{II}$.

Recall that the bar-resolution of a {\em commutative\/} DG-algebra has in fact the structure of a bialgebra where the coalgebra structure is defined by
\begin{equation}
\label{eq14}
\nabla([a_1|\dots|a_n])=\sum_{i=0}^n[a_1|\dots|a_i]\otimes[a_{k+1}|\dots|a_n],
\end{equation}
and the algebras structure is given by the shuffle-products see \cite{McL}, Ch.8, for details. Moreover, in papers \cite{GJP, Chen} it is shown that the kernel of the Chen's iterated integral map $\sigma$ contains the subcomplex $\ker'\sigma$ of $B(\Omega_{DR}(X))$, generated by the elements
\begin{equation}\label{eq22}
\{[w_1|\dots|w_n]\in B(\Omega_{DR}(X))|\exists i=1,\dots,n,\
w_i\in\Omega^0_{DR}(X)=C^\infty(X)\}.
\end{equation}
Under the assumptions we have adopted in this paper, this subcomplex is acyclic, see \cite{GJP} for details. We put
\begin{equation}\label{eq17 6/7}
N(\Omega_{DR}(X))\eqdef B(\Omega_{DR}(X))\bigl/\ker'\sigma.
\end{equation}
The factor-complex complex $N(\Omega_{DR}(X))$ is a Hopf algebra and it is called the {\em normalized bar-resolution\/}, and $\ker'\sigma$ is the normalization kernel. 

\subsection{Homological Monodromy Map}
\label{sect2}
In what follows we shall primarily deal with the normalized bar-resolution of $\DR{P}$. However we shall not make any difference in the notation between it and the usual bar-resolution, unless it is necessary. We shall also formulate some of our results in the terms of usual (not normalized) bar-resolution, when it is possible.

Our goal is to give a description of the inverse image map, induced by the holonomy of a principal bundle. First of all, taking inverse image of the Chen's iterated integral map $\sigma$ we can regard $M_A^*$ as a map from $\Omega_{poly}^0(G)=\ac(G)$ to the bar-resolution of the base $X$ of the trial bundle $P$. In particular, we can take the bundle $\pi^*(P)\to P$, where $P\xrightarrow[G]{\pi}X$ is an arbitrary principal bundle. Then
\begin{equation}
\label{eq13}
f\mapsto\widetilde M_A(f)=\sum_{n,(\bar k)}(X_{k_1}\dots X_{k_n}(f)|_e)[A_{k_1}|\dots|A_{k_n}].
\end{equation}
\begin{prop}
\label{prop1}
The map $\widetilde M_A$ is a homomorphism of bialgebras.
\end{prop}
\begin{Proof}
First of all, recall (refer to the cited book by McLane,) that comultiplication and shuffle product induce the commutative DG bialgebra structure on the bar-resolution of a commutative DG algebra. Consider the following $\DR{X}$-valued differentiation on $\Omega_{poly}(G)$ (one can regard it as a vector field on $G$ with coefficients  in $\DR{X}$ given by the global gauge potential $A$):
\begin{equation}
\label{eq15}
F_A=\sum_{i=1}^N A_i\otimes X_i.
\end{equation}
We can consider $\Omega_{DR}(X)$ as a subspace of its bar-resolution. In order to distinguish between the map to $\Omega_{DR}(X)$ and the corresponding map to the bar resolution, we shall write $\widehat F_A$ in the latter case. Then it is clear that $\widehat F_A$ is also a differentiation on $\Omega_{poly}^0(G)$ with values in the bar-resolution, this time with respect to the shuffle products. Just recall, that $1\in\mathbb C=\Omega_{DR}(P)^{\otimes 0}\subseteq B(\Omega_{DR}(P))$ is unit with respect to the shuffle products.

The formula \eqref{eq13} can now be rewritten as
\begin{equation}
\label{eq16}
\widetilde M_A(f)=\sum_n \widehat F_A^{\otimes n}(f)|_e.
\end{equation}
Here, $\widehat F_A^{\otimes n}$ is a $B(\Omega(P))\otimes\ac(G)$-valued differential operator
on $\ac(G)$:
\begin{equation}\label{eq16 1/2}
\widehat F_A^{\otimes n}(f)=\sum_{(\bar k)}[A_{k_1}|\dots|A_{k_n}]\otimes X_{k_1}\dots X_{k_n}(f).
\end{equation}

In effect for any differential graded algebra $\Omega$ and any linear map $\alpha:\g\to\Omega$, one can consider an $\Omega$-valued vector field on $G$ by putting
\begin{equation}
\label{eq17 1/2}
F_\alpha\eqdef\sum_{i=1}^N \alpha(X_i)\otimes X_i.
\end{equation}
Once again, one can consider the corresponding maps to the bar-resolution, which we shall denote by $\widehat F_\alpha$. Now, observe, that (for $\alpha=A$, the connection form)
\begin{equation}\label{eq16 2/3}
\begin{split}
\sum_{k_1,\dots,k_n=1}^N[A_{k_1}|\dots|A_{k_n}]&=[\sum_{k_1=1}^NA_{k_1}|\dots|\sum_{k_n=1}^NA_{k_n}]\\
                                               &=\frac{1}{n!}{\rm sh}(\sum_{k_1=1}^NA_{k_1},\dots,\sum_{k_n=1}^NA_{k_n})\\
                                               &=\frac{1}{n!}\sum_{k_1,\dots,k_n=1}^N{\rm sh}(A_{k_1},\dots,A_{k_n}).
\end{split}
\end{equation}
In the virtue of this formula, we can rewrite the equation \eqref{eq16} as
\begin{equation}\label{eq16 3/4}
\widetilde M_A(f)=\sum_n \frac{1}{n!}\widehat F_A^{\,{\rm sh}(n)}(f)|_e,
\end{equation}
where
\begin{equation}\label{eq16 4/5}
{\rm sh}(\widehat F_{\alpha_1},\dots,\widehat
F_{\alpha_n})(f)\eqdef\sum_{\sigma\in\Sigma_n,(\bar k)}(-1)^{\tilde\sigma}[\alpha_{\sigma^{-1}(1)}(h_{k_1})|\dots|\alpha_{\sigma^{-1}(n)}(h_{k_n})]\otimes
X_{k_1}\dots X_{k_n}.
\end{equation}
Now the equation $\widetilde M_A(fg)={\rm sh}(\widetilde M_A(f),\,\widetilde M_A(g))$ follows from
the Leibnitz rule for $\widehat F_A$ (with respect to the shuffle product).

In addition it is easy to show, that $\widehat F_A|_e$ verifies the Leibnitz rule with
respect to the coalgebra structure (comultiplication and counit):
\begin{equation}\label{eq17}
  (\widehat F_A|_e\otimes\epsilon+\epsilon\otimes\widehat F_A|_e)\circ\Delta=\nabla\circ\widehat F_A|_e.
\end{equation}
We shall use the Sweedler's notation for the coproduct: $\Delta(f)=\sum_{(f)}f_{(1)}\otimes f_{(2)}$ and omit the sign $|_e$. We compute:
\begin{equation*}
  \begin{split}
    (\widehat F_A\otimes\epsilon+\epsilon\otimes\widehat F_A)(\Delta(f))&=(\widehat F_A\otimes\epsilon+\epsilon\otimes\widehat F_A)(\sum_{(f)} f_{(1)}\otimes
f_{(2)})\\
                                                           &=\sum_{(f)}\Bigl\{\widehat F_A(f_{(1)})\otimes\epsilon(f_{(2)})+\epsilon(f_{(1)})\otimes
\widehat F_A(f_{(2)})\Bigr\}\\
                                                           &=\widehat F_A(f)\otimes1+1\otimes
\widehat F_A(f)=\nabla (\widehat F_A(f)).
  \end{split}
\end{equation*}
We use the counit relation $\sum_{(f)}f_{(1)}\epsilon(f_{(2)})=f=\sum_{(f)}\epsilon(f_{(1)})f_{(2)}$. As above, we conclude that $\widetilde M_A$ is a homomorphism of coalgebras too.
\end{Proof}\\ \\ \vspace{1mm}
Now, since the map $\widetilde M_A$ is a homomorphism of commutative bialgebras, we can extend it to the algebra $\Omega_{poly}(G)$ of poilynomial differential forms on $G$. Clearly, the latter map will also be a homomorphism of (this time) graded algebras. In order to investigate the behaviour of this homomorphism with respect to the differentials we shall consider a little bit more general situation.

Namely, let $\alpha$ and $\beta$ be two homogeneous maps (with respect to the degree of image) $\alpha,\beta:\g\to\Omega$ of degrees $|\alpha|,|\beta|>0$ respectively (here $\Omega$ is an arbitrary differential graded algebra as above.) Then in the notation, explained above, we have the following formula (see \eqref{eq16 4/5}):
\begin{equation}\label{eq17 2/3}
\begin{split}
(d_B\otimes1)\bigl({\rm sh}(&\widehat F_\alpha,\,\widehat F_\beta)(f)\bigr)=\\ &=\bigl({\rm sh}(\widehat F_{d\alpha},\,\widehat F_\beta)+(-1)^{|\alpha|-1}{\rm sh}(\widehat F_\alpha,\,\widehat F_{d\beta})-(-1)^{|\alpha|-1}\widehat{[F_\alpha,F_\beta]}\bigr)(f),
\end{split}
\end{equation}
for all $f\in\ac(G)$. Here $[F_\alpha,\,F_\beta]=\sum_{i,k=1}^N\alpha(h_i)\beta(h_k)\otimes[X_i,\,X_k]$ is the commutator of vector fields $F_\alpha$ and $F_\beta$.

In fact, first two terms on the right side of \eqref{eq17 2/3} are obvious. They correspond to the De Rham part ($d_I$, see above) of
the bar-resolution's differential. Let's show, that the multiplication part of this differential gives the third term. We compute
\begin{equation}\label{eq17 3/4}
\begin{split}
(d_{II}&\otimes1)\bigl({\rm sh}(\widehat F_{\alpha},\,\widehat F_\beta)(f)\bigr)\\
                                                 &=(d_{II}\otimes1)\Bigl(\sum_{i,k=1}^N[\alpha(h_i)|\beta(h_k)]\otimes
X_iX_k(f)\Bigr)\\
                                                 &\quad+(-1)^{(|\alpha|-1)(|\beta|-1)}(d_{II}\otimes1)\Bigl(\sum_{i,k=1}^N[\beta(h_i)|\alpha(h_k)]\otimes
X_iX_k(f)\Bigr)\\
                                                 &=\sum_{i,k=1}^N d_{II}\Bigl([\alpha(h_i)|\beta(h_k)]+(-1)^{(|\alpha|-1)(|\beta|-1)}[\beta(h_i)|\alpha(h_k)]\Bigr)\otimes X_iX_k(f)\\
                                                 &=\sum_{i,k=1}^N\bigl((-1)^{|\alpha|-1}\alpha(h_i)\beta(h_k)+(-1)^{|\alpha|}\alpha(h_k)\beta(h_i)\bigr)\otimes X_iX_k(f)\\
                                                 &=(-1)^{|\alpha|-1}\sum_{i,k=1}^N\alpha(h_i)\beta(h_k)\otimes(X_iX_k(f)-X_kX_i(f))=\widehat{[F_\alpha,F_\beta]}(f).
\end{split}
\end{equation}
This formula easily generalizes to the $n$-fold shuffle product of the aforesaid described
fields of strictly positive degrees:
\begin{equation}\label{eq17 4/5}
\begin{split}
(d_B\otimes1)&{\rm sh}(\widehat F_{\alpha_1},\dots,\widehat
F_{\alpha_n})\\
             &=\sum_{k=1}^n(-1)^{\varepsilon_k}{\rm sh}(\widehat F_{\alpha_1},\dots,\widehat F_{d\alpha_k},\dots,\widehat
F_{\alpha_n})\\
             &\quad+\sum_{i<j}(-1)^{\varepsilon_{ij}}{\rm
sh}(\widehat{[F_{\alpha_i},\,F_{\alpha_j}]},\,F_{\alpha_1},\dots,\underbrace{\widehat{}}_i\,,\dots,\underbrace{\widehat{}}_j\,,\dots,F_{\alpha_n}).
\end{split}
\end{equation}
These formulae should be slightly modified, when some of these maps take values in $\Omega^0$. Namely we shall assume that the algebra $\Omega$ is equipped with an augmentation $\Omega\to\mathbb C$. In the case of de Rham algebras of a manifold $X$, this augmentation is given by sending a function $f$ to its vakue in a point $x_0$. Then, if $|\alpha|=0$ one should add the terms $\alpha|_{x_0}(\widehat F_\beta(f))-\widehat F_\beta(\alpha|_{x_0}(f))=[\alpha|_{x_0},\,\widehat F_\beta](f)$ to the right hand side of the formula \eqref{eq17 2/3}. Here we regard $\alpha|_{x_0}$ as a $\mathbb C\subseteq B(\Omega_{DR}(X))$-valued vector field on $G$. Respectively, if $|\alpha_1|=0$, then one should add to the formula \eqref{eq17 4/5}, the following commutator of differential operators: $[\alpha_1|_{x_0},\,{\rm sh}(\widehat F_{\alpha_2},\dots,\widehat F_{\alpha_n})]$.

Besides this, the shuffle product of vector fields, defined by equation \eqref{eq16 4/5} is graded commutative. Indeed let's compute in the case of two vector fields for example:
\begin{equation*}
\begin{split}
{\rm sh}(\widehat F_{\alpha},\,\widehat F_\beta)(f)&=\sum_{i,k=1}^N[\alpha(h_i)|\beta(h_k)]\otimes X_iX_k(f)\\
                               &\quad+(-1)^{(|\alpha|-1)(|\beta|-1)}\sum_{i,k=1}^N[\beta(h_i)|\alpha(h_k)]\otimes
X_iX_k(f)\\
                               &=(-1)^{(|\alpha|-1)(|\beta|-1)}\Bigl(\sum_{i,k=1}^N[\beta(h_i)|\alpha(h_k)]\otimes
X_iX_k(f)\\
                               &\quad+(-1)^{(|\alpha|-1)(|\beta|-1)}\sum_{i,k=1}^N[\alpha(h_i)|\beta(h_k)]\otimes X_iX_k(f)
\Bigr)\\
&=(-1)^{(|\alpha|-1)(|\beta|-1)}{\rm sh}(\widehat F_{\beta},\,\widehat
F_\alpha)(f),
\end{split}
\end{equation*}
for all $f\in\ac(G)$. We shall use this below. However, it is not correct to think of the formula \eqref{eq16 4/5} as defining an associative product of differential operators. For instance, it is not clear, how to define the product of operators of higher degrees.

\subsection{Homological monodromy map and the gauge transformations}
\label{sect3}
Let the group $G$ act on itself by adjoint action and consider the induced action on the corresponding algebra of smooth functions $\Omega_{poly}^0(G)=\ac(G),\ f\mapsto f^g,\ f^g(x)=f(g^{-1}xg)$. Compose this action with $\widetilde M_A$. The purpose of this section is to determine the equivariance properties of the map $\widetilde M_A$ with respect to this action.

Since we assumed that $G$ is connected, it is enough to consider the infenitesimal part of this action. So we identify $\g$ with the space of left-invariant vector fields on $G$ and let it act on $\ac(G)$ by Lie derivatives. We ask about the equivariance of $\widetilde M_A$ under this action.

To answer this question we first recall, that any map $g:X\to G$ induces an endomorphism of the trivial principal bundle $X\times G$ (in effect, the group of endomorphisms of this bundle can be identified with the set of such maps with poinwise multiplication.) In particular, if $X=P$ (for an arbitrary principal bundle $P$) and we work with the pullback bundle $\pi^*(P)=P\times G$, then one can regard the transformations of $P$ as tose transformations of $\pi^*(P)$, which can be pulled down to the bundle $P$. This amounts to considering the subgroup of $Ad_G-$equivariant maps $g:P\to G,\ g(ph)=\inv{h}g(p)h$. Thus the general case can be reduced to the trivial one, which we consider here. (Observe, however, that our construction has so far been given only for trivial bundles. The last remark should be considered in the light of the next sections, see e.g. paragraph \ref{sect17}.)

The map $g$ induces the gauge transformation of the connection form on $X$:
\begin{equation}\label{eq17 5/6}
A\mapsto A^g=g^{-1}Ag+g^{-1}dg.
\end{equation}
We claim that the homological monodromy map $\widetilde M_A$ intertwines the adjoint action of $G$ on itself with the action, induced on the image of $\widetilde M_A$ by the gauge transformations. Here we identify an element $g\in G$ with the gauge transformation, induced by constant map $g:X\to G$.

The statement will follow from the explicit description of the action, induced by the gauge transformations on the image of $\widetilde M_A$ in the normalized bar-resolution (see \eqref{eq28}). As we have said above it is enough to consider the infinitesimal action of gauge transformations. Namely, let the $g=\exp(ta),\ t\to0$, where $a:P\to\g,\ a(ph)=Ad_h(a(p))$, be an infenitesimal gauge transformation. Then one can write down the formula \eqref{eq17 5/6} as
\begin{equation}\label{eq18}
A\mapsto A+ t\delta_a(A)+\dots =A+ t([A,a] + da)+\dots,
\end{equation}
where $+\dots$ denotes the sum of terms of quadratic and higher orders. If $a(p)=\sum_{k=1}^N a_k(p)\,h_k,\ a_i\in C^\infty(P)$, then
\begin{equation}\label{eq19}
 \delta_a(A)_k=\sum_{i,j=1}^N C_{ij}^k A_i a_j +da_k.
\end{equation}
Here $C_{ij}^k$ are the structure constants of the Lie algebra $\g$:
\begin{equation}\label{eq20}
[h_i,h_j]=C_{ij}^kh_k.
\end{equation}

Under this transformation, the vector field $F_A$ transforms into
$$F_{A+t\delta_a(A)+\dots}=F_A+ t(F_{[A,a]}+F_{da})+\dots=F_A + t([F_A,\,F_a]+F_{da})+\dots.$$
Now, using the definition \eqref{eq13} and equation \eqref{eq16 2/3}, we compute:
\begin{equation}\label{eq21}
\begin{split}
&\widetilde M_{A+t\delta_a(A)+\dots}(f)=\sum_{n=0}^\infty \frac{1}{n!} \widehat F^{\,{\rm sh}(n)}_{A+t\delta_a(A)+\dots}(f)|_e\\
                      &\qquad=\sum_{n=0}^\infty \frac{1}{n!}\Bigl(\sum_{k=0}^n\genfrac{(}{)}{0pt}{}{n}{k} t^{n-k}{\rm sh}(\underbrace{\widehat F_A,\dots,\widehat F_A}_{k\ {\rm times}},\widehat F_{\delta_a(A)},\dots,\widehat F_{\delta_a(A)})\Bigr)(f)|_e+\dots\\
                      &\qquad=\widetilde M_A(f)+t\sum_{n=1}^\infty\frac{1}{(n-1)!}{\rm sh}(\widehat F_{\delta_a(A)},\widehat F_A,\dots,\widehat F_A)+\dots,
\end{split}
\end{equation}
where $\dots$ consists of quadratic and higher order terms. Here we used the graded commutativity of the shuffle product of vector fields. Thus, the linear part of the gauge transformation action on the image of $\widetilde M_A$ is equal to
\begin{equation}\label{eq22 1/2}
\widetilde M^a_A=\frac{d}{dt}(\widetilde M_{A^g})|_{t=0}=\sum_{n=0}^\infty\frac{1}{n!}{\rm sh}(\widehat F_{\delta_a(A)},\underbrace{\widehat F_A,\dots,\widehat F_A}_n).
\end{equation}

Let now $\alpha_1,\dots,\alpha_n:\g\to\Omega_{DR}(X)$ be a collection of linear maps. For the sake of simplicity we shall assume, that they all are of odd degrees. Then using \eqref{eq17 4/5} we obtain the following formula:
\begin{equation}\label{eq23}
\begin{split}
&d_B({\rm sh}(\widehat F_a,\widehat F_{\alpha_1},\dots,\widehat
F_{\alpha_n})(f)|_e)\\&={\rm sh}(\widehat F_{da},\widehat F_{\alpha_1},\dots,\widehat
F_{\alpha_n})(f)|_e+\sum_{i=1}^n{\rm sh}(\widehat F_{a},\widehat F_{\alpha_1},\dots,\widehat
F_{d\alpha_i},\dots,\widehat F_{\alpha_n})(f)|_e\\
                    &\quad+\sum_{i=1}^n{\rm sh}(\widehat F_{[\alpha_i,\,a]},\widehat F_{\alpha_1},\dots,
\underbrace{\widehat{}}_i\,,\dots,\widehat F_{\alpha_n})(f)|_e\\
&\quad+\sum_{1\le i<j\le n}{\rm sh}(\widehat F_a,\widehat F_{[\alpha_i,\,\alpha_j]},\widehat F_{\alpha_1},\dots,\widehat
F_{\alpha_n})(f)|_e +[a|_{x_0},\,{\rm sh}(\widehat F_{\alpha_1},\dots,\widehat F_{\alpha_n})](f)|_e.
\end{split}
\end{equation}
Now, taking into consideration, the definition of $\ker'\sigma$, the kernel of the Chen's iterated integral map, we have the
following relation
\begin{equation}\label{eq24}
\begin{split}
{\rm sh}(\widehat F_{da},\widehat F_{\alpha_1},\dots,\widehat
F_{\alpha_n})(f)|_e&+\sum_{i=1}^n{\rm sh}(\widehat F_{[\alpha_i,\,a]},\widehat F_{\alpha_1},\dots,
\underbrace{\widehat{}}_i\,,\dots,\widehat F_{\alpha_n})(f)|_e\\
                   &\equiv[{\rm sh}(\widehat F_{\alpha_1},\dots,\widehat F_{\alpha_n}),\,a|_{x_0}](f)|_e({\rm mod}\/(\ker'\sigma))
\end{split}
\end{equation}
Let's apply this formula to the element $\widetilde M^a_A(f)$. We compute:
\begin{align}
\label{eq25}
{\rm sh}(\widehat F_{da},\underbrace{\widehat F_A,\dots,\widehat
F_A}_n(f)|_e&\equiv-n\,{\rm sh}(\widehat F_{[A,\,a]},\underbrace{\widehat F_A,\dots,\widehat F_A}_{n-1})(f)|_e
                   +[{\rm sh}(\underbrace{\widehat F_A,\dots,\widehat F_A}_n),\,a|_{x_0}](f)|_e,\\
\intertext{(modulo $\ker'\sigma$) and hence}
\begin{split}\label{eq26}
\sum_n\frac{1}{n!}{\rm sh}(\widehat F_{da},\underbrace{\widehat F_A,\dots,\widehat
F_A}_{n})&(f)|_e\equiv-\sum_n\frac{1}{(n-1)!}{\rm sh}(\widehat F_{[A,\,a]},\underbrace{\widehat F_A,\dots,\widehat
F_A}_{n-1})(f)|_e\\
                   &\qquad\qquad\qquad\qquad\quad+[\sum_n\frac{1}{n!}{\rm sh}(\underbrace{\widehat F_A,\dots,\widehat F_A}_n),\,a|_{x_0}](f)|_e.
\end{split}
\end{align}
Finally, we conclude, that modulo $\ker'\sigma$
\begin{equation}\label{eq27}
\widetilde M_A^a(f)\equiv[\widetilde M_A,\,a|_{x_0}](f).
\end{equation}
Exponentiating this equality, we conclude, that the normalized image of $\widetilde
M_A$ verifyies the following equation (recall, that $g:P\to G$ -- gauge transformation):
\begin{equation}\label{eq28}
\widetilde M_{A^g}(f)=\widetilde M_A(f^{g^{-1}(x_0)}),
\end{equation}
where for any $g\in G,\ f\mapsto f^{g}$ is the action of $G$ on \ac, induced by the right adjoint action of $G$ on itself, $f^g(x)=f(g^{-1}x\/g)$. This is also true for the induced action on $\Omega(G)$. We see, that $\widetilde M_{A^g}(\omega)=\widetilde M_{A}(\omega)$, if $\omega$ is an invariant (with respect to the adjoint action of $G$ on itself) differential form on $G$. For a constant map $g$, this formula can be written down as
$$
\widetilde M_A(f^{g^{-1}})=\widetilde M_{A^g}(f),
$$
precisely as we claimed before.

\subsection{Homological monodromy of gauge bundles}
\label{sect17}
In this papragraph we biefly discuss the possible extension of our previous construction to the gauge bundles. We omit most details and proofs here, postponing the discussion to sections 3 and 4.

Let $P$ be an arbitrary principal bundle. Consider the "`gauge bundle"' $P_{Ad}$? associated with $P$ (i.e. the bundle, whose sections induce the transformations of the principal bundle $P$.) Recall that $P_{Ad}=P\times_{Ad} G$ (we shall denote by $\pi_{Ad}$ the corresponding projection.) In particular, there is trivial (identical) transformation of $P$, associated with a ``unit'' section of this bundle. Let us now consider the following differential algebra associated with $P_{Ad}$:
$$
\Omega^*_{vpoly}(P_{Ad})=\DR{P}\otimes_G \Omega^*_{poly}(G)=\{\omega\otimes_\varphi\in \DR{P}\otimes \Omega^*_{poly}(G)|\omega^g\otimes\varphi=\omega\otimes\varphi^{g^{-1}},\ g\in G\}.
$$
Here we let $G$ act on $\Omega^*_{poly}(G)$ via the adjoint action. One should regard it as the algebra of differential forms on $P_{Ad}$ which are polynomial in the direction of the fibre. There is a map $1^*:\Omega^*_{vpoly}(P_{Ad})\to\DR{X}$, induced by the counit of $\Omega^*_{poly}(G)$ (equivalently by the inclusion of the unit in $G$.) Now let $A=\sum_k A_k\otimes h_k$ be the (global) connection form on the trivial bundle $P\times G$, induced from the connection on $P$. Take Lie derivatives of the elements in $\Omega^*_{poly}(G)$ with respect to $h_k\in\g$. In this way we associate to $A$ a vector field on $G$ with values in $\DR{P}$. We can in an evident way extand it to $P\times G$ and take the corresponding restriction to the subspace of invariant forms $\Omega^*_{vpoly}(P_{Ad})=\DR{P}\otimes_G\Omega^*_{poly}(G)$. Moreover, the equivariance properties of $A$ show that the image of the induced map still belongs to $\Omega^*_{vpoly}(P_{Ad})$. Taking the iterations and then evaluating the result at the unit section, gives us a map $\widetilde T_A:\Omega^0_{vpoly}(P_{Ad})\to NH(\DR{X})$
\begin{equation}
\widetilde T_A(a\otimes b)=\sum_{n=0}^\infty (\underbrace{Id\otimes\dots\otimes Id}_n\otimes 1^*)A^{\otimes n}(a\otimes b),\ a\in\Omega_{DR}^0(P),\ b\in\Omega_{poly}^0(G).
\end{equation}
Here $NH(\DR{X})$ denotes the \textit{normalized (or reduced) Hochschild complex} of $\DR{X}$. It is defined in the same manner as $N(\Omega_{DR}(X))$, i.e. all the elements of degree $0$ in $\DR{X}$ are factored out. The fact that the image of this map is indeed in $NH(\DR{X})$, follows from the next observations: one can identify $\DR{X}$ with the subalgebra of basic forms in $\DR{P}$. Then the invariance property of the differentiation, induced by $A$:
$$
A^g(b)=A(b^{g^{-1}})=A(b)\ \mbox{for all}\ g\in G
$$
(the last equality follows from the definitions of $1^*$ and the adjoint action) shows that the image of $T_A$ consists of $G$\/-invariant elements. And the second condition, that define the basic elements, follows from the fact that we consider reduced Hochschild complex $NH(\Omega_{DR}(X))$. This complex is a model for the algebra of differential de Rham forms on the free loop space of $X$, when $X$ is 1-connected, see for example \cite{GJP, Chen, Chen2}. A reasoning, similar to the one we used in paragraphs \ref{sect11} and \ref{sect2}, starting with an explicit formula of the time-ordered exponent, shows that the map $\widetilde T_A$ is in fact a homological counterpart of the holonomy map $T_A$, that associates a gauge transformation of $P$ to a connection $A$ and a section of the infinite-dimensional bundle $\mathcal LX\to X$. One can now repeat all the reasoning from previous sections to show that $\widetilde T_A$ can be extended to a homomorphism of differential graded algebras (we shall denote it by the same symbol) and that it intertwines coproduct structures on both sides (on the left it is induced from the product $P_{Ad}\times_X P_{Ad}\to P_{Ad}$ and on the right it is a generalization of the bar-construction coproduct). We shall not discuss this subject now, postponing it to sections 3 and 4, where this will be explained in full detail. In particular one can show, that if $g:X\to G$ is a gauge transformation, then
\begin{equation}
\label{thochmap}
T_{A^g}(\alpha)=T_A(\alpha^{g^{-1}}),\ \text{for any}\ \alpha\in\Omega_{vpoly}^*(P_{Ad}).
\end{equation}
Here as before the action of $g$ on the algebra $\Omega_{vpoly}^*(P_{Ad})$ is induced from its action on functions, which is given by the product of the agrument by $g$ from the left (recall, that the bundle $P_{Ad}$ is in fact a group fibration, in particular, one can multiply its elements on its sections on both sides).

This construction is quite convenient, but not very explicit, since it is not very easy to see explicitly why the image of $\widetilde T_A$ belongs indeed to $NH(\DR{X})$. One can try to amend it by considering trivial bundles $P=X\times G$. Observe, that in this case $P_{Ad}$ is again trivial, $P_{Ad}\cong X\times G$, but it shall not be regarded as a principal bundle any more, (i.e. $G$ on the right hand side of this homeomorphism plays rather the role of non-linear representation of the group $G$.) It is very tempting then to consider the following construction. Let $P$ be \emph{an arbitrary} principal bundle, $P_{Ad}$, defined as above -- the associated gauge bundle. Consider a trivializing open cover of $\mathcal{U}=\{U_i\}$ of $P$ (here $U_i\subseteq X$ are open subsets, such that $\pi^{-1}(U_i)\cong U_i\times G$, and $\cup_iU_i=X$). We can assume that this cover is finite. Let $\mathcal V=\{V_i\}$, where $V_i=\pi_{Ad}^{-1}(U_i)$, be the corresponding open cover of $P_{Ad}$. One then can choose the partition of unity $\varphi_i$, associated with $\mathcal V$.

Now take an arbitrary differential form $\alpha\in\Omega_{DR}(P_{Ad})$, decompose it into a sum of differential forms on $V_i$s, $\alpha=\sum\alpha_i,\ \alpha_i=\varphi_i\alpha$ and apply the previos construction to every $\alpha_i$, which is possible since $V_i={P_{Ad}}_{|_{U_i}}\cong U_i\times G$ is a trivial bundle over $U_i$.

Applying the map $T_{A_{U_i}}$ to every $\alpha_i$ (where $A_{U_i}$ is the local gauge potential in $U_i$ of a globally defined connection on $P$), we obtain a collection of elements $T_{A_{U_i}}(\alpha_i)\in NH(\Omega_{DR}(U_i))$. After this we can regard them as elements of $NH(\Omega_{DR}(X))$ (using another partition of unity, if necessary). Then, since the glueing transformations of the fibration $P_{Ad}$ are given by the conjugations of $G$ by the cocycle $g_{ij}:U_i\cap U_j\to G$ of the principal bundle $P$, we conclude (using formula \eqref{thochmap}) that these elements are compatible on the intersections of $U_i$s and hence define an element of $NH(\Omega_{DR}(X)$.

More accurately, we can consider the \v{C}ech complex of the open cover $\uc$ with values in the presheaf of the normalized Hochschild complexes of the differential forms, $\check{C}(\mathcal U,\,NH(\Omega_{DR}))$. Then the elements $T_{A_{U_i}}(\alpha_i)$ determine a \v Cech cocycle in $\check{C}(\mathcal U,\,NH(\Omega_{DR}))$ and hence we obtain a chain map from $\Omega_{DR}(P_{Ad})$ to $\check{C}(\mathcal U,\,NH(\Omega_{DR}))$. The latter complex is regarded with the total differential equal to the sum of all three natural differentials: de Rham differential on the forms, Hochschild differential and \v Cech differential.

However, these constructions are hardly of any use if one wants to obtain an element in the cohomology of the free loop space. In effect, suppose that the cover \uc\ is fine enough (in particular all the intersections of its elements should be contractible). Consider then the filtration of $\check{C}(\mathcal U,\,NH(\Omega_{DR}))$ induced by the degrees of the elements in $NH(\Omega_{DR}(U))$ and use the corresponding spectral sequence to calculate its cohomology. It turns out that the $E^1$ term of this sequence is equal to the \v Cech complex of \uc\ with values in the sheaf of the Hochschild cohomology of $U_i$s. But since we assumed that all the intersections were contractible, this sheaf is constant and concentrated in degree 0, thus the $E^2$ term of this sequence is equal to the cohomology of the nerve of \uc, which is known to be equal to $H^*(X)$ and the spectral sequence collapses after the second term.

Thus the complex $\check{C}(\mathcal U,\,NH(\Omega_{DR}))$ is not equivalent to $NH(\Omega_{DR}(X))$ and one has to find another construction to incorporate the nontrivial bundles in this picture. This will be done in the next sections.


\section{Twisting maps and flat connections}
In this section we discuss a more abstract approach to the holonomy, based on the idea of the twisting cochains and their suitable generalizations. It will lead us in the next section to the abstract characterization of the Getzler-Jones characteristic class (also known as \emph{cyclic Chern character}) in the terms of the homological monodromy map (Theorem \ref{omagn} and section \ref{sectnnext}).

\subsection{Twisting maps. Examples}
\label{sectnext}
Let $K$ be an arbitrary cochain complex and let $\Omega$ be a differential graded algebra. The following definition is a generalization of the well-known notion of \emph{twisting cochain} (see definition \ref{twcochdf} below).
\begin{df}
\label{twmapdf}
We shall call a linear function
$$
\psi: K^*\to (K\otimes\Omega)^{*+1}
$$
{\em an $\Omega$\/-valued twisting map on $K$\/}, if the map $d_\psi:K\otimes\Omega\to K\otimes\Omega$, which we shall call \emph{twisted differntial, associated with $\psi$}, given by the formula
\begin{equation}
\label{twmap}
d_\psi(\alpha\otimes\beta)=d_K(\alpha)\otimes\beta+(-1)^{|\alpha|}\alpha\otimes
d_\Omega\beta+(1\otimes m)(\psi\otimes1)(\alpha\otimes\beta),
\end{equation}
where $\alpha\in K,\ \beta\in\Omega$ and $m$ is the multiplication in $\Omega$, verifies the equation $d_\psi^2=0$.
\end{df}
Observe, that one can dualize the definition \ref{twmapdf} as follows: let $K$ be a coalgebra with coproduct $\Delta$, $\Omega$ an arbitrary cochain complex, and consider maps $\varpi:K\otimes\Omega\to\Omega$ instead of maps $K\to K\otimes\Omega$. Then we demand that the equation $d_\varpi^2=0$ holds for the map $d_\varpi$, defined by
$$
d_\varpi(\alpha\otimes\beta)=d_K(\alpha)\otimes\beta+(-1)^{|\alpha|}\alpha\otimes
d_\Omega\beta+(1\otimes\varpi)(1\otimes\Delta)(\alpha\otimes\beta).
$$
We shall use this dual construction below. One can also consider the most general type of such maps ($\varsigma:K\otimes\Omega\to K\otimes\Omega$ such that the corresponding twisted differential $d_\varsigma$ verifies the equation $d_\varsigma^2=0$) but we shall not use this most general construction here, because we shall need the bar-resolution of $\Omega$ and/or cobar-resolution of $K$.

If we write down the map $\psi(x)$ as $\psi^{(1)}(x)\otimes\psi^{(2)}(x),\ \psi^{(1)}(x)\in K,\ \psi^{(2)}(x)\in\Omega$, then
$$
d_\psi(\alpha\otimes\beta)=d_K(\alpha)\otimes\beta+(-1)^{|\alpha|}\alpha\otimes
d_\Omega\beta+\psi^{(1)}(\alpha)\otimes\psi^{(2)}(\alpha)\beta,
$$
and the condition $d_\psi^2=0$ is equivalent to
\begin{equation}
\label{twmap1}
\begin{split}
&\psi^{(1)}(d_K\alpha)\otimes\psi^{(2)}(d_K\alpha)=\\
                            &=d_K\psi^{(1)}(\alpha)\otimes\psi^{(2)}(\alpha)+(-1)^{|\psi^{(1)}(\alpha)|}\left(\psi^{(1)}(\alpha)\otimes
d_\omega\psi^{(2)}(\alpha)+\psi^{(1)}(\psi^{(1)}(\alpha))\otimes\psi^{(2)}(\psi^{(1)}(\alpha))\psi^{(2)}(\alpha)\right)
\end{split}
\end{equation}
for all $\alpha\in K$. In what follows we shall usually omit the subscripts $K$ and $\Omega$ in the differentials, if it doesn't cause confusion.

Recall now the definition of \emph{twisting cochains}. Observe that according to our general assumptions, all differentials increase degrees by $1$, and hence we should use cohomological version of the standard definition, which amounts to inverting the direction of all maps (c.f. \cite{Br59}):
\begin{df}
 \label{twcochdf}
A twisting cochain on a coaugmented coalgebra $K$ over field $Bbbk$ with values in a differential algebra $\Omega$ is a linear map $\phi:\Omega^\ast_G\to\Omega^{*+1}$, verifying the condition $\phi\circ\eta=0$ ($\eta:\Bbbk\to K$ is the coaugmentation) and the following relation
\begin{equation}
\label{flatc}
d_\Omega\phi-\phi d_G-\phi\ast\phi=0,
\end{equation}
where $\phi\ast\phi$ is the linear map $\Omega^\ast_G\to\Omega^{*+2}$, given by $\phi\ast\phi(\omega)=\phi(\omega^{(1)})\phi(\omega^{(2)})$ (we use the standard notation for the comultiplication $\omega\mapsto\omega^{(1)}\otimes\omega^{(2)}$).
\end{df}
Recall that a coaugmentation in differential coalgebras is defined as a coalgebra homomorphism $\eta$ from the 1-dimensional coalgebra $\Bbbk$ (with zero differential) to $K$. Equivalently one can say, that coaugmentation amounts to chosing a group-like element $e,\ de=0,\ \eta(1)=e$ in $K$. For instance below we shall consider $K=\Omega_{poly}^*(G)$, then we shall always assume $e$ to be equal to the constant function $e(g)=1$ on $G$.

Let now $K$ be a differential graded coalgebra and $\tilde\psi$ a linear map on $K$ with values in a differential algebra $\Omega$. It follows from the definition we gave that the map
$$
\psi(k)=\sum (-1)^{|k^{(1)}|}k^{(1)}\otimes\psi(k^{(2)})
$$
is twisting map, iff $\tilde\psi$ is a twisting cochain.
Given a twisting cochain $\phi$ one introduces the \emph{twisted differential} $d_\phi$ on the tensor product $K\otimes\Omega$ as follows:
\begin{equation}
\label{eqtwdiff}
d_\phi=d_K\otimes1+1\otimes d_\Omega-\phi\cap1,\
\phi\cap1(\omega\otimes\psi)=\omega^{(1)}\otimes\phi(\omega^{(2)})\psi.
\end{equation}
This is precisely the twisted differential, associated with the twisting map, constructed from $\phi$ as above. We shall call the tensor product, equiped with this differential, \emph{twisted tensor product} and denote it by $K\otimes_\phi\Omega$. Similar notation will be used for the tensor product, equiped with a twisted differential $d_\psi$ for a twisting map $\psi$.

The following example is one of our principal sources of inspiration in this paper.
\begin{ex}\rm
\label{example1}
Let $\g$ be a (differential graded) Lie algebra acting on a chain complex $K$. For example, one can suppose that $K$ is a differential graded coalgebra, then a natural choice for $\g$ is the dg Lie algebra of codifferentiations in $K$, i.e. of linear maps $X:K\to K$ (changing the degree of elements), such that
$$
(Xk)^{(1)}\otimes (Xk)^{(2)}=X(k^{(1)})\otimes k^{(2)}+(-1)^{|X||k^{(1)}|}k^{(1)}\otimes X(k^{(2)}),
$$
for all $X\in\g,\ k\in K$. Similarly, one can take $K$ to be an algebra and consider the Lie algebra of its graded differentiations. Another important example procures the standard Cartan-Chevalley cohomology complex of the Lie algebra \g. It is well-known, that unless \g\ is commutative, this complex is not a differential coalgebra. However, \g\ always couples with its own complex in two different ways: by Lie derivatives and by the inner products with elements of \g. One can use these couplings and the construction given in this example below to obtain twisting maps with the Chevalley complex as their domain.

So choose a linear basis $X_i,\ i=0,\dots,n$ in $\g$. One can now define a map $K^*\to(K\otimes\Omega)^{*+1}$ by the formula
$$
\psi(k)=\sum_i (-1)^{\varepsilon_i}X_i(k)\otimes a^i,
$$
where $a^i$ is a set of elements in $\Omega,\ deg\,a^i=-deg\,X_i+1$ and $\varepsilon_i=|k||a^i|$.

In this case the condition \eqref{twmap1} takes the form
\begin{equation}
\label{eq5102}
\sum_i \{dX_i\otimes a^i+(-1)^{|X_i|}X_i\otimes d a^i\}=\sum_{i,j}(-1)^{\varepsilon_{ij}}X_jX_i\otimes a^j a^i,
\end{equation}
where $\varepsilon_{ij}=(|a^i|+1)|a^j|+1$.

Assume that $\Omega$ is a graded commutative algebra, i.e. $ a^j a^i=(-1)^{|a^j||a^i|}a^i a^j$. Then we can rewrite the equation \eqref{eq5102} as follows
$$
\sum_i d(X_i\otimes a^i)-\sum_{i<j}[X_i,\,X_j]\otimes a^i a^j=0,
$$
or, for short $dA+\frac12[A,\,A]=0$, where $A$ is the formal expression $\sum_i X_i\otimes a^i$. Thus we see that twisting map in this particular case is a variant of the flat connection on a principal bundle.
\end{ex}

\vspace{1mm}
The general purpose of this and next subsections is to outline a way in which such slightly more general twisting maps (flat connections) can be used as a substitute for twisting cochains. But before we go in this direction, let us give a natural construction, leading to such twisting map.

Let $P$ be a trivial (but not trivialized) principal bundle with base $X$ and structure group $G$. Let $\omega$ be an arbitrary connection on $P$ (non flat, in general). By abuse of language we shall not distinguish between a connection on $P$ and the corresponding connection $1$\/-form. Since $P$ is trivial, there exists a global section $s:X\to P$. We can use this section to pull back $\omega$ from $P$ to $X$. Thus we obtain a $\g$\/-valued 1-form $A=\sum A^i\otimes X_i$ on $X$ (where $\{X_i\}$ is a base of $\g$ and $A^i$ are 1-forms on $X$). It is clear, that if we exchange the order of $X_i$ and $A^i$ and interprete $X_i$ as left-invariant vector fields on $G$ (and hence as degree $0$ differentiations of the de Rham algebra of $G$, we shall obtain a map of the sort, described in example \ref{example1}. However, this is not necessarilly a twisting map, since the connection may fail to be flat.

So we modify it a little bit. Let $F=dA-\frac12[A,\,A]$ be the curvature form of the connection $A$. We can write $F$ down as a sum
$$
F=\sum_iF^i\otimes X_i,\ \quad F^i=dA^i-\sum_{i<j}C_{jk}^iA^j\wedge A^k,
$$
where $C^i_{jk}$ are the structure constants of the Lie algebra $\g$. Put now
$$
\tilde A=\sum\{A^i\otimes X_i-F^i\otimes I_i\},
$$
where $I_i$ is the internal product by $X_i$ (contraction of differential forms with vector field $X_i$) in the de Rham algebra of $G$. It is a degree $-1$ differentiation of $\Omega^*_{DR}(G)$, thus we obtain a map of the given type again. Now, with the help of Cartan's formulas $X_i=dI_i$, $[X_i,\,I_j]=C_{ij}^kI_k$ and $[I_i,\,I_j]=0$ and Bianchi's identity $dF^i=\sum_{j,k}C^i_{jk}A^j\wedge F^k$ we compute
\begin{equation*}
\begin{split}
d\tilde A-\frac12[\tilde A,\,\tilde A]&=\sum_i\{dA^i\otimes X_i-dF^i\otimes I_i-F^i\otimes dI_i\}\\
                                              &\quad-\sum_{j<k}\{A^j\wedge A^k\otimes [X_j,\,X_k]+F^j\wedge F^k\otimes [I_j,\,I_k]\}\\
                                              &\quad+\sum_{j,k}\{A^j\wedge F^k\otimes X_jI_k-F^j\wedge A^k\otimes I_jX_k\}\\
                                              &=\sum_i\{dA^i\otimes X_i-dF^i\otimes I_i-F^i\otimes X_i\}\\
                                              &\quad-\sum_{j<k}C^i_{jk}A^j\wedge A^k\otimes X_i+\sum_{j,k}C^i_{jk}A^j\wedge F^k\otimes I_i=0
\end{split}
\end{equation*}
In a similar way, if the group $G$ acts smoothly on a manifold $F$, we obtain a twisting map on the algebra of differential de Rham forms on $F$: one should use the differentiations of this algebra by the vector fields induced by 1-parameter groups $\exp(tX),\ X\in\g$.

Unfortunately, the twisting map, given by this construction is trivial in some sense (we shall discuss this below), which conforms to the fact that the bundle under consideration is trivial. One can amend this default by considering this construction only on the local level (i.e. choosing the trivializing cover of the bundle, etc.) and then glueing the resluting maps in a global map with the help of the cocycle $g_{\alpha\beta}$, as we shall show in the next chapter. Another approach is to consider the gauge bundle and construct a variant of twisting map for it, as we have discussed it in paragraph \ref{sect17} (see also section \ref{sect43}.)

\begin{ex}\rm
\label{expoly}
One can try to extand the theory of (flat) connections and associated twisting maps to the case when $\Omega$ is not a commutative algebra. This leads to a kind of connections with values in universal enveloping algebra.

Namely let $\g$ be a differential graded Lie algebra, acting on the complex $K$. Let $U\g$ be the universal enveloping algebra of $\g$. It is a (noncommutative) differential graded Hopf algebra, acting on $K$. Choose a linear basis $X_I$ of $U\g$ (e.g. we can assume that a Poincare-Birhoff-Witt type theorem holds for $U\g$, then one can choose $X_I$ to be monomials in variables $X_i$.) Consider the following map
\begin{equation}
\label{polyflat}
\psi(k)=\sum_I X_I(k)\otimes A^I.
\end{equation}
Although this sum is infinite in a general case, we shall not investigate the problem of convergence in here. Once again, the condition \eqref{twmap1} reduces to
$$
\sum_I\{dA^I\otimes X_I+(-1)^{|A^I|}A^I\otimes dX_I\}=\sum_{I,J}A^IA^J\otimes X_IX_J.
$$
As in the previous example, we shall give a method to obtain such maps in a particular case, where $\g$ is the Lie algebra of a group $G$ or its generalization, and $K$ is the de Rham algebra of $G$ (or any other complex, on which $\g$ acts, such as the Cartan-Eilenberg complex). In the case we extand $\g$ to a differential graded Lie algebra of Cartan operators $\tilde\g=\g_{-1}\oplus\g_0,\ \g_0=\langle X_1,\dots,X_n\rangle,\ \g_{-1}=\langle I_1,\dots,I_n\rangle,\ dI_i=X_i$. By abuse of notation we shall write $U\g$ instead of $U\tilde\g$.

One can ask, whether there exist such generalized twisting maps in this context at all. The following lemma is intended to provide a wide class of such maps in a particular case when $\g$ is a nilpotent Lie algebra (the nilpotency condition is necessary to circumvene the convergence problems.)
\begin{lemma}
\label{polyex}
Let $X_i$ be the basis of a nilpotent Lie algebra $\g$ and $A^i$ -- a collection of 1-forms (elements of degree 1) in $\Omega$ ($\Omega$ is a non-commutative differential graded algebra). Then there exists a twisting map $\psi$ of type \eqref{polyflat}, beginning with $\psi(k)=X_1(k)\otimes A^1+\dots+X_n(k)\otimes A^n+\dots$. Moreover, the coefficients of $\psi$ can be expressed in the terms of the 1-forms $A^i$.
\end{lemma}
\begin{Proof}
Let us remind, that in this section we do not take care of the convergeance questions (we can do it here since $\g$ is assumed to be nilpotent). So, we can base our proof on the induction with respect to the natural filtration in $U\g$, $F^nU\g=p(\bigoplus_{k\ge n}\g^{\otimes k}),\ n\ge1$, where $p:T\g\to U\g$ is the natural projection. More accurately, we form the Rees' ring of the filtered ring $U\g$. It is a graded ring
\begin{equation*}
rU\g:=\mathbb C\oplus\bigoplus_{i}t^iF^iU\g,
\end{equation*}
where $t$ is a formal parameter. The grading is given by the degree of $t$. One can put $deg\,t=1$ or $2$ (the latter is preferrable, since it is usually presumed that $t^2=0$, if $deg\,t=1$). There's a family of formal homomorphisms $p_t$ from $rU\g$ to $U\g$, given by prescribing some complex value to $t$ (formality is due to the appearance of infinite sums in the formulae). If $\g$ is nilpotent, these maps are actually homomorphisms, so we can limit ourselves with producing such formal maps. Also observe, that the differential $d$ on $U\g$, which comes from $\g$, survives on $rU\g$, since it doesn't change the filtration of an element.

The action of $U\g$ on $K$ generates the $\mathbb C[t]$\/-linear action of $rU(\g)$ on $K[t]$. Thus instead of the twisting maps $\psi:K\to K\otimes\Omega$ we can consider the $\mathbb C[t]$\/-linear twisting maps $\psi_t:K[t]\to(K\otimes\Omega)[t]$. In effect, one can try to consider the non-$\mathbb C[t]$\/-linear maps as well, extending $\mathbb C[t]$ to $\Omega(\mathbb C^1)$ and changing $U\g$ suitably. This seems to be quite an intersting and non-vacuous theory. But we shall not need it here.

Thus, our aim is to show that there exists a series
\begin{equation}
\label{series}
\psi_t=\sum_n t^n\bigl(\sum_{X_I\in F^nU\g}X_I\otimes A^I\bigr),
\end{equation}
which satisfies the equation
\begin{equation}
\begin{split}
\sum_n t^n\bigl(\sum_{X_I\in F^nU\g}&\{dA^I\otimes X_I+(-1)^{|A^I|}A^I\otimes
dX_I\}\bigr)\\
&=\sum_n t^n\bigl(\sum_{\substack{X_I\in F^pU\g,\ X_J\in F^qU\g\\ p+q=n}}A^IA^J\otimes
X_IX_J\bigr),
\end{split}
\end{equation}
and begins with $tX_1\otimes A^1+\dots+tX_n\otimes A^n+\dots$.

Now we can use the induction in the degree of the parameter $t$. That is, we shall assume that there exist the $n$\/th degree part of \eqref{series}: $\psi_t^n=\sum_{k=1}^nt^n\sum_{I\in F^kU\g}X_I\otimes A^I$, such that
\begin{equation}
\label{stepn}
d\psi_t^n-\psi_t^n\circ\psi_t^n=O(t^{n+1}).
\end{equation}
We shall then show that there exist elements $A^J\in\Omega,\ X_J\in F^{n+1}U\g$, such that $\psi^{n+1}_t=\psi_t^n+t^{n+1}\sum_{X_J\in F^{n+1}U\g}X_J\otimes A^J$ verifies the equation \eqref{stepn} with $n+1$ instead of $n$.

It is enough to do this in the case, when the algebra $\Omega$ is free differential graded algebra $\Omega(n)$, generated by $n$ degree $1$ elements $e^1,\dots,e^n$. Then sending $e^k$ to $A^k$ we obtain the formula in the case of a arbitrary $\Omega$.

The first step of induction is evident: it follows from the formula $dI_k=X_k$, that $\psi^1_t=t(X_1\otimes e^1+\dots+X_n\otimes e^n-I_1\otimes de^1-\dots-I_n\otimes de^n)$ verifies the equation \eqref{stepn} for $n=1$. The following lemma is important for the induction step. Its proof is a simple consequence of the definitions.
\begin{lemma}
Differential algebra $\Omega(n)$ is acyclic. If differential graded Lie algebra $\g$ is contractible so that the contracting homotopy is a (graded) derivative of \g (in particular, if $\g$ is composed of elements $X_i$ and $I_i$ as above) then so is $rU\g$.
\end{lemma}
It follows from the Lemma, that in the case we consider, the algebra $rU\g\otimes\Omega(n)$ is acyclic.

Now, for general $n$ let us write $\psi_t^n$ as $\sum_{k=1}^n t^k\psi_k,\ deg\,\psi_k=1$, then comparing the elements of the same power of $t$ in \eqref{stepn}, we obtain the following set of equations:
\begin{align*}
d\psi_1&=0,\\
d\psi_2&=\psi_1\circ\psi_2,\\
d\psi_3&=\psi_1\circ\psi_2+\psi_2\circ\psi_1,\\
&\dots\\
d\psi_n&=\sum_{p+q=n}\psi_p\circ\psi_q,
\end{align*}
and the coefficient at the $n+1$\/st degree of $t$ is equal to $\sum_{p+q=n+1}\psi_p\circ\psi_q$. Its degree is equal to $2$. Now we compute
\begin{equation*}
\begin{split}
d(\sum_{p+q=n+1}\psi_p\circ\psi_q)&=\sum_{p+q=n+1}d\psi_p\circ\psi_q-\sum_{p+q=n+1}\psi_p\circ d\psi_q\\
                                  &=\sum_{p+q=n+1}\bigl(\sum_{r+s=p}\psi_r\circ\psi_s\bigr)\circ\psi_q\\
                                  &\quad-\sum_{p+q=n+1}\psi_p\circ\bigl(\sum_{r+s=q}\psi_r\circ\psi_s\bigr)=0.
\end{split}
\end{equation*}
Since $rU\g\otimes\Omega(n)$ is acyclic, we can choose an element $\psi_{n+1}\in rU\g\otimes\Omega(n)$ of degree $1$, such that $d\psi_{n+1}=\sum_{p+q=n+1}\psi_p\circ\psi_q$. Now we can put $\psi_t^{n+1}=\psi_t^n+t^{n+1}\psi_{n+1}$.
\end{Proof}\\ \\ \vspace{1mm}
However the twisting maps, defined in this lemma are not very useful, since they are all equivalent to trivial ones, just like the map of the previous example. This shall be proven in the a following subsection.
\end{ex}

\subsection{Characteristic map of a twisting map}
\label{sectnext2}
Now we are going to explain how the twisting maps, defined in previous section, can be used to define maps from $K$ to the bar resolution of the algebra $\Omega$. This construction is a slight generalization of the usual map one can associate to a twisting cochain, see for example \cite{SmSOp}.

We suppose now that $K$ is an augmented complex, i.e. that there exists a morphism $\epsilon:K\to\Bbbk$ (the field $\Bbbk$ can in the most part of this and previous section be replaced by an arbitrary ground ring) commuting with differentials. That means in particular that $\epsilon(K^{\ge1})=0$. Then any twisting map $\psi:K\to K\otimes\Omega$ determines a map $\psi_{(1)}:K\to\Omega$, given by the formula
$$
\psi_{(1)}(\alpha)=(\epsilon\otimes 1)\psi(\alpha).
$$
This map changes the dimension of elements by $1$, since $\epsilon(\alpha)=0$ when $\mathrm{dim}\,\alpha>0$. The following statement is a straightforward consequence of the relation \eqref{twmap1}:
\begin{equation}
\label{psi1}
\psi_{(1)}(d\alpha)-d\psi_{(1)}(\alpha)=\psi_{(1)}(\psi^{(1)}(\alpha))\psi^{(2)}(\alpha).
\end{equation}
Let us now iterate this construction: consider the maps $\psi_{(k)}:K\to\Omega^{\otimes k}$, defined by the formula
\begin{equation}
\psi_{(k)}(\alpha)=(\epsilon\otimes 1\dots\otimes 1)(\psi\otimes 1\dots\otimes 1)\dots(\psi\otimes 1)\psi(\alpha),
\end{equation}
where in the leftmost bracket there are $k$ tensors $1$ ($1$ denotes the identity map here.)

\begin{lemma}
\label{lembar}
Assume, that for every $alpha\in K$ there exists a natural number $n=n(\alpha)$, such that $\psi_{(k)}(\alpha)=0$ for all $k\ge n$. Then the map
$$
\tilde\psi(\alpha)=\sum_{k=1}^\infty \psi_{(k)}(\alpha)
$$
is well-defined and determines a homomorphism of (co)chain complexes (i.e. commutes with the differentials).
\end{lemma}
\begin{df}
The map $\tilde\psi$ will be called the \emph{characteristic map of $\psi$}.
\end{df}
\begin{Proof}
The well-definedness of $\tilde\psi$ is a direct consequence of conditions of the lemma. So, we must prove only the fact that this map commutes with the differentials.

First of all in the view of the notation introduced above we can rewrite the equation \eqref{psi1} as
$$
\psi_{(1)}(d\alpha)=d\psi_{(1)}(\alpha)+b\psi_{(2)}(\alpha),
$$
where $b$ stands for the standard differential in the bar-resolution of an algebra. Thus, the map $\tilde\psi$ verifies the equation at the first level of tensor product.

Further let us consider the equation \eqref{twmap1}. Let us apply the map $\psi$ to the left leg of the tensor product on both sides of the equality. We obtain:
\begin{equation*}
\begin{split}
(\psi\otimes 1)\psi(d\alpha)&=\psi^{(1)}(\psi^{(1)}(d\alpha))\otimes\psi^{(2)}(\psi^{(1)}(d\alpha))\otimes\psi^{(2)}(d\alpha)\\
                                                 &=\psi^{(1)}(d\psi^{(1)}(\alpha))\otimes\psi^{(2)}(d\psi^{(1)}(\alpha))\otimes\psi^{(2)}(\alpha)\\
                                                 &\quad+(-1)^{|\psi^{(1)}(\alpha)|}\bigl(\psi^{(1)}(\psi^{(1)}(\alpha))\otimes\psi^{(2)}(\psi^{(1)}(\alpha))\otimes d\psi^{(2)}(\alpha)\\
                                                 &\quad+\psi^{(1)}(\psi^{(1)}(\psi^{(1)}(\alpha)))\otimes\psi^{(2)}(\psi^{(1)}(\psi^{(1)}(\alpha)))\otimes\psi^{(2)}(\psi^{(1)}(\alpha))\psi^{(2)}(\alpha)\bigr)\\
                                                  &=d\psi^{(1)}(\psi^{(1)}(\alpha))\otimes\psi^{(2)}(\psi^{(1)}(\alpha))\otimes\psi^{(2)}(\alpha)\\
                                                  &\quad+(-1)^{|\psi^{1}(\psi^{(1)}(\alpha))|}\bigl(\psi^{(1)}(\psi^{(1)}(\alpha))\otimes d\psi^{(2)}(\psi^{1}(\alpha))\otimes\psi^{(2)}(\alpha)\\
                                                  &\quad+\psi^{(1)}(\psi^{(1)}(\psi^{(1)}(\alpha)))\otimes\psi^{(2)}(\psi^{(1)}(\psi^{(1)}(\alpha)))\psi^{(2)}(\psi^{(1)}(\alpha))1mm\otimes\psi^{(2)}(\alpha)\bigr)\\
                                                  &\quad+(-1)^{|\psi^{(1)}(\alpha)|}\bigl(\psi^{(1)}(\psi^{(1)}(\alpha))\otimes\psi^{(2)}(\psi^{(1)}(\alpha))\otimes d\psi^{(2)}(\alpha)\\
                                                  &\quad+\psi^{(1)}(\psi^{(1)}(\psi^{(1)}(\alpha)))\otimes\psi^{(2)}(\psi^{(1)}(\psi^{(1)}(\alpha)))\otimes\psi^{(2)}(\psi^{(1)}(\alpha))\psi^{(2)}(\alpha)\bigr)
\end{split}
\end{equation*}
If we now apply $\epsilon$ to the leftmost leg of the tensor product on both sides of this equation, we obtain the equality:
$$
\psi_{(2)}(d\alpha)=d\psi_{(2)}(\alpha)+b\psi_{(3)}(\alpha),
$$
i.e. our equality holds at the second degree of tensor algebra too. Repeating this reasoning by induction we obtain the conclusion of this lemma.
\end{Proof}\\

\vspace{1mm}
If one wants to extend this map to the case when the conditions of the proposition is not satisfied, one should consider the completion $\widehat B(\Omega)$ of $B(\Omega)$ with respect to the natural filtration (by tensor powers.) That is, we ought to allow infinite sums of tensor products of elements of $\Omega$, provided there are only finite number of them in any given tensor degree. In this way, every twisting map of an augmented complex determines a homomorphism into the bar-resolution of the algebra, in which this twisting map takes values.

On the other hand it is possible to extend the construction of the characteristic map of a twisting map $\psi$ in such a way that one obtaines a map from the twisted tensor product $K\otimes_\psi\Omega$ to the bar-complex of $\Omega$ with coefficients in $\Omega$, $B(\Omega,\,\Omega)$. To this end, we just apply the construction of $\tilde\psi$ to the elements of $K$ inside the tensor product $K\otimes\Omega$. Below we shall use a generalization of this observation in the case of a bitwisted tensor product (see definition \ref{bitwtenpr}.)

Let us discuss briefly another possible variation of this theory, namely when $K$ is a coalgebra and we define dual twisting map as $\varphi:K\otimes\Omega\to \Omega$, such that the differential
$$
d_\varphi(\alpha\otimes\omega)=d\alpha\otimes\omega+(-1)^{|\alpha|}\alpha\otimes d\omega+(-1)^{|\alpha^{(1)}|}\alpha^{(1)}\otimes\varphi(\alpha^{(2)}\otimes\omega)
$$
satisfies the relation $d_\varphi^2=0$ (see discussion after the definition \ref{twmapdf}.) Here $\alpha^{(1)}\otimes\alpha^{(2)}$ denotes the coproduct of $\alpha\in K$. If $\Omega$ is a coaugmented complex, i.e. there exist a chain map $\eta:\Bbbk\to\Omega$, then similarly to $\tilde\psi$ one can define the map
\begin{equation}
\label{tldphi}
\tilde\varphi:F(K)\to\Omega,\ \tilde\varphi=\sum_{n\ge1}\varphi_{(n)},
\end{equation}
where
$$
\varphi_{(n)}(\alpha_1\otimes\dots\alpha_n)=\varphi(\alpha_1\otimes\varphi(\alpha_2\otimes\varphi(\dots\varphi(\alpha_n\otimes1)))\dots).
$$
We put $\eta(1)=1\in\Omega$. It is an immediate consequence of the definitions that $\tilde\varphi$ is well-defined on $F(K)$, (and even on its completed version $\widehat F(K)$ if the series of elements in $\Omega$ on the right hand side converges.)

Suppose now that $K$ is a coalgebra and $\Omega$ is an algebra. One can ask, for which $\psi$ and $\varphi$ the two constructions of twisted tensor product $K\otimes_\psi\omega$ and $K\otimes_\varphi\Omega$ coincide. It follows that this happens in the case when $\psi$ and $\varphi$ are generated by the usual twisting cochain $\phi$ and only in this case. Indeed, in this situation the maps $\psi$ and $\varphi$ given by $\psi(\alpha)=\alpha^{(1)}\otimes\phi(\alpha^{(2)})$ and $\varphi(\alpha\omega)=\phi(\alpha)\omega$ verify the equations, imposed on the twisting maps above and give the same differential on $K\otimes\Omega$. And vice-versa, if $\psi:K\to K\otimes \Omega$ and $\varphi:K\otimes\Omega\to\Omega$ give the same differential $d_\psi\equiv d_\varphi$, then we have
$$
\psi^{(1)}(\alpha)\otimes\psi^{(2)}(\alpha)\omega=\alpha^{(1)}\otimes\varphi(\alpha^{(2)}\otimes\omega)
$$
for all $\alpha\in K$ and $\omega\in\Omega$. Substituting $1$ for $\omega$, we see that $\psi^{(1)}(\alpha)\otimes\psi^{(2)}(\alpha)=\alpha^{(1)}\otimes\phi(\alpha^{(2)})$, where $\phi(\alpha)=\varphi(\alpha\otimes1)$. Now, applying $\epsilon$ to the left leg in the tensor product we obtain $\varphi(\alpha\otimes\omega)=\phi(\alpha)\omega$. Finally, an easy computation shows that the map $\phi$ thus obtained verifies the usual equation of twisting cochains.

Due to this property, the characteristic maps, associated to twisting cochains have an important advantage:
\begin{lemma}
If $\psi$ is given by the formula $\psi(\alpha)=\alpha^{(1)}\phi(\alpha^{(2)})$ for a twisting cochain $\phi$, then the map $\tilde\psi$ is a homomorphism of differential graded coalgebras. Dually, if $\varphi(\alpha\otimes\omega)=\phi(\alpha)\omega$, then the map $\tilde\varphi$ is a homomorphism of differential graded algebras.
\end{lemma}
\begin{Proof}
It is enough to observe, that these maps are given by the formulae
\begin{align*}
\tilde\psi&(\alpha)=\sum_{n=0}^\infty\phi(\alpha^{(1)})\otimes\dots\otimes\phi(\alpha^{(n)}),\\
\tilde\varphi&(\alpha_1\otimes\dots\otimes\alpha_n)=\phi(\alpha_1)\dots\phi(\alpha_n).
\end{align*}
\end{Proof}\\ \\ \vspace{1mm} 
In general, if $K$ is not a coalgebra, but only an $A_\infty$\/-coalgebra and $\Omega$ is only $A_\infty$\/-algebra, one can ask for the analogs of this properties (see \cite{Smirn76, Smirn77, SmSOp} and references therein.) However, we shall not diverge in this direction in our paper. 

\begin{ex}
Let a Lie algebra $\g$ act on the coalgebra $K$. We can consider the following construction, dual to the construction of flat connections, see example \ref{example1}. Let $X_i,\ i=1,\dots,n$ be the linear base of $\g$. We let $\g$ act on the dual algebra of $K$ by conjugation, i.e. $X(\alpha^*)(\alpha)=\alpha^*(X(\alpha))$ for $\alpha^*\in K^*,\ \alpha\in K$ and $X\in\g$. Observe, that for any coalgebra $K$ its dual $K^*$ is a unital algebra (its unit is given by the counit $\epsilon$ on $K$.) Let $f^i:\Omega^*\to\Omega^{*+\mathrm{deg}\,X_i}$ be a collection of linear maps. We define a map
\begin{equation}
\label{eqduflcon}
\varphi(\alpha\otimes\omega)=\sum_i X_i(\epsilon)(\alpha)f^i(\omega)=\sum_i\epsilon(X_i(\alpha))f^i(\omega).
\end{equation}
(here we omit the multiplication sign in $\Omega$.) One can write down an equation, which is equivalent to $d_\varphi^2=0$: if we write down $\varphi$ as $X_i(\epsilon)\otimes f^i$, it will take the form
$$
\sum_i\{dX_i(\epsilon)\otimes f^i+(-1)^{|X_i|}X_i(\epsilon)\otimes df^i\}=\sum_{i,j}X_i(\epsilon)\cup X_j(\epsilon)\otimes f^i\circ f^j.
$$
\end{ex}

\vspace{1mm}
To conclude this paragraph, let us describe the conditions, under which the map \eqref{eqduflcon} will coincide with a map, defined in example \ref{example1}. Comparing the formulas, we see that first of all the maps $f^i$ should be equal to the multiplication by elements $\omega^i\in\Omega$ and second that the action of $\g$ should verify the following condition
\begin{equation}
\label{condg}
X(\alpha)=\alpha^{(1)}\epsilon(X(\alpha^{(2)})).
\end{equation}
In particular this condition is verified, when the action of $\g$ is induced by an action of a Lie group $G$, $K\times G\to K$, satisfying the equation
$$
(k^g)^{(1)}\otimes (k^g)^{(2)}=k^{(1)}\otimes (k^{(2)})^g.
$$
Here $k^g$ denotes the action of an element $g\in G$ on $k\in K$. It is easy to check that this is the case, when $K$ is algebra of differential (K\"ahler) forms on an algebraic group $G$: the action of $G$ is given by left multiplication on the argument of a function, and the comultiplication is determined by the group product:
$$
\alpha(g_1g_2)=\sum\alpha^{(1)}(g_1)\alpha^{(2)}(g_2).
$$
Thus we obtain the following important proposition:
\begin{prop}
If $K=\Omega_{DR}(G)$, $X_i,\ i=1,\dots,n$ is a linear base of $\g=lie(G)$ then for a commutative algebra $\Omega$ and arbitrary $\omega^i\in\Omega^1$ the map
$$
\psi(\alpha)=\sum_i\{\mathcal{L}_{X_i}(\alpha)\omega^i+\imath_{X_i}(\alpha)d\omega^i\}
$$
is twisting map in both senses. (Here $\lc_X$ is the Lie derivative and $\imath_X$ the contraction with the vector field X on the group.) In particular, the map $\tilde\psi$, defined by this $\psi$ is a coalgebra homomorphism.
\end{prop}
On the other hand, we have already defined a bialgebra (in fact, a Hopf algebra) homomorphism $\widetilde M_A:\Omega_{DR}(\ac(G))\to B(\Omega_{DR}(X))$ for any connection $A$ on the trivial bundle $P=G\times X$ over $X$ (more precisely, the image of $\widetilde M_A$ should lie in the normalized bar-resolution $N(\Omega_{DR}(X))$). The following proposition relates the maps of the present section with $\widetilde M_A$:
\begin{prop}
\label{compmono}
Let $A=\sum_i A^i\otimes X_i$ be a connection on the trivial bundle $G\times X$, where $X_i$ are the basis elements in the Lie algebra $\g$ of a group $G$, and $A^i\in\Omega_{DR}(X)$. Then the map $\widetilde M_A$ coincides with the map $\tilde\psi$ for the twisting map $\psi$ determined by $A$.
\end{prop}
The proof is obtained by mere inspection of the formulas and we omit it. Observe that in this way we obtain a new proof of the coalgebraic part of proposition \ref{prop1}. In order to show that $\tilde\psi$ is a homomorphism of algebras too, it is enough to notice that $\psi$, defined above, is a differentiation of $\Omega_{DR}(\ac(G))$ with values in $B(\Omega_{DR}(X))$.

\subsection{Gauge transformations of twisting maps}
\label{sectnnext}
Let us now discuss the equivalence relations on the set of twisting maps (in particular, on the generalized flat connections.) The relation we need should verify an evident property: two equivalent twisting maps $\psi$ and $\psi'$ should give equivalent structures of twisted tensor product on $K\otimes\Omega$ i.e. the structures, that are in some sense isomorphic. Here we don't explain, what precisely are the structures on the twisted tensor products that we need. We shall do it more accurately and with all the necessary details only in the case of twisting cochains.

So we now suppose that $K$ is a (co)chain complex, and $\Omega$ -- a DG algebra. Suppose, we are given two twisting maps, $\psi$ and $\psi'$ from $K$ to $\Omega$. In order to define equivalence between them, let us begin with a morphism $K\otimes\Omega\to K\otimes\Omega$, given by the formula
\begin{equation}
\label{gaugetr}
\alpha\otimes\beta\mapsto c^{(1)}(\alpha)\otimes c^{(2)}(\alpha)\beta,
\end{equation}
where $c:K\to K\otimes\Omega,\ c(\alpha)=c^{(1)}(\alpha)\otimes c^{(2)}(\alpha)$ is a homogeneous linear map of degree 0. There are two evident conditions, that one should impose on $c$:
\begin{enumerate}
\item{the map \eqref{gaugetr} should intertwine the differentials, induced by the twisting maps $\psi$ and $\psi'$, i.e. it should give us a cochain map from $K\otimes_{\psi}\Omega$ to $K\otimes_{\psi'}\Omega$;}
\item{the map \eqref{gaugetr} should be invertible in the class of all maps of the same kind, i.e. there should exist another map $c':K\to K\otimes\Omega$ such that the composition of the transformations \eqref{gaugetr}, determined by $c$ and $c'$ gives identity.}
\end{enumerate}
It is not difficult to write down a formula, expressing the first condition in the terms of components of $c$: one should have for all $\alpha\in K$
\begin{equation}
\label{gaugec}
\begin{split}
dc^{(1)}(\alpha)\otimes c^{(2)}(\alpha)&+(-1)^{c^{(1)}(\alpha)}c^{(1)}(\alpha)\otimes dc^{(2)}(\alpha)+\psi^{(1)}(c^{(1)}(\alpha))\otimes\psi^{(2)}(c^{1}(\alpha))c^{(2)}(\alpha)\\
                                       &c^{(1)}(d\alpha)\otimes c^{(2)}(d\alpha)+c^{(1)}({\psi'}^{(1)}(\alpha))\otimes c^{(2)}({\psi'}^{(1)}(\alpha)){\psi'}^{(2)}(\alpha),
\end{split}
\end{equation}
(we don't distinguish the differentials in $K$ and $\Omega$ with subscripts here.) One can abbreviate this to
\begin{equation}
\label{gaugeshort}
dc=c\smile\psi'-\psi\smile c,
\end{equation}
where the product $\smile$ is determined by the formula
\begin{equation}
\label{smile}
(a\smile b)(\alpha)=a^{(1)}(b^{(1)}(\alpha))\otimes a^{(2)}(b^{(1)}(\alpha))b^{(2)}(\alpha),
\end{equation}
for any two linear maps $a$ and $b:K\to K\otimes\Omega$ (observe that $a\smile b$ is also a map of this sort), and $dc=(d\otimes 1+1\otimes d)\circ c-c\circ d$. It is also worth to note that the invertibility condition can be expressed in the terms of the product \eqref{smile} as follows: $c$ is invertible, if there exists a map $c':K\to K\otimes\Omega$ for which one has $c\smile c'=c'\smile c =Id_K\otimes 1$ ($1$ is the unit in the algebra $\Omega$).

In general it is not clear how one can check the invertibility of a map $c$ in the case of arbitrary $K$ and $\Omega$. However, if we regard for example the twisting maps, defined by the action of a Lie algebra of a (compact) Lie group or their generalizations (see examples \ref{example1} and \ref{expoly} above,) it is natural to use the gauge transformations, induced by the action of the Lie group too. In this case we define $c$ to be the map
$$
\alpha\mapsto\sum_i\alpha^{g^i}\otimes f_i,
$$
where $g^i\in G$ and $f_i\in (\Omega^0)^\times$, and the symbol $A^\times$ denotes the invertible part of an algebra $A$. In particular, one can take $f_i=1$. More generally, if $\Omega$ is the algebra of differential forms on a manifold $X$, one can take $c$ to be an infinite sum of the elements of the form $g\otimes f$, provided it converges in a suitable norm to a well-defined map $X\to G$. In this case equation \eqref{gaugeshort} takes the form of the usual gauge transformation.

However, in example \ref{expoly} it is also possible to consider the gauge transformations, induced by the action of the enveloping algebra of the Lie algebra. Then $c$ will be equal to the sum of the form
\begin{equation}
\label{formgauge}
c=\sum_I X_I\otimes c^I,
\end{equation}
where $\{X_I\}$ is a basis in the $U\g$ and $c^I\in\Omega$ are such elements that the total degree of the tensor product is equal to $0$. Then the invertibility of $c$ is the same as its invertibility as an element in the tensor product of two algebras. Observe that if we don't care about the convergence of the series defining \eqref{formgauge} (e.g. under the conditions of the lemma \ref{polyex},) then invertibility holds for example, if the first term in this series is equal to $1\otimes 1$.

The collection of all such gauge transformations gives an equivalence relation on the set of all twisting maps for the pair $K,\ \Omega$ (one easily checks that a composition of the maps \eqref{gaugetr} induced by $c$ and $c'$ is equal to the transformation, induced by the map $(c\otimes 1)\smile c'$). It is the set of equivalence classes of this relation, which will be the subject of our attention in the rest of this paper. Consider the following example.

\begin{ex}\rm
\label{counterex}
Let us show that all the twisting maps from example \ref{example1} and from lemma \ref{polyex} (example \ref{expoly}) are gauge equivalent to the trivial twisting map $\psi'$ which sends everything to zero.

In the case of \ref{example1}, it is enough to consider the map
\begin{equation*}
\begin{split}
c(k)&=k\otimes 1+\sum_i(-1)^{\epsilon_i}I_i(k)\otimes A^i+\frac12\sum_{i,j}(-1)^{\epsilon_j+\epsilon_i}I_jI_i(k)\otimes A^j\wedge A^i\\
&\quad\qquad+\frac16\sum_{i,j,l}(-1)^{\epsilon_l+\epsilon_j+\epsilon_i}I_lI_jI_i(k)\otimes A^l\wedge A^j\wedge A^i+\dots.
\end{split}
\end{equation*}
In a more conceptual coordinate-free form, this formula can be written as
$$
c=\exp{(\sum_i I_i\otimes A^i)},
$$
for the formal power series exponent. It is clear, that the map $c$ being exponent of an element, is invertible. A direct computation, based on the formulas
$$
d(I_i)=X_i,\ \quad\ d\omega^i=-C^i_{kl}\omega^k\wedge\omega^l,\ \mbox{ and }\ F^i=dA^i-C^i_{kl}A^k\wedge A^l,
$$
(here $\omega^i$ are the coefficients of the Maurer-Cartan form and the formula, involving them is just the Maurer-Cartan equation), shows that
$$
-\tilde A\smile c=dc.
$$
That is $c$ is a gauge transformation, connecting flat connection $\tilde A$ of example \ref{example1} with a trivial twisting map. A still more conceptual way of looking onto this phenomenon is as follows. Suppose that the principal bundle $P$ is trivial. Consider the following automorphism of the algebra of differential forms on $P$, $\DR{P}\cong\DR{G}\otimes\DR{X}$:
$$
f\omega^{i_1}\wedge\dots\wedge\omega^{i_p}\otimes\alpha\mapsto f\exp{(\sum_i I_i\otimes A^i)}(\omega^{i_1}\wedge\dots\wedge\omega^{i_p})\otimes\alpha.
$$
It is easy to show that this map intertwines the usual differential in $\DR{P}$ with the "`twisted covariant derivation"', determined by the twisting cochain: $d+\sum_i\{A^i\otimes X_i-F^i\otimes I_i\}$.

Now we come to the twisting maps given by the lemma \ref{polyex}. Assume that all the conditions of this lemma are fulfilled. We should find a twisting cochain $c$ of the form \eqref{formgauge} such that
$$
dc=-\psi\smile c.
$$
Observe that in this case the map $\psi\smile c$ is equal (up to a sign) to the product of the series which determine $c$ and $\psi$. Then we can reason by induction as we did in the proof of the lemma \ref{polyex}. So we put $c_0=1\otimes 1$ and use the acyclicity of the tensor product $rU\g\otimes\Omega(n)$ to solve the further equations. Invertibility of the element we obtain in this way can be proved in a usual way as the invertibility of a formal power series beginning with $1$.
\end{ex}

\vspace{1mm}
An important property of the introduced notion of the gauge transformation is the following
\begin{prop}
\label{someprop1}
Let the twisting maps $\psi$ and $\psi'$ be gauge equivalent (i.e. there exists a gauge transformation, which relates $K\otimes_\psi\Omega$ and $K\otimes_{\psi'}\Omega$). Suppose that the conditions of the lemma \ref{lembar} are satisfied, then the maps $\tilde\psi,\tilde{\psi'}:K\to B(\Omega)$ are homotopic. The same is true for arbitrary $\psi$ and $\psi'$, if we consider the completion of $B(\Omega)$ with respect to the natural filtration by degree.
\end{prop}
\begin{Proof}
Let $c$ be the gauge transformation, which relates $\psi$ and $\psi'$. We shall use it to define the chain homotopy, relating $\tilde\psi$ and $\tilde\psi'$. First of all, let us rewrite the equation \eqref{gaugeshort}:
$$
\psi-\psi'=\psi'\smile(c-1)-(c-1)\smile\psi+d(c-1).
$$
Here the map $1:K\to K\otimes\Omega$ is given by $1(\alpha)=\alpha\otimes1$. We used the evident properties of this map $d1=0$ and $1\smile a=a\smile1$ for all $a:K\to K\otimes\Omega$. Next we put $\hat c=c-1$ and define the maps $h_{(1)}$ and $h_{(2)}$ by the formulas
\begin{align*}
h_{(1)}&=(\epsilon\otimes1)\circ\hat c\\
\intertext{and}
h_{(2)}&=(\epsilon\otimes 1\otimes 1)\circ((\hat c\otimes 1)\circ\psi+(\psi'\otimes 1)\circ\hat c))\\
\intertext{or, in more explicit terms}
h_{(1)}(\alpha)&=\epsilon(\hat c^{(1)}(\alpha))\hat c^{(2)}(\alpha)\\
\begin{split}
h_{(2)}(\alpha)&=
\epsilon\bigl(\hat c^{(1)}({\psi}^{(1)}(\alpha))\bigr)\hat c^{(2)}({\psi}^{(1)}(\alpha))\otimes{\psi}^{(2)}(\alpha)\\
&\quad+\epsilon\bigl({\psi'}^{(1)}(\hat c^{(1)}(\alpha))\bigr){\psi'}^{(2)}(\hat c^{(1)}(\alpha))\otimes\hat c^{(2)}(\alpha).
\end{split}
\end{align*}
It follows from \eqref{smile} (the definitions of $\psi_{(1)}$ and $\psi'_{(1)}$ are given above, see the discussion preceding lemma \ref{lembar}):
\begin{equation*}
\begin{split}
h_{(1)}(d\alpha)&-\bigl(d(h_{(1)}(\alpha))+b(h_{(2)}(\alpha))\bigr)=\psi_{(1)}(\alpha))\\
                &=(\epsilon\otimes1)(-d\hat c(\alpha)-(\psi'\smile\hat c)(\alpha)+(\hat c\smile\psi)(\alpha))\\
                &=(\epsilon\otimes1)(\psi'(\alpha)-\psi(\alpha))=\psi'_{(1)}(\alpha)-\psi_{(1)}(\alpha).
\end{split}
\end{equation*}
Now put
\begin{equation*}
h_{(3)}=(\epsilon\otimes 1^{\otimes 3})\circ((\hat c\otimes 1^{\otimes2})\circ(\psi\otimes1)\circ\psi
+(\psi'\otimes 1^{\otimes2})\circ(\hat c\otimes1)\circ\psi+(\psi'\otimes 1^{\otimes2})\circ(\psi'\otimes1)\circ\hat c).
\end{equation*}
It is easy to check that the following equation is true
$$
h_{(2)}(d\alpha)-\bigl(d(h_{(2)}(\alpha))+b(h_{(3)}(\alpha))\bigr)=\psi'_{(2)}(\alpha)-\psi_{(2)}(\alpha).
$$
Finally, we put $H=\sum_{k=0}^\infty h_{(k)}$, where
$$
h_{(k)}=(\epsilon\otimes1^{\otimes k})\Bigl(\sum_{i=1}^{k-1}(\psi'\otimes1^{\otimes k})\circ\dots\circ(\psi'\otimes1^{\otimes(k-i)})\circ(\hat c\otimes1^{\otimes(k-i-1)})\circ(\psi\otimes1^{\otimes(k-i-2)})\circ\dots\circ\psi\Bigr).
$$
Then by induction one proves that $H\circ d-(d+b)\circ H=\tilde\psi'-\tilde\psi$, i.e. $H$ is the chain homotopy, connecting $\tilde\psi$ and $\tilde\psi'$.
\end{Proof}\\ \\ \vspace{1mm}
\begin{rem}\rm
\label{remgaugetw}
Of course, one can develop an equivalent of the theory of gauge transformations in the case when $K$ is a coalgebra, $\Omega$ -- an arbitrary cochain complex and the twisting is defined with the help of linear maps $\phi:K\otimes\Omega\to\Omega$. One can show that the maps $\tilde\phi$ and $\tilde\phi'$ are homotopic, when $\phi$ and $\phi'$ are gauge equivalent (see equations \eqref{tldphi} and the corresponding discussion.) However, we shall not do it here. We shall confine ourselves to few remarks on the theory of gauge transformations in the case of twisting cochains.

In effect when the twisting map $\psi$ is induced by a twisting cochain, one can give a more detailed description of the gauge transformations. Let $K$ be a connected differential graded coalgebra (i.e. coalgebra, whose degree $0$ part is isomorphic to $\Bbbk$) and $\Omega$ - an algebra. Following the exposition of \cite{Smirn76} we consider the space $\hc(K,\Omega)=\Hom(K,\Omega)$ of linear maps $c:K^*\to \Omega^*$ of degree $0$, verifying the condition $c(K^0)=0$. There's a natural multiplication $\cup$ on $\hc(K,\Omega)$, induced from the comultiplication on $K$ and the multiplication on $\Omega$ (the convolution product.) This product, together with the usual differential of maps (commutator with the differentials of the domain and the range) turns the space $\hc(K,\,\Omega)$ into a differential graded algebra. We shall twist the product in this algebra a little bit, in order to make it more suitable for our purposes. Put
\begin{equation}
\label{eqcirc}
c\circ c'=c+c'-c'\cup c.
\end{equation}
Now since $c(K^0)=0$, there always exists $c'$ such that $c\circ c'=0$ (one can use \eqref{eqcirc} to define $c'$ recursively, see \cite{Smirn76}.)

There's a natural action of $\hc(K,\Omega)$ with the little circle product on it on the space of twisting cochains $\mathcal T(K,\Omega)$. It is given by $\phi\circ c=\phi'$ where $\phi'$ is determined by the following equation
\begin{equation}
\label{eqcphi}
\phi'=\phi-c\cup\phi+\phi'\cup c+dc.
\end{equation}
Once again the condition $c(K^0)=0$ allows one to solve the equation \eqref{eqcphi} recursively (see \cite{Smirn76} again.) It is a matter of direct computations to show that $\phi'$ verifies the equation
\eqref{flatc} i.e. that it is a twisting cochain.

Instead of the condition $c(K^0)=0$ and connectivity of $K$ one can ask for the existence of inverse elements with respect to the $\cup$\/-product for all maps $1-c$, $c\in\hc(K,\Omega)$, where $1(\alpha)=\epsilon(\alpha)\cdot1_\Omega$. This is the case, e.g. if the product in $\hc(K,\Omega)$ is nilpotent. Or else this happens for the set of algebra homomorphisms from a Hopf algebra $K$ to an algebra $\Omega$. \hfill$[]$\end{rem}

\subsection{The gauge bundle and twisting cochains}
\label{sectnext4}
This paragraph is devoted to a description of a map $\tilde{\hat\phi}$ similar to the homological monodromy map of the previous section, section \ref{sect11} (see the final remark of the paragraph \ref{sect17}.) Thus it should be a map from (a model of) the algebra of differential forms on the gauge bundle to the normalized Hochschild complex of a differential graded algebra $\Omega$, which we should regard as the "`model"' of the base. Our purpose is to derive a construction of $\tilde{\hat\phi}$ from the twisting map $\phi$. For the sake of simplicity we shall restrict our attention to the case where the twisting map is in fact induced by a twisting cochain. The main difficulty that one has to overcome is that the natural model of the gauge bundle (see the discussion after the proof of proposition \ref{propleftright} and the corresponding footnote) can hardly be used as the domain of this map: one has to introduce a more convenient complex, one we denote by $K\hat{\otimes}_\phi\Omega$ below, and then prove that it can be regarded as a model for the gauge bundle.

So, let $K$ (resp. $\Omega$) be a differential graded coalgebra (resp. algebra), let $K$ be coaugmented and $\phi:K^*\to\Omega^{*+1}$ be a twisting cochain. Instead of the usual twisted tensor product $K\otimes_\phi\Omega$ one can consider the following generalization thereof, which one can call \emph{bitwisted}:
\begin{df}
\label{bitwtenpr}
The tensor product $K\hat\otimes_{\phi}\Omega$ is as before linearly isomorphic to the graded tensor product of $K$ and $\Omega$ and the differential is given by the formula:
\begin{equation*}
\begin{split}
\hat d_{\phi}(\alpha\otimes\beta)&=d(\alpha)\otimes\beta+(-1)^{|\alpha|}\alpha\otimes d\beta\\
                                                            &\quad+\alpha^{(1)}\otimes\phi(\alpha^{(2)})\beta+(-1)^{(|\alpha^{(1)}|+1)(|\alpha^{(2)}|+|\beta|)}\alpha^{(2)}\otimes\beta\phi(\alpha^{(1)}).
\end{split}
\end{equation*}
\end{df}
An easy calculation now shows that $\hat d_\phi^2=0$. The following statement  is the main theorem of this paragraph. It is a direct analogy of the Brown's theorem, (see \cite{Br59}.) We give a brief proof of it for the sake of the self-containedness of our paper:
\begin{prop}
\label{propleftright}
Let the usual $\phi$\/-twisted tensor product of $K$ and $\Omega$, $K\otimes_\phi\Omega$, be a model of the principal bundle $P$ (i.e. its cohomology is isomorphic to that of $P$ as $H^*(X)$\/-module and $H^*(G)$\/-comodule, see discussion below.)  Suppose that $\Omega$ is unitary and $K$ is connected (i.e. $K^0$ is isomorphic to the ground field). Then the cohomology of the complex $K\hat\otimes_\phi\Omega$ is isomorphic to the cohomology of the associated gauge bundle $P\times_{Ad}G$.
\end{prop}
\begin{Proof}
First of all, consider the following intermediate complex:
$$
C(K,\,\Omega;\,\phi)=\bigoplus_{p,q=0}^\infty (K\otimes \bar K^{\otimes^p}\otimes K)\otimes(\Omega\otimes\bar\Omega^{q}\otimes\Omega),
$$
where $\bar K=K/1$ (here $1$ is the image of $1\in\Bbbk$ under the coaugmentation of $K$,) $\bar\Omega=\Omega/1$, the grading is given by
$$
\mathrm{deg}\,a[k_1|\dots|k_p]b\otimes\alpha[\omega_1|\dots|\omega_q]\beta=\mathrm{deg}\,a+\mathrm{deg}\,b+\sum\mathrm{deg}\,k_i-p+\mathrm{deg}\,\alpha+\sum\mathrm{deg}\,\omega_j+q+\mathrm{deg}\,\beta,
$$
(here we use $|$ instead of the tensor signs.) Observe, that since we have assumed that $K$ is connected, there's no problem with negative degrees here. The differential in this complex shall be given by the formula
\begin{equation*}
\begin{split}
&d(a[k_1|\dots|k_p]b\otimes\alpha[\omega_1|\dots|\omega_q]\beta)\\
  &=da[k_1|\dots|k_p]b\otimes\alpha[\omega_1|\dots|\omega_q]\beta+\sum_{i=1}^p(-1)^{\epsilon_i}a[k_1|\dots|dk_i|\dots|k_p]\otimes\alpha[\omega_1|\dots|\omega_q]\beta\\
  &\qquad\qquad\qquad+(-1)^{\epsilon_p}a[k_1|\dots|k_p]db\otimes\alpha[\omega_1|\dots|\omega_q]\beta\\
  &\quad+(-1)^\epsilon\Bigl(a[k_1|\dots|k_p]b\otimes d\alpha[\omega_1|\dots|\omega_q]\beta+\sum_{j=1}^q(-1)^{\eta_j}a[k_1|\dots|k_p]\otimes\alpha[\omega_1|\dots|d\omega_j|\omega_q]\beta
\end{split}
\end{equation*}
\begin{equation*}
\begin{split}
  &\qquad\qquad\qquad+(-1)^{\eta_q}a[k_1|\dots|k_p]\otimes\alpha[\omega_1|\dots|\omega_q]d\beta\Bigr)\\
  &\quad+a^{(1)}[a^{(2)}|k_1|\dots|k_p]b\otimes\alpha[\omega_1|\dots|\omega_q]\beta+\sum_{i=1}^p(-1)^{\epsilon_i}a[k_1|\dots|k_i^{(1)}|k_i^{(2)}|\dots|k_p]b\otimes\alpha[\omega_1|\dots|\omega_q]\beta\\
  &\qquad\qquad\qquad+(-1)^{\epsilon_p}a[k_1|\dots|k_p|b^{(1)}]b^{(2)}\otimes\alpha[\omega_1|\dots|\omega_q]\beta\\
  &\quad+(-1)^\epsilon\Bigl(a[k_1|\dots|k_p]b\otimes1[\alpha\omega_1|\dots|\omega_q]\beta+\sum_{j=1}^{q-1}(-1)^{\eta_j}a[k_1|\dots|k_p]b\otimes\alpha[\omega_1|\dots|\omega_j\omega_{j+1}|\omega_q]\beta\\
  &\qquad\qquad\qquad\qquad+(-1)^{\eta_q}a[k_1|\dots|k_p]b\otimes\alpha[\omega_1|\dots|\omega_q\beta]1\Bigr)\\
  &\qquad\qquad\qquad+(-1)^{\chi_1}a^{(2)}[k_1|\dots|k_p]b\otimes\alpha[\omega_1\dots|\omega_q]\beta\phi(a^{(1)})\\
  &\qquad\qquad\qquad\qquad\qquad\quad+(-1)^{\chi_2}a[k_1|\dots|k_p]b^{(1)}\otimes\phi(b^{(2)})\alpha[\omega_1|\dots|\omega_q]\beta.
\end{split}
\end{equation*}
Here $\epsilon_i=|a|+\sum_{l=1}^{i-1}|k_l|-i+1,\ \epsilon=|a|+\sum_{i=1}^p|k_i|-p+|b|,\ \eta_j=|\alpha|+\sum_{m=1}^{j-1}|\omega_m|+j-1,\ \chi_1=(|a^{(1)}+1|)(|a^{(2)}|+\sum_{i=1}^p|k_i|-p+|b|+|\alpha|+\sum_{j=1}^q|\omega_j|+q+|\beta|)$ and $\chi_2=|a|+\sum_{i=1}^p|k_i|=\epsilon_p$.

To prove the proposition, it is enough to show that the following two statements hold:
\begin{itemize}
\item\emph{The complex $C(K,\,\Omega;\,\phi)$ models the gauge bundle $P\times_{Ad}G$ (in particular, its cohomology is isomorphic to the cohomology of the gauge bundle).}
\item\emph{Complexes $C(K,\,\Omega;\,\phi)$ and $K\hat\otimes_\phi\Omega$ are homotopy equivalent.}
\end{itemize}

First, let us calculate the cohomology of the complex $C(K,\,\Omega;\,\phi)$. Let $\Omega{}_\phi\!\otimes K$ be the \emph{left}-twisted tensor product of $K$ and $\Omega$, i.e.  $\Omega{}_\phi\!\otimes K\cong\Omega\otimes K$ as linear spaces, and the differential is given by the formula
$$
{}_\phi d(\omega\otimes k)=d\omega\otimes k+(-1)^{|\omega|}(\omega\otimes dk+\omega\phi(k^{(1)})\otimes k^{(2)}).
$$
It is clear (see Brown's paper \cite{Br59} and Smirnov's book \cite{SmSOp}) that $\Omega{}_\phi\!\otimes K$ is a model for the cohomology of the principal bundle $\hat P=P\times_{R^{-1}} G$, where $R^{-1}$ is the left action of $G$ on itself, given by $R^{-1}_g(h)=hg^{-1}$. (We shall prove this fact in detail in the next section, see remark \ref{leftright}.) Observe, that $\hat P$ is diffeomorphic to $P$ as a differentialble manifold, while the structure group $G$ acts on $\hat P$ from the right (we assumed that it acted on $P$ from the left): $\hat p\cdot g=\widehat{g^{-1}\cdot p}$. So we conclude that $C(K,\,\Omega;\,\phi)$ is homotopy-equvalent to the realization of the \emph{cosimplicial simplicial DG algebra} $\DR{\mathscr C_\bullet(G,\,X;\,P)}$, where $\mathscr C_\bullet(G,\,X;\,P)$ is the simplicial cosimplicial space (in fact, manifold,) given by
$$
\mathscr C_{p,q}(G,\,X;\,P)=\hat P\times\underbrace{G\times\dots\times G}_p\times P\times\underbrace{X\times\dots\times X}_q.
$$
Here $p$ denotes the simplicial dimension and $q$ is the cosimplicial dimension of this space. The simplicial structure maps are given by the usual bar-construction morphisms, so that for a fixed $q$ we shall have $\mathscr B_\bullet(\hat P,\,G,\,P)\times X^{\times q}$, and the cosimplicial coface (resp. codegeneracy) maps will be given by diagonal embedings (resp. projections onto subfactors) in the cartesian powers of $X$. So according to the Bott-Segal theorem (see \cite{BottSeg},) cohomology of $C(K,\,\Omega;\,\phi)$ is isomorphic to the cohomology of the suitable realization of $\mathscr C_\bullet(G,\,X;\,P)$. On the other hand, realizing this space first in the simplicial direction, and then in the cosimplicial, and using the fact that $P$ is a free $G$\/-space we obtain first the cosimplicial space
$$
\mathscr F_q(P,\,X,\,\hat P)/G=(P\times\underbrace{X\times\dots\times X}_q\times\hat P)/\sim,
$$
where the relation $\sim$ is given by $(g\cdot p,x_1\dots,x_q,\hat p)\sim(p,x_1,\dots,x_q,\hat p\cdot g)$. Since the action of $G$ is free, we can commute the geometric realization and the factorization procedure, thus we conclude (using Anderson's theorem, see \cite{Anders},) that
$$
|\mathscr C_\bullet(G,\,X;\,P)|=\mathscr P(X;\,\hat P\times_G P).
$$
Here $\hat P\times_G P=\{(\hat p,\,p)\}/(\hat p\cdot g,p)\sim(\hat p,g\cdot p)$. If $\pi_1,\ \pi_2$ are the natural projections of $\hat P\times_G P$ on $X$. Then we have
$$
\mathscr P(X;\,\hat P\times_G P)=\{(\gamma,\,[\hat p,\,p])|\gamma:[0;\,1]\to X,\ \gamma(0)=\pi_1([\hat p,\,p]),\ \gamma(1)=\pi_2([\hat p,\,p])\}.
$$
Now it is enough to observe that $\mathscr P(X;\,\hat P\times_G P)$ is homotopy equivalent to the pullback of $\pi_1\times\pi_2:\hat P\times_G P\to X\times X$ under the diagonal map $diag:X\to X\times X$ (it is enough to pull all the paths $\gamma$ to their origins and use the homotopy lifting property of $P$,) which is isomorphic to $P\times_{Ad} G$.

In order to prove the second statement, let us consider the complex $F(K,\,\Omega;\,\phi)$:
$$
F(K,\,\Omega;\,\phi)=\bigoplus_p K\otimes \bar K^{\otimes p}\otimes K\otimes A,
$$
with differential
\begin{equation*}
\begin{split}
d(a[k_1|\dots|k_p]b\otimes\omega)&=da[k_1|\dots|k_p]b\otimes\omega+\sum_{i=1}^p(-1)^{\epsilon_i}a[k_1|\dots|dk_i|\dots|k_p]\otimes\omega\\
                                                               &\quad+(-1)^{\epsilon_p}a[k_1|\dots|k_p]db\otimes\omega+(-1)^\epsilon a[k_1|\dots|k_p]b\otimes d\omega\\
                                                               &\quad+a^{(1)}[a^{(2)}|k_1|\dots|k_p]b\otimes\omega+\sum_{i=1}^p(-1)^{\epsilon_i}a[k_1|\dots|k_i^{(1)}|k_i^{(2)}|\dots|k_p]b\otimes\omega\\
                                                               &\quad+(-1)^{\epsilon_p}a[k_1|\dots|k_p|b^{(1)}]b^{(2)}\otimes\omega+(-1)^{\epsilon_p}a[k_1|\dots|k_p]b^{(1)}\otimes\phi(b^{(2)})\omega\\
                                                               &\quad+(-1)^{\chi_3}a^{(2)}[k_1|\dots|k_p]b\otimes\omega\phi(a^{(1)}).
\end{split}
\end{equation*}
Here all the gradings, signs and notation are taken from the definition of $C(K,\,\Omega;\,\phi)$, and $\chi_3=(|a^{(1)}|+1)(|a^{2)}|+\sum_{i=1}^p|k_i|-p+|b|+|\omega|)$. There is a map of complexes $C(K,\,\Omega;\,\phi)\to F(K,\,\Omega;\,\phi)$, given by projection of $C(K,\,\Omega;\,\phi)$ onto the $q=0$ part followed by multiplication $\Omega\otimes\Omega\to\Omega$. We claim, that this map induces an isomorphism in cohomology. To this end we introduce filtrations on both complexes: by degree $q$ in $C(K,\,\Omega;\,\phi)$ and a trivial filtration (concentrated in $0$) in $F(K,\,\Omega;\,\phi)$. The map we introduced respects these filtrations and induces an isomorphism at the first level of the corresponding spectral sequences: recall, that we assumed, that $\Omega$ has unit, hence its bar-resolution is contractible.

Finally, there's a map of complexes $K\hat{\otimes}_\phi\Omega\to F(K,\,\Omega;\,\phi)$, given by the comultiplication in $K$ on the left followed by emedding on the $p=0$ part. By the arguments, similar to what we just used (consider spectral sequences,) we prove that this embedding induces an isomorphism in cohomology.
\end{Proof}\\

\vspace{-4mm}
\begin{rem}\rmfamily
\label{reminvprop}
One can easily show by a slight generalization of the reasoning in proposition \ref{someprop1} that if $\phi$ and $\phi'$ are gauge equivalent twisting cochains, then the corresponding maps $\tilde{\hat\phi}$ and $\tilde{\hat{\phi'}}$ are chain homotopic. In particular, we conclude that the homological monodromy map $\widetilde T_A$ associated to the gauge bundle of a trivial bundle, described in the previous section (paragraph \ref{sect17},) gives the same map on the level of gauge bundles and free loops, as the morphism, induced by the twisting map $\phi_A$, described in example \ref{example1}, since they both are equivalent to trivial maps (this follows from the proposition \ref{compmono} and its natural generalization to the gauge bundles.) Below (section \ref{sect40} and theorem \ref{omagn}) we shall show that similar statement holds for the gauge bundles of arbitrary bundles.
\hfill$[]$\end{rem}
Let us finally make a couple of observations concerning the situation, when $\Omega$ is not connected. In this case it is necessary to consider the normalized variants of the complexes $C(K,\,\Omega;\,\phi)$, etc (i.e. factorized by the subcomplex, generated by the 0-dimensional part of $\Omega$.) In fact, since we assume that the base $X$ is connected (even 1-connected) and $\Omega$ is a model for $X$, it follows that $\Omega$ is homotopy-equivalent to a connected algebra, and hence the $0$\/-dimensional part of $\Omega$ can be factorized. It is now evident that the formula
\begin{equation}
\label{eqloopmap}
k\otimes\omega\mapsto\sum_n\phi(k^{(1)})\otimes\dots\otimes\phi(k^{(n)})\otimes\omega
\end{equation}
(or if we use $|$ instead of $\otimes$, the right hand side will look as $\sum_n[\phi(k^{(1)})|\dots|\phi(k^{(n)})]\omega$) determines a chain map  $\tilde{\hat\phi}:K\hat\otimes_\phi\Omega\to NH(\Omega)$. It is also evident that all the homotopy equivalence properties of the ordinary maps, described in remark \ref{reminvprop}, are valid in this setting as well.
\begin{rem}\rm
\label{actiong}
One can introduce additional structures on both complexes, $K\hat\otimes_\phi\Omega$ and $NH(\Omega)$. Recall that the gauge bundle $P\times_{Ad}G$ is the bundle of groups, so that one can define the product of two (global) sections of $P\times_{Ad}G$. Geometrically, this structure is determined by the map
$$
(P\times_{Ad}G)\times_X(P\times_{Ad}G)\to P\times_{Ad}G,
$$
verifying evident associativity conditions. Similar map can be defined for the free loop space $\mathcal LX$ of $X$, when we regard it as a fibre bundle over $X$ with respect to the evaluation map $e(\gamma)=\gamma(0)$ (this time associatiovity is replaced with homotopy associativity properties.)

If we want to find an analog of this map on the level of algebraic models we can first of all demand that the model $A$ of $P\times_{Ad}G$ that we consider should be a differential graded module over a suitable algebra $\Omega$, corresponding to $X$, and second, that there be a map
$$
A\to A\otimes_\Omega A,
$$
verifying the usual coassociativity conditions (and similarly for the model of $\mathcal LX$.) More genrally one can speak about homotopy analogs of all these maps, i.e. about the $A_\infty$\/-algebras and coalgebras, about their modules etc., but this will bring us far beyond the modest purposes of the present paper.

In our case one can easily introduce the necessary maps on both $K\hat\otimes_\phi\Omega$ and $NH(\Omega)$, \textit{when $\Omega$ is commutative}. Namely we put
\begin{align}
\label{coprodom1}
k\hat\otimes\omega&\mapsto (-1)^{|k^{(2)}||\omega|}(k^{(1)}\hat\otimes\omega)\otimes_\Omega(k^{(2)}\hat\otimes1)\\
\intertext{and}
\label{coprodom2}
[\omega_1|\dots|\omega_n]a&\mapsto\sum_{i=0}^n(-1)^{\eta_i}([\omega_1|\dots|\omega_i]a)\otimes_\Omega([\omega_{i+1}|\dots|\omega_n]1).
\end{align}
Here $\eta_i=|a|(|\omega_{i+1}|+\dots+|\omega_n|+n-i)$ and we let $\Omega$ act on $K\hat\otimes_\phi\Omega$ by multiplication in the $\Omega$ factor, and on $NH(\Omega)$ by multiplication at the last tensor place. In case, when $\Omega$ is not commutative, but only homotopy-commutative, one should use the higher homotopy maps. In effect, as it will be shown in the next section (see also \cite{SmSOp},) if $\Omega$ and $\Omega'$ are homotopy-equivalent algebras, one can use the $A_\infty$\/-map $\mathcal P:\Omega\Rightarrow\Omega'$, which establishes this quasi-isomorphism, to determine a new twisting cochain $\mathcal P\circ \phi:K\to\Omega'$ (see \cite{SmSOp} and section \ref{sectinfty}, remark \ref{inftymaps} below,) which corresponds to $\phi$ under the homotopy equivalence. Thus we can always assume that $\Omega$ is commutative. It is now a matter of simple calculation to show that the map $\tilde{\hat\phi}$ intertwines the coproduct structures \eqref{coprodom1} and \eqref{coprodom2}.

Moreover, if $\phi,\ \phi'$ are two equivalent cochains and $c:K\to\Omega$ (with commutative $\Omega$) is an invertible map, that establishes this equivalence, i.e. the following equality holds $\phi'\smile c-c\smile\phi=dc$, then we can consider the following isomorphism of linear spaces:
\begin{align*}
\hat c:K\hat\otimes_\phi\Omega&\to K\hat\otimes_{\phi'}\Omega\\
k\otimes\omega&\mapsto (-1)^{(|k^{(1)}|+1)(|k^{(2)}|+|k^{(3)}|+|\omega|)}k^{(2)}\otimes c(k^{(3)})\omega c^{-1}(k^{(1)}),
\end{align*}
the inverse of this map being equal to $\widehat{c^{-1}}$. A straightforward calculation shows that this map commutes with the differentials. If the algebra $\Omega$ is commutative, then the map $\hat c$ commutes with the coproducts in $K\hat\otimes_{\phi}\Omega$ and $K\hat\otimes_{\phi'}\Omega$ determined by equation \eqref{coprodom1}.  Similarly one can show that the map $\tilde{\hat\phi}$ and $\tilde{\hat{\phi'}}\circ\hat c$ are homotopic. The homotopy is given by the map
\begin{equation}
\label{hombad}
\hat H(k\otimes a)=\sum_{n=0}^\infty\sum_{l=1}^{n-1} \phi(k^{(1)})\otimes\dots\otimes\phi(k^{(l-1)})\otimes c^{-1}(k^{(l)})\otimes \phi'(k^{(l+1)})\otimes\dots\otimes\phi'(k^{(n-1)})\otimes c(k^{(n)})a.
\end{equation}
\hfill$[]$\end{rem}

In fact, if $K$ is a Hopf algebra (i.e. if it is equipped with the multiplication and antipodal map,) it would be quite illuminating to investigate the relation of the complex $K\hat{\otimes}_\phi\Omega$ with the complex, induced by the same twisting cochain $\phi$ from the adjoint coaction of $K$ on itself: we let $ad(\alpha)=\alpha^{(2)}\otimes S(\alpha^{(1)})\alpha^{(3)}$ and consider the twisted tensor product of $\Omega$ and $K_{ad}$, which is $K$ regarded as a right $K$\/-comodule with respect to the adjoint coaction. The latter complex is a natural model for the gauge bundle in the sense of Brown's paper: we just replace the group by a map on which it acts. Another important problem is to construct a map from $K_{ad}\otimes_{\phi}\Omega$ to $NH(\Omega)$ and compare this map with the map $\tilde{\hat\phi}$. Observe that the direct generalization of the map $\tilde\phi$ to $\phi_{ad}$ gives trivial result; in effect even the range of the map is different from what the $\tilde{}$\/-construction naturally gives: it should be equal to the the normalized Hochschild complex, not the bar-resolution.

Another important question is concerned with the possibility to introduce a multiplication in $K\hat{\otimes}_{\phi}\Omega$, if $K$ is a Hopf algebra. This question is related to the properties of the twisting cochain with respect to the multiplication in $K$.


\section{Twisting cochains and characteristic classes of principal bundles}
\label{sect30}
In this section we give few constructions of the twisting cochains (more generally, twisting maps,) associated with an \emph{arbitrary} principal bundle $P$ (in general, non-trivial) with base $X$. Unless the opposite is explicitly stated we assume that $X$ is a closed connected manifold. We shall also show the relation of our constructions and the usual definitions of Chern classes as given by the Chern-Weil homomorphism, associated with a connection on the principal bundle. This chapter is closely related to the paper \cite{mypapLF}, where we give few details, omitted here and also discuss the relation of our construction and some previously known results.

\subsection{Twisting cochain and \v Cech cohomology}
\label{sect31}
The domain of the twisting cochain we want to construct now is the coalgebra (even Hopf algebra) of (polynomial) functions on the structure group $G$, we assume that this group is algebraic, so that we can avoid the problems with completion of tensor products in the definition of diagonal. This cochain takes values in the algebra of \v{C}ech cochains on $X$ with values in the sheaf of de Rham differential forms on $X$. Other constructions of the
twisting cochains and the relation of these constructions and the notion of the twisting maps, see previous section, will be discussed in the following paragraphs. Our ideas are based on the work of Shih, \cite{Shih}.

Let $P\stackrel{\pi}{\to}X$ be a principal bundle with structure group $G$ and base $X$. We suppose that $G$ is a compact (algebraic) Lie group and $X$ --- a closed $C^\infty$\/-manifold and that $G$ acts on $P$ from the \emph{left}. Let $\mathcal U=\{U_\alpha\}_{\alpha\in A}$ be an open cover of $X$, for an ordered set $A$. We can even suppose that $A$ is finite, for instance, if $X$ is compact. We assume that $P_{|U_\alpha}$ is trivial and fix once and forever the isomorphism of principal bundles $\varphi_\alpha:\pi^{-1}(U_\alpha)=P_{|U_\alpha}\cong G\times U_\alpha$. Then the cocycle $\{g_{\alpha\beta}\}_{\alpha\preceq\beta}$, that determines $P$, can be identified with the map $\varphi^{-1}_\alpha\circ\varphi_\beta$ (here by abuse of notation we identify the isomorphism $\varphi_\alpha$ and its restriction to $U_{\alpha\beta}=U_\alpha\bigcap U_\beta$.)

Consider the open cover of $P$ with cylindrical sets $V_\alpha=\pi^{-1}(U_\alpha)$. These subsets form an open cover $\mathcal V$ on $P$. Then one can use the complex of \v{C}ech cochains on $P$, subdued to $\mathcal V$:
$$
\check{C}^n(\mathcal V,\,\Omega_{DR}(V))=\{h_{\alpha_0,\dots,\alpha_p}\in\Omega^q_{DR}(V_{\alpha_0,\dots,\alpha_p})|\alpha_0\prec\dots\prec\alpha_p,\ \alpha_i\in A,\
i=0,\dots,p,\ p+q=n\},
$$
where $V_{\alpha_0,\dots,\alpha_{p}}=V_{\alpha_0}\bigcap\dots\bigcap V_{\alpha_p}$. The differential on $\check{C}^n(\mathcal V,\,\Omega_{DR}(V))$ is given by the combination of de Rham differential on $V_\alpha$ and the \v{C}ech differential
$$
\delta(\{h\})_{|V_{\alpha_0,\dots,\alpha_{p+1}}}=\sum_{i=0}^{p+1} (-1)^{i}h_{\alpha_0,\dots,\widehat{\alpha_i},\dots,\alpha_{p+1}}.
$$
It is easy to show, that the cohomology of $\check{C}^n(\mathcal V,\,\Omega_{DR}(V))$ is isomorphic to the de Rham cohomology of $P$. For instance one can use induction on the number of elements in the cover $\mathcal V$ and apply Mayer-Vietoris exact sequence. Moreover, since $\Omega_{DR}^\ast(U)$ is a sheaf of algebras, there's an algebra structure on $\check{C}^n(\mathcal V,\,\Omega_{DR}(V))$ (see the formula \eqref{prod} below) and the isomorphism of cohomology, described above, commutes with the algebra structures.

Isomorphisms $\varphi_\alpha$ allow one identify the algebras of differential forms on $V_\alpha$ with tensor products $\Omega_{DR}^\ast(G)\otimes\Omega_{DR}^\ast(U_\alpha)$. We restrict the isomorphisms $\varphi_{\alpha_0}$ to $\Omega^\ast_{DR}(V_{\alpha_0,\dots,\alpha_p})$; recall, that the set $A$ is ordered and every time we consider an intersection $U_{\alpha_0}\cap\dots\cap U_{\alpha_k}$ we suppose $\alpha_0\prec\dots\prec\alpha_p$. Combining these isomorphisms we obtain an isomorphism of graded spaces
$$
\nabla:\check{C}^n(\mathcal V,\,\Omega_{DR}(V))\cong\Omega_{DR}^\ast(G)\otimes\check{C}^n(\mathcal
U,\,\Omega_{DR}(U)).
$$
However, this isomorphism is neither an isomorphism of algebras, nor does it commute with differentials. In fact, one can consider two differentials on $\Omega_{DR}^\ast(G)\otimes\check{C}^n(\mathcal U,\,\Omega_{DR}(U))$: one is induced from $\check{C}^n(\mathcal V,\,\Omega_{DR}(V))$ with the help of the isomorphism $\nabla$, and the other one is the standard differential that one always has on the tensor product of differential complexes. Let us denote the former differential by $d_P$ and the latter by $d$. Then $d=1\otimes d_U+1\otimes\delta'+d_G\otimes1$, where $d_U,\ d_G$ are the de Rham differentials on $U\subseteq X$ and $G$ respectively, and $\delta'$ is the \v{C}ech differential on $\check{C}^n(\mathcal U,\,\Omega_{DR}(U))$.

We shall suppose that the algebra of de Rham forms on $G$, $\Omega_G=\Omega_{DR}^\ast(G)$ is equiped with the differential bialgebra structure, with comultiplication induced from the product of matrices. This conjecture allows one to speak about the twisting cochains on $\Omega_G$. See example \ref{example1} of the previous section for the definition of twisting cochains. Recall that given a twisting cochain, one can introduce a twisted differential $d_\phi$ on the tensor product $\Omega_G\otimes\Omega$:
$$
d_\phi=d_G\otimes1+1\otimes d_\Omega-\phi\cap1,\
\phi\cap1(\omega\otimes\psi)=\omega^{(1)}\otimes\phi(\omega^{(2)})\psi.
$$
Our purpose is to define a twisting cochain $\phi$ on $\Omega_G$ with values in $\check{C}^n(\mathcal U,\,\Omega_{DR}(U))$ such that $d_P=d_\phi$ on $\Omega_{DR}^\ast(G)\otimes\check{C}^n(\mathcal U,\,\Omega_{DR}(U))=\Omega_G\otimes\check{C}^n(\mathcal U,\,\Omega_{DR}(U))$. To this end consider the following map (c.f. Ch. 2, \S1 of \cite{Shih}):
\begin{equation}
\label{twistcoch}
\Omega_G\stackrel{e}{\longrightarrow}\Omega_G\otimes\check{C}^n(\mathcal
U,\,\Omega_{DR}(U))\stackrel{d_p-d}{\longrightarrow}\Omega_G\otimes\check{C}^n(\mathcal
U,\,\Omega_{DR}(U))\stackrel{\epsilon\otimes1}{\longrightarrow}\check{C}^n(\mathcal
U,\,\Omega_{DR}(U)).
\end{equation}
Here $e(\omega)=\omega\otimes1$ and $\epsilon:\Omega_G\to\mathbb C$ is the counit in bialgebra $\Omega_G$. We shall denote the map \eqref{twistcoch} by $\phi_P$.
\begin{prop}
The map $\phi_P$ is a twsting cochain, i.e. the equation \eqref{flatc} holds.
\end{prop}
\begin{Proof}
Since the comultiplication on $\Omega_G$ is induced from the comultiplication on $\Omega_G^0=\ac(G)$, and the isomorphisms induced from $\varphi_\alpha$ (as well as the maps $e$ and $epsilon\otimes1$) intertwine the de Rham differentials, it is enough to check that equation \eqref{flatc} holds for $\omega\in\Omega_G^0$, i.e. for
functions on $G$.

To this end we recall that $1\in\check{C}^n(\mathcal U,\,\Omega_{DR}(U))$ is given by the degree $0$ \v{C}ech cochain that is equal to the constant function $1_\alpha\equiv1$ on every open subset $U_\alpha\in\mathcal U$. So for arbitrary $f\in\ac(G)$ we have
$$
d(e(f))=d(f\otimes\{1_\alpha\})=d_Gf\otimes\{1_\alpha\},
$$
since $d_U(\{1_\alpha\})=\delta'(\{1_\alpha\})=0$. On the other hand,
$$
d_P(e(f))=d_P(f\otimes\{1_\alpha\})=(\{\varphi^*_\alpha\})^{-1}(d_V+\delta)(\{\varphi^*_\alpha(f\otimes1_\alpha)\})
$$
But since $d_V\varphi^*_\alpha=\varphi^*_\alpha(d_G\otimes1+1\otimes d_U)$, we conclude, that
\begin{equation}
\begin{split}
d_P(e(f))-d(e(f))&=(\{\varphi^*_\alpha\})^{-1}\delta(\{\varphi^*_\alpha(f\otimes1_\alpha)\})\\
                 &=\{(\varphi^*_\alpha)^{-1}(\varphi^*_\alpha(f\otimes1_\alpha)-\varphi^*_\beta(f\otimes1_\beta))\}_{\alpha\prec\beta}\\
                 &=\{f\otimes1_\alpha-(\varphi^*_\alpha)^{-1}\varphi^*_\beta(f\otimes1_\beta)\}_{\alpha\prec\beta}\\
                 &=\{f\otimes1_\alpha-(\varphi_\beta\varphi_\alpha^{-1})^*(f\otimes1_\beta)\}_{\alpha\prec\beta}.
\end{split}
\end{equation}
Now recall, that the homeomorphism $\varphi_\beta\varphi_\alpha^{-1}:G\times U_{\alpha\beta}\to G\times U_{\alpha\beta}$ is given by $(g,\,x)\mapsto(gg_{\alpha\beta}(x),\,x)$ and $\epsilon(f)=f(e)$, where $e$ is the unit in $G$. So we conclude that
$$
\phi_P(f)=\{f(e)1_{\alpha\beta}-f\circ g_{\alpha\beta}\}_{\alpha\prec\beta}.
$$
In case when $f$ is replaced with a differential form $\omega$ of degree greater than $0$, one should write $g_{\alpha\beta}^*\omega$ instead of $f\circ g_{\alpha\beta}$.

Now we are able to check the equation \eqref{flatc}. Recall that the comultiplication in $\Omega_G$ is induced from the product of matrices, so that the following relation holds
\begin{equation}
\label{comult}
\sum\omega^{(1)}(g)\omega^{(2)}(h)=\omega(gh),\ \omega\in\Omega_G,\ g,\,h\in G.
\end{equation}
Also recall, that the multiplication in $\check{C}(\mathcal U,\,\Omega_{DR}(U))$ is given by the formula
\begin{equation}
\label{prod}
(h'\cup
h'')_{\alpha_0,\dots,\alpha_{p+q}}=(-1)^{p|h''|_2}(h'_{\alpha_0,\dots,\alpha_p})_{|U_{\alpha_0,\dots,\alpha_{p+q}}}\cdot(h''_{\alpha_p,\dots,\alpha_{p+q}})_{|U_{\alpha_0,\dots,\alpha_{p+q}}},
\end{equation}
where $|h''|_2$ is the second (de Rham) degree of $h''$. Let us denote the differential in $\check{C}(\mathcal U,\,\Omega_{DR}(U))$ by $d_B=d_U+\delta'$. We compute
\begin{equation}
\begin{split}
d_B\phi_P(f)-\phi_P(d_Gf)&=(d_U+\delta')\phi_P(f)-\phi_P(d_Gf)\\
                         &=d_U\{f(e)1_{\alpha\beta}- f\circ g_{\alpha\beta}\}_{\alpha\prec\beta}+\delta'\{f(e)1_{\alpha\beta}- f\circ g_{\alpha\beta}\}_{\alpha\prec\beta}-\{0-g_{\alpha\beta}^*d_Gf\}_{\alpha\prec\beta}\\
                         &=\delta'\{f(e)1_{\alpha\beta}- f\circ g_{\alpha\beta}\}_{\alpha\prec\beta}.
\end{split}
\end{equation}
Here we used the equation $d_Uf(g_{\alpha\beta})=g^*_{\alpha\beta}d_Gf$. Further,
\begin{equation}
\begin{split}
\delta'\{f(e)1_{\alpha\beta}&-f\circ g_{\alpha\beta}\}_{\alpha\prec\beta}=\{f(e)1_{\alpha\beta\gamma}-(f\circ g_{\beta\gamma})_{|U_{\alpha\beta\gamma}}\}_{\alpha\prec\beta\prec\gamma}\\
                            &\ \ -\{f(e)1_{\alpha\beta\gamma}-(f\circ g_{\alpha\gamma})_{|U_{\alpha\beta\gamma}}\}_{\alpha\prec\beta\prec\gamma}+ \{f(e)1_{\alpha\beta\gamma}- (f\circ g_{\alpha\beta})_{|U_{\alpha\beta\gamma}}\}_{\alpha\prec\beta\prec\gamma}\\
                            &=\{f(e)1_{\alpha\beta\gamma}- (f\circ g_{\alpha\beta})_{|U_{\alpha\beta\gamma}}+(f\circ g_{\alpha\gamma})_{|U_{\alpha\beta\gamma}}-(f\circ g_{\beta\gamma})_{|U_{\alpha\beta\gamma}}\}_{\alpha\prec\beta\prec\gamma}
\end{split}
\end{equation}
Next, we compute $\phi_P\cup\phi_P(f)$:
\begin{equation}
\begin{split}
\phi_P\cup\phi_P(f)&=\{f^{(1)}(e)1_{\alpha\beta}-f^{(1)}\circ g_{\alpha\beta}\}_{\alpha\prec\beta}\cup\{f^{(2)}(e)1_{\alpha\beta}-f^{(2)}\circ g_{\alpha\beta}\}_{\alpha\prec\beta}\\
                   &=\{(f^{(1)}(e)1_{\alpha\beta}-f^{(1)}\circ g_{\alpha\beta})(f^{(2)}(e)1_{\beta\gamma}-f^{(2)}\circ g_{\beta\gamma})\}_{\alpha\prec\beta\prec\gamma}\\
                   &=\{f^{(1)}(e)f^{(2)}(e)1_{\alpha\beta\gamma}-f^{(1)}\circ g_{\alpha\beta}f^{(2)}(e)1_{\beta\gamma}-f^{(1)}(e)1_{\alpha\beta}f^{(2)}\circ g_{\beta\gamma}\\
                   &\qquad\qquad+f^{(1)}\circ g_{\alpha\beta}f^{(2)}\circ g_{\beta\gamma}\}_{\alpha\prec\beta\prec\gamma}\\
                   &=\{f(e)1_{\alpha\beta\gamma}-(f\circ g_{\alpha\beta})_{|U_{\alpha\beta\gamma}}-(f\circ g_{\beta\gamma})_{|U_{\alpha\beta\gamma}}+(f\circ g_{\alpha\gamma})_{|U_{\alpha\beta\gamma}}\}_{\alpha\prec\beta\prec\gamma}.
\end{split}
\end{equation}
Here we have used the relations \eqref{comult} and \eqref{prod} (and also omitted the evident restrictions for the sake of brevity). Now it is evident that \eqref{flatc} holds for $\phi_P$.
\end{Proof}\\ \\ \vspace{1mm}
The next theorem is the main result of this paragraph. It is similar to that of the Th\'eorem 2, of Ch. 2 \cite{Shih}.
\begin{theorem}
The following formula holds
$$
d_P=d_{\phi_P},
$$
i.e. $d_P=d+\phi_P\cap1$, where $d$ is the usual differential on the tensor product $\Omega_G\otimes\check{C}^\ast(\mathcal U,\,\Omega_{DR}(U))$.
\end{theorem}
\begin{Proof}
The simplest way to demonstrate this fact is a direct inspection of formulas. We compute
\begin{equation}
\label{eqdphip}
\begin{split}
(d+&\phi_P\cap1)(\omega\otimes\{h_{\alpha_1\dots\alpha_n}\})=d_G\omega\otimes\{h_{\alpha_1\dots\alpha_n}\}+(-1)^{|\omega|}\omega\otimes\{d_Uh_{\alpha_1\dots\alpha_n}\}\\
   &+(-1)^{|\omega|}\omega\otimes\delta_U\{h_{\alpha_1\dots\alpha_n}\}-(-1)^{|\omega|}\omega^{(1)}\otimes\{\omega^{(2)}(e)1_{\alpha\beta}-g_{\alpha\beta}^*\omega^{(2)}\}\cup\{h_{\alpha_1\dots\alpha_n}\}
\end{split}
\end{equation}
Observe further, that $\delta'\{h_{\alpha_1\dots\alpha_n}\}=\{1_{\alpha\beta}\}\cup\{h_{\alpha_1\dots\alpha_n}\}$
and $\omega^{(1)}\cdot\omega^{(2)}(e)=\omega$, while $\omega^{(1)}\cdot g_{\alpha\beta}^*\omega^{(2)}=(\varphi_\beta\varphi_\alpha^{-1})^*\omega$. Thus the expression on the right side of \eqref{eqdphip} is equal to
$$
d_G\omega\otimes\{h_{\alpha_1\dots\alpha_n}\}+(-1)^{|\omega|}\omega\otimes\{d_Uh_{\alpha_1\dots\alpha_n}\}+(-1)^{|\omega|}(\varphi_\beta\varphi^{-1}_\alpha)^*\omega\otimes\{h_{\alpha_1\dots\alpha_n}\}
$$
which is precisely equal to $d_P$.
\end{Proof}\\ 

\vspace{-2mm}
\begin{rem}\rm
It is not difficult to show that if $g_{\alpha\beta}$ and $g'_{\alpha\beta}$ are two cohomologous noncommutative $1$\/-cocycles with values in the group $G$, then the twisting cochains, determined by these cocycles are gauge-equivalent. In fact recall that two cocycles $g_{\alpha\beta}$ and $g'_{\alpha\beta}$ are cohomologous, iff there exists a collection of maps $h_\alpha:U_\alpha\to G$ such that $(h_\alpha)_{|_{U_{\alpha\beta}}} g_{\alpha\beta}=g'_{\alpha\beta}(h_{\beta})_{|_{U_{\alpha\beta}}}$. Then it is easy to check by a direct calculation that the following formula defines a gauge transformation $c_h$ connecting the twisting cochains, associated with $g_{\alpha\beta}$ and $g'_{\alpha\beta}$:
$$
\Omega_G\ni k\stackrel{c_h}{\mapsto}\{h_{\alpha}^*(k)\}.
$$
\hfill$[]$\end{rem}
\begin{rem}\rm
\label{leftright}
It is now easy to prove the fact we used in proposition \ref{propleftright}: that the left twisted tensor product $\Omega{}_{\phi}\otimes K$ models the fibre bundle $\hat P=P\otimes_{R^{-1}}G$. To this end it is enough to prove this only for one particular choice of twisting cochain in the equivalence class. Similarly, one can manipulate with coalgebra $K$ and algebra $\Omega$, see section \ref{sectinfty}. So let us take $\phi=\phi_P$\ -- the cochain, we constructed in this section, $K$ and $\Omega$ will be equal to $\DR{G}$ and  $\check{C}^\ast(\mathcal U,\,\Omega_{DR}(U))$ respectively. A direct computation now shows that in this case differential in $\Omega{}_\phi\otimes K$ coincides with the differential induced on $\check{C}^\ast(\mathcal U,\,\Omega_{DR}(U))\otimes\DR{G}$ by the construction, described above, from the following principal bundle: 
$$
\hat{\hat P}=\coprod_\alpha U_\alpha\times G/\sim.
$$
Here the equivalence relation $\sim$ is given by $U_{b}\ni(x,\,g)\sim(x,\,g_{\alpha\beta}g)\in U_\alpha$ for $x\in U_{\alpha\beta}$. Comparing this with the analoguous construction of $P$, where $(g,\,x)\in U_\alpha$ is equivalent to $(gg_{\alpha\beta},\,x)\in U_\beta$ we conclude that $\hat{\hat P}=\hat P$.
\hfill$[]$\end{rem}

\subsection{Twisting cochains and Chern-Weil classes}
\label{sect32}
Now we are going to discuss the relation of the twisting cochain defined above and the Chern-Weil construction of characteristic classes of principal bundles. A more detailed overview of this subject one can find in the paper \cite{mypapLF}, in which we also compare the construction, presented here with previously-known ones (that of Bott and Dupont, see \cite{BottTu,Dupont}.)

So consider a general twisting cochain $\phi$ on an coaugmented differential graded coalgebra $K$ with values in a differential graded algebra $\Omega$. Recall that coaugmentation is a map $\eta:\Bbbk\to K$, where $\Bbbk$ is regarded as the trivial $1$\/-dimensional differential coalgebra, generated by a group-like element $1$ vanishing under differential. We shall identify $1\in\Bbbk$ and its image $\eta(1)\in K$. Thus coaugmenting $K$ is equivalent to choosing a closed group-like element $1\in K$. Recall that we assume that a twisting cochain $\phi$ sends $1$ to $0$.

One calls an element $\omega$ in a coaugmented coalgebra $K$ primitive, if $\Delta(\omega)=\omega\otimes1+1\otimes\omega$, where $\Delta$ stands for the comultiplication in $K$. Let $Pr^*(K)$ denote the space of primitive elements in $K$ (with grading, induced from $K$.) Clearly, $Pr^*(K)$ is a subcomplex in $K$. From the definition of twisting cochain it now follows that $\phi$ defines a degree $1$ map of chain complexes $\phi:Pr^*(K)\to A^{*+1}$. We shall call this map {\em transgression, determined by $\phi$\/}.

It follows from the discussion of the gauge transformations in the end of the previous section, that the maps from the space $Pr^*K$ to $A$, induced by $\phi$ and $\phi\circ c$ (see the definition of the gauge transformations $c$ and of the action $\circ$ in the remark \ref{remgaugetw}) are homotopic, in particular they induce the same map on cohomology. More generally, the twisting cochain $\phi$ defines a homomorphism of differential graded algebras $\tilde\phi$ (see section \ref{sectnext2}) from the cobar-resolution $F(K)$ of $K$ to $A$. The chain homotopy class of this map does not change with the action of gauge transformations group $\hc(K,A)$, see section \ref{sectnnext}. The map we considered above is equal to the composition of $\phi^\otimes$ and the natural inclusion of $Pr^*(K)$ to $F(K)$ (which commutes with differential, if we use the reduced version of cobar-construction, i.e. factor out the entries with $1$ in the tensor powers of $K$.)

From our previous discussion (see the cited sections and E.~Brown's paper \cite{Br59}) it follows that the map $\tilde\phi$ is a cohomological counterpart of the classifying map $f:X\to BG$ of the bundle $P$ in the sense that $F(K)$ is a cochain complex, whose cohomology coincides with $H^*(BG)$ (here we consider cohomology with coefficients in the characteristic zero field $\Bbbk$.) On the other hand, characteristic classes of $P$ are given by the pullbacks $f^*(c)$, where $c$ is a class in $H^*(BG)$. Thus we see that the problem of finding the characteristic class, associated with an element $c\in H^*(BG)$ by the map $f^*$ is reduced to finding a closed cocycle in $F(K)$ which would correspond to $c$ under the isomorphism of cohomology.

Putting for a while aside the problem of finding such a representing cocycle in a general situation, we shall give here a simple example of such a formula. Let $K=\Omega_G$ for an algebraic matrix group $G$, and let $\omega=g^{-1}dg$ be the canonical left-invariant $\g$\/-valued Maurer-Cartan form on $G$, and $\tilde\omega=(dg)g^{-1}$ -- its right invariant counterpart. Here $g$ is the generic element on $G$. One should regard $g$ as a matrix-valued function on $G$, which associates to $x\in G$ its matrix form with respect to a fixed basis. If $u_{ij}$ are the generating functions of the algebra $\ac(G)$, one can identify $g$ with the element $(u_{ij})$ of $Mat_{n}\,(\ac(G))$. Then both $\omega$ and $\tilde\omega$ are matrix-valued differential $1$\/-forms and we can consider the element $\omega_3=Tr(\omega\wedge\omega\wedge\omega)=\omega_i^k\wedge\omega_k^l\wedge\omega_l^i$ (the summation is taken over the repeating indices). This element is closed, but not primitive. It is easy to check this with the help of formulas (the first formula is just the Maurer-Cartan equation, and the second one follows directly from the definition of comultiplication on the level of $u_{ij}$)
\begin{equation}
\label{MauCar}
d\omega_i^j=-\omega_i^k\wedge\omega_k^j
\end{equation}
and
\begin{equation}
\label{primit}
\Delta(\omega_i^j)=\omega_k^l\otimes g_i^kg_l^j+1\otimes\omega_i^j.
\end{equation}
However, in spite of its being not primitive, this differential form ($\omega_3$) is biinvariant with respect to the action of $G$. As one knows (see e.g. \cite{Onisch}) such differential forms generate the cohomology of $G$. Moreover, one can show that on the level of cohomology, the corresponding element is primitive with respect to the Pontrjagin coproduct, thus it can be transgressed to a cocycle in $BG$. However, $\omega_3$ being not primitive with respect to the usual comultiplication in $\Omega_G$, it doesn't give a cocycle in $F(K)$, if we embed $\Omega_G=K$ into $F(K)$ as the space of $1$\/-tensors. However one can modify $\omega_3$ a little bit inside $F(K)$ by adding a correction term which will have tensor degree 2 so that the result will be closed: take $\hat\omega_3\in F(K)$ equal to $\omega_3-3Tr(\omega\otimes\tilde\omega)$, where
$$
Tr(\omega\otimes\tilde\omega)=\omega_i^k\otimes\tilde\omega_k^i,
$$
the summation is again taken over the repeating indices. Once again formulas \eqref{MauCar} and \eqref{primit} (and similar formulas for the form $\tilde\omega$) are used to check that $\hat\omega_3$ is a cocycle in $F(K)$. If we apply the map $\phi^\otimes$ (where $\phi=\phi_P$ is the twisting cochain from the previous paragraph) to this element, we obtain the \v{C}ech form
$$
\{Tr(g_{\alpha\beta}^{-1}dg_{\alpha\beta}\wedge g_{\alpha\beta}^{-1}dg_{\alpha\beta}\wedge
g_{\alpha\beta}^{-1}dg_{\alpha\beta})\}_{\alpha\prec\beta}-3\{Tr(g_{\alpha\beta}^{-1}dg_{\alpha\beta}\wedge
dg_{\beta\gamma}g_{\beta\gamma}^{-1})\}_{\alpha\prec\beta\prec\gamma}.
$$
It is easy to check that this form is closed and corresponds to the second Chern class of the vector bundle, associated to our principal bundle. In effect, one first can check that this cochain is closed. To this end we compute:
$$
d(g_{\alpha\beta}^{-1}dg_{\alpha\beta})=dg_{\alpha\beta}^{-1}\wedge dg_{\alpha\beta}=dg_{\alpha\beta}^{-1}g_{\alpha\beta}\wedge g_{\alpha\beta}^{-1}dg_{\alpha\beta}=-g_{\alpha\beta}^{-1}dg_{\alpha\beta}\wedge g_{\alpha\beta}^{-1}dg_{\alpha\beta},
$$
and similarly $d(dg_{\alpha\beta}g_{\alpha\beta}^{-1})=dg_{\alpha\beta}g_{\alpha\beta}^{-1}\wedge dg_{\alpha\beta}g_{\alpha\beta}^{-1}$.
Hence
\begin{equation}
\label{dtrace}
\begin{split}
d\bigl(\{Tr(g_{\alpha\beta}^{-1} &dg_{\alpha\beta}\wedge g_{\alpha\beta}^{-1}dg_{\alpha\beta}\wedge
g_{\alpha\beta}^{-1}dg_{\alpha\beta})\}_{\alpha\prec\beta}+3\{Tr(g_{\alpha\beta}^{-1}dg_{\alpha\beta}\wedge
dg_{\beta\gamma}g_{\beta\gamma}^{-1})\}_{\alpha\prec\beta\prec\gamma}\bigr)\\
                                                      &=-\{Tr(g_{\alpha\beta}^{-1}dg_{\alpha\beta}\wedge g_{\alpha\beta}^{-1}dg_{\alpha\beta}\wedge
g_{\alpha\beta}^{-1}dg_{\alpha\beta}\wedge g_{\alpha\beta}^{-1}dg_{\alpha\beta})\}_{\alpha\prec\beta}\\
                                                      &\quad+3\bigl(\{Tr(g_{\alpha\beta}^{-1}dg_{\alpha\beta}\wedge g_{\alpha\beta}^{-1}dg_{\alpha\beta}\wedge
dg_{\beta\gamma}g_{\beta\gamma}^{-1})\}_{\alpha\prec\beta\prec\gamma}\\
                                                      &\qquad\qquad\qquad+\{Tr(g_{\alpha\beta}^{-1}dg_{\alpha\beta}\wedge
dg_{\beta\gamma}g_{\beta\gamma}^{-1}\wedge dg_{\beta\gamma}g_{\beta\gamma}^{-1})\}_{\alpha\prec\beta\prec\gamma}\bigr)\\
                                                      &=3\bigl(\{Tr(g_{\alpha\beta}^{-1}dg_{\alpha\beta}\wedge g_{\alpha\beta}^{-1}dg_{\alpha\beta}\wedge
dg_{\beta\gamma}g_{\beta\gamma}^{-1})\}_{\alpha\prec\beta\prec\gamma}\\
                                                      &\qquad\qquad\qquad+\{Tr(g_{\alpha\beta}^{-1}dg_{\alpha\beta}\wedge
dg_{\beta\gamma}g_{\beta\gamma}^{-1}\wedge dg_{\beta\gamma}g_{\beta\gamma}^{-1})\}_{\alpha\prec\beta\prec\gamma}\bigr).
\end{split}
\end{equation}
Here we used the equality $Tr(\omega^{\wedge 4})=0$, which is due to the graded cyclical invariance of the trace.
Now from the cocycle relation $g_{\alpha\gamma}=g_{\alpha\beta}g_{\beta\gamma}$ for $\alpha\prec\beta\prec\gamma$ such that the open maps $U_\alpha,\ U_\beta$ and $U_\gamma$ intersect, we obtain the formula
$$
g_{\alpha\gamma}^{-1}dg_{\alpha\gamma}=(g_{\alpha\beta}g_{\beta\gamma})^{-1}d(g_{\alpha\beta}g_{\beta\gamma})=g_{\beta\gamma}^{-1}(g_{\alpha\beta}^{-1}dg_{\alpha\beta})g_{\beta\gamma}+g_{\beta\gamma}^{-1}dg_{\beta\gamma}
$$
and similarly $dg_{\alpha\gamma}g_{\alpha\gamma}^{-1}=dg_{\alpha\beta}g_{\alpha\beta}^{-1}+g_{\alpha\beta}(dg_{\beta\gamma}g_{\beta\gamma}^{-1})g_{\alpha\beta}^{-1}$. Applying the first equality to $Tr\bigl((g_{\alpha\gamma}^{-1}dg_{\alpha\gamma})^{\wedge 3}\bigr)$, and using the (graded) cyclical invariance of trace, we get
\begin{equation*}
\begin{split}
Tr\bigl((g_{\alpha\gamma}^{-1}dg_{\alpha\gamma})^{\wedge 3}\bigr)&=Tr\bigl((g_{\alpha\beta}^{-1}dg_{\alpha\beta})^{\wedge 3}\bigr)+3Tr\bigl((g_{\alpha\beta}^{-1}dg_{\alpha\beta})^{\wedge 2}\wedge dg_{\beta\gamma}g_{\beta\gamma}^{-1}\bigr)\\                                                                              &\quad+3Tr\bigl(g_{\alpha\beta}^{-1}dg_{\alpha\beta}\wedge (dg_{\beta\gamma}g_{\beta\gamma}^{-1})^{\wedge 2}\bigr)+Tr\bigl((g_{\alpha\gamma}^{-1}dg_{\alpha\gamma})^{\wedge 3}\bigr).
\end{split}
\end{equation*}
Applying this formula, we obtain
\begin{equation}
\label{deltatrace}
\begin{split}
\delta(\{Tr\bigl((g_{\alpha\beta}^{-1}dg_{\alpha\beta})^{\wedge 3}\bigr)\})_{\alpha\prec\beta\prec\gamma}&=\Bigl\{Tr\bigl((g_{\alpha\beta}^{-1}dg_{\alpha\beta})^{\wedge 3}\bigr)-Tr\bigl((g_{\alpha\gamma}^{-1}dg_{\alpha\gamma})^{\wedge 3}\bigr)+Tr\bigl((g_{\beta\gamma}^{-1}dg_{\beta\gamma})^{\wedge 3}\bigr)\Bigr\}\\
                &=\Bigl\{-3Tr\bigl((g_{\alpha\beta}^{-1}dg_{\alpha\beta})^{\wedge 2}\wedge dg_{\beta\gamma}g_{\beta\gamma}^{-1}\bigr)-3Tr\bigl(g_{\alpha\beta}^{-1}dg_{\alpha\beta}\wedge (dg_{\beta\gamma}g_{\beta\gamma}^{-1})^{\wedge 2}\bigr)\Bigr\}
\end{split}
\end{equation}
(here in order to make the notation more readable we do not distinguish differential forms on $U_\alpha\bigcap U_\beta$ and their restrictions to the triple intersections of open sets.)

Similarly, we can apply the cocyclce relations and the cyclicity of trace to $Tr(g_{\alpha\gamma}^{-1}dg_{\alpha\gamma}\wedge dg_{\gamma\delta}g_{\gamma\delta}^{-1})$ 
and obtain the formula
\begin{equation*}
Tr(g_{\alpha\gamma}^{-1}dg_{\alpha\gamma}\wedge dg_{\gamma\delta}g_{\gamma\delta}^{-1})=Tr(g_{\alpha\beta}^{-1}dg_{\alpha\beta}\wedge dg_{\beta\delta}g_{\beta\delta}^{-1})-Tr(g_{\alpha\beta}^{-1}dg_{\alpha\beta}\wedge dg_{\beta\gamma}g_{\beta\gamma}^{-1})+Tr(g_{\beta\gamma}^{-1}dg_{\beta\gamma}\wedge dg_{\gamma\delta}g_{\gamma\delta}^{-1}),
\end{equation*}
which is clearly equivalent to the equation $\delta(\{Tr(g_{\alpha\beta}^{-1}dg_{\alpha\beta}\wedge dg_{\beta\gamma}g_{\beta\gamma}^{-1})\}_{\alpha\prec\beta\prec\gamma})=0$. Now, summing up the equations \eqref{dtrace} and \eqref{deltatrace}, we obtain the necesary equality:
$$
(d+\delta)(\{Tr(g_{\alpha\beta}^{-1}dg_{\alpha\beta}\wedge g_{\alpha\beta}^{-1}dg_{\alpha\beta}\wedge
g_{\alpha\beta}^{-1}dg_{\alpha\beta})\}_{\alpha\prec\beta}-3\{Tr(g_{\alpha\beta}^{-1}dg_{\alpha\beta}\wedge
dg_{\beta\gamma}g_{\beta\gamma}^{-1})\}_{\alpha\prec\beta\prec\gamma})=0.
$$
Let us denote this cocycle by $\tilde c_2$. There exist many more or less abstract ways to show that this element corresponds to a characteristic class of the principal bundle $P$. But instead of this we shall just specify a cochain $a=a(g_{\alpha\beta}, A, F)$, where $A,\ F$ are a connection in the principal bundle $P$ and its curvature, such that
\begin{equation}
\label{chersim}
(d+\delta)(a)=3\{Tr(F\wedge F)\}_\alpha-\tilde c_2.
\end{equation}
Namely, it is well-known that $Tr(F^{\wedge 2})=dTr(F\wedge A-\frac13 A^{\wedge 3})$ (the Chern-Simons form on the right is not invariant under the gauge transformations in general.) So, we put
$$
a(g_{\alpha\beta},A,F)=\{Tr(3F_{\alpha}\wedge A_{\alpha}- A_{\alpha}^{\wedge 3})\}_{\alpha}+b(g_{\alpha\beta},A,F),
$$
where $A_{\alpha}$ and $F_\alpha$ are the local gauge potential and curvature respectively (i.e. gauge potential in the open set $U_\alpha$) and $b$ is a \v Cech 1-cochain. An easy calculation, based on the Bianchi identities
\begin{align*}
dA_\alpha&=F_\alpha+A_\alpha\wedge A_\alpha,\\
dF_\alpha&=F_\alpha\wedge A_\alpha-A_\alpha\wedge F_\alpha
\end{align*}
and the gauge transformation formulas
\begin{align}
\label{transform3}
A_\beta&=g_{\alpha\beta}^{-1}A_\alpha g_{\alpha\beta}+ g_{\alpha\beta}^{-1}dg_{\alpha\beta},\\
\label{transform4}
F_\beta&=g_{\alpha\beta}^{-1}F_\alpha g_{\alpha\beta}
\end{align}
gives the formula
\begin{equation*}
\begin{split}
(d+\delta)(\{Tr(3F_{\alpha}\wedge A_{\alpha}- A_{\alpha}^{\wedge 3})\}_{\alpha})&=3\{Tr(F_\alpha\wedge F_\alpha)\}_\alpha-\{Tr\bigl((g_{\alpha\beta}^{-1}dg_{\alpha\beta})^{\wedge 3}\bigr)\}_{\alpha\prec\beta}\\                                                                                       &\quad+3(-\{dA_\alpha\wedge dg_{\alpha\beta}g_{\alpha\beta}+A_{\alpha}\wedge(dg_{\alpha\beta}g_{\alpha\beta})^{\wedge 2}\}_{\alpha\prec\beta})
\end{split}
\end{equation*}
The third term on the right hand side is a \v{C}ech cochain, which is equal to exact forms on all open sets where it is defined, so we can guess the formula for $b$: put $b=-\{A_\alpha\wedge dg_{\alpha\beta}g_{\alpha\beta}^{-1}\}_{\alpha\prec\beta}$. Now one obtains the formula \eqref{chersim} by a direct calculation.

In fact, it is possible to show in general that there always exists a cochain, relating the forms we have described here (coming from the twisted cochain described in this section, applied to the cobar resolution) and the usual Chern-Weil forms. To this end, let us briefly recall the classical Chern-Weil theory (see e.g. \cite{Onisch} and \cite{mypapLF} for details of the discussion that follows.)

To every (graded) Lie algebra \g\ one can associate a so-called \emph{completed} Weil algebra $W(\g)=\widehat S^*(\g^*[2])\hat\otimes\Lambda^*(\g^*[1])$, where $\widehat S^*$ and $\Lambda^*$ denote the \emph{completed} symmetric and the exterior algebras of the vector space $\g^*$ (the dual space of \g) respectively, and $[1]$ and $[2]$ denote the shift of dimension, i.e. $\g^*[1]^k=(\g^*)^{k-1}$ and $\g^*[2]^k=(\g^*)^{k-2}$ for all $k\in\mathbb Z$. The adjective ``completed'' means that one can consider infinite sums of elements in the symmetric (or exterior) algebra, provided their tensor powers grow so that there remain only a finite sum in every finite tensor product. The elements of $\g^*[1]$ are denoted by $a_{X^*}$ and the elements of $\g^*[2]$ -- by $f_{X^*}$, where $X^*$ are the corresponding elements of $\g^*$. The differential $\partial$ in $W(\g)$ is defined in the following way: let $a_i$ stand for $a_{X^*_i}$ and $f_i=f_{X^*_i}$, where $X^*_i,\ i=1,\dots, n$ is the basis of $\g^*$. Then one puts:
\begin{align}
\partial a_i&=f_i+ \frac12C_i^{jk}a_ja_k,\\
\partial f_i&=C_i^{jk}f_ja_k,
\end{align}
where $C_i^{jk}$ are the structure constatnts of the Lie algebra $\g$ (i.e. if $X^i$ is the dual to $X_i^*$ basis of \g, then $[X^j,\,X^k]=C^{jk}_iX_i$.) One can show that $W(\g)$ is the universal commutative differential algebra  in the following sense: let $\Omega^*$ be an arbitrary commutative differential algebra (which either has finite grading, or is completed in the graded sense) and $v:\g^*\to\Omega^1$ -- an arbitrary linear map. Then there exists a unique homomorphism of differential algebras $v^*:W(\g)\to\Omega^*$, such that $v^*_{|_{\g^*[1]}}=v$.

The Lie algebra \g\ couples with $W(\g)$ in two different ways: first, by adjoint action on $\g^*$ (i.e. $\theta_X(a_{Y^*})=a_{ad^*_X(Y^*)}$ and similarly for $f_{Y^*}$) and second by inner multiplications (i.e. $\imath_X(a_{Y^*})=Y^*(X)$ and $\imath_X(f_Y^*)=0$.) One easily checks the usual Cartan's identities for these actions, e.g. $[\partial,\,\imath_X]=\theta_X$, etc.. One calls an element $\omega\in W(\g)$ basic, if $\theta_X(\omega)=\imath_X(\omega)=0$ for all $X\in\g$. Clealy, the set of all basic elements in $W(\g)$ is a differential subalgebra in $W(\g)$. It is a well-known fact, that if \g\ is the Lie algebra of a compact Lie group $G$, the cohomology of the basic subalgebra $W(\g)^{bas}$ of its Weil algebra is isomorphic to the real-valued cohomology of the classifying space $BG$ of $G$. On the other hand, the cohomology of the algebra $W(\g)$ is trivial.

Let us now inetgrate $ad^*$ to the coadjoint action $Ad^*$ of the Lie group $G$ on $\g^*$. We shall assume that $G$ is an algebraic group, more precisely that it is a subgroup of the general linear group. Then this action of $G$ on $\g^*$ is algebraic in the following sense: there is a homomorphism of algebras $\Delta_W:\g^*\to \mathcal \ac(G)\otimes\g^*$, such that
\begin{equation}
\label{transform1}
g(X^*)=\Delta_W(X^*)(g).
\end{equation}
Here we regard $\Delta_W(X^*)$ as a $\g^*$\/-valued function on $G$. We extend this coaction to a map $\tilde\Delta_W:\g^*\to\ac(G)\otimes\g^*\oplus \Omega_G^1\otimes 1\subset (\Omega_G\otimes W(\g))^1$ by the formula
\begin{equation}
\label{transform2}
\tilde\Delta_W(X^*)=\Delta_W(X^*)+\omega_{X^*}\otimes 1,
\end{equation}
where $\omega_{X^*}$ is the unique right-invariant differential $1$\/-form on $G$, which coincides with $X^*$ on the tangent space in the unit element of $G$. Now, since $W(\g)$ is the universal (in the above sense) commutative differential algebra, there exists a unique extension $\tilde\Delta_W^*:W(\g)\to\Omega^*_{DR}(G)\otimes W(\g)$ of $\tilde\Delta_W$ to a homomorphism of differential algebras. It is easy to check, that the restriction of $\tilde\Delta_W^*$ to $W(\g)^{bas}$ is given by the formula $\tilde\Delta_W^*(\omega)=1\otimes\omega$, moreover, an element $\omega$ of $W(\g)$ is basic if and only if $\tilde\Delta_W^*(\omega)=1\otimes\omega$.

Thus $W(\g)$ is a (left) $\Omega_{DR}^*(G)$ differential comodule. We can now form the cobar resolution of $W(\g)$ as $\DR{G}$\/-comodule,
$$
F(W)=F(W(\g),\,\DR{G})=\bigoplus_{n\ge0}\DR{G}^{\otimes n}\otimes W(\g).
$$
The differential in $F(W)$ is equal to the sum of the differentials $\partial$ and $d$ in $W(\g)$ and $\DR{G}$ and the alternating sum of the comultiplications of the tensor factors. Observe that both $F(\DR{G})$ and $W(\g)^{bas}$ can be embedded into $F(W)$ as its cochain subcomplexes, namely: in the former case we use the map
\begin{align}
\label{incl1}
i_1:[a_1|\dots]a_n]&\mapsto1[a_1|\dots|a_n],\\
\intertext{and in the latter case the map}
\label{incl2}
i_2:\omega&\mapsto\omega[\,].
\end{align}
 The following proposition is an important technical step to the proof of the claim above:
\begin{prop}
\label{prop33}
The cohomology of $F(W)$ is equal to the real-valued cohomology of $BG$. Moreover, the inclusions \eqref{incl1} and \eqref{incl2} give isomorphisms on cohomologies.
\end{prop}
\begin{Proof} This statement seems to be generally known but since we failed to find its proof in literature, we give some hint to it here, see \cite{mypapLF} for further details.

First we consider the following filtration of $F(W)$:
$$
F'_p(F(W))=\bigl\{[x_1|\dots|x_n]\omega\in F(W)\ {\Bigl|}\ x_i\in\DR{G},\ \omega\in W(\g), |\omega|\le p\bigr\}
$$
for $p\ge0$, i.e. this is just the filtration, associated to the degree of $F(W)$ ``in the $W(\g)$ direction.'' Then, the associated spectral sequence can be easily computed: $E'_0=grF'_p(F(W))=T(\DR{G})\otimes W(\g)$ as linear space (here $T(V)$ is the tensor algebra of a vector space $V$) and the differential is equal to $1\otimes\partial$. But since $W(\g)$ is acyclic in degrees greater than $0$, we conclude that ${E'}_1^{p0}=T(\DR{G})^p$, and $E_1^{pq}=0,\ q\ge1$ and the differential in $T(\DR{G})$ is equal to the usual differential in the cobar construction, where we identify $T(\DR{G})$ with $F(\DR{G})$ in an obvious way (the differentials coincide because the unique element of degree $0$ in $W(\g)$ is $1$.) Hence, the spectral sequence collapces at $E_2$\/-term, and its cohomology is equal to the cohomology of $F(\DR{G})$, moreover the inclusion \eqref{incl1} establishes an isomorphism in cohomology.

Now consider another filtration $F''$ on $F(W)$, defined as follows
$$
F''_p(F(W))=\bigl\{[x_1|\dots|x_n]\omega\in F(W)\ {\Bigl|}\ n\le p\bigr\}.
$$
The $E_1$\/-term of the associated spectral sequence can be identified with the cohomology of $W(\g)$ as the (left) $\DR{G}$\/-comodule, non-differential.

In order to calculate this cohomology, consider the exponential map $\exp:\g\to G$. As it is well-known, if we identify \g\ with the space of left-invariant vector fields on $G$, this map will intertwine the (right) action of $G$ on itself by translations and its adjoint action on \g. Then the inverse image of $\exp$ determines a homomorphism from $\ac(G)$ (the algebra of polynomial or power series functions on $G$) to the algebra of formal power series functions on \g, compatible with the action of $G$. But this latter algebra is naturally isomorphic to $\widehat S^*(\g^*)$. So, ignoring the grading, we obtain a homomorphism $\ac(G)\to \widehat S^*(\g^*[2])$, compatible with the coaction of $\ac(G)$ and hence, as it is easy to check, with coaction of $\DR{G}$. If we combine this map with the evident identification of the differential forms on $G$, regarded as a $G$\/-module, with $\ac(G)\otimes\Lambda[\omega_{X^*_1},\dots,\omega_{X^*_n}]\cong\ac(G)\otimes\Lambda^*(\g^*[1])$ (which is due to the fact that every group is parallelisable,) we obtain a homomorphism $\DR{G}\to W(\g)$, compatible with coaction of $\DR{G}$ on both sides. Hence, $W(\g)$ turns into a \emph{Hopf-module} over the Hopf algebra $\DR{G}$ (see \cite{Sweedler} for the definition and properties of Hopf-modules.)

It follows from a well-known property of Hopf-modules, that $W(\g)\cong\DR{G}\otimes W(\g)^{bas}$ as $\DR{G}$\/-comodules. Hence, the cohomology of $W(\g)$ as of $\DR{G}$\/-comodule is equal to $W(\g)^{bas}$ in tensor degree $0$ and $0$ otherwise, recall that here we neglect the grading in $W(\g)$. So ${E''}_1^{p0}=(W(\g)^{bas})^p,\ {E''}_1^{pq}=0,\ q\ge1$. Hence this spectral sequence also collapses at $E_2$\/-term and converges to the cohomology of $W(\g)^{bas}$. Moreover, the inclusion of \eqref{incl2} establishes the isomorphism in cohomology. Observe that we have also proved that the cohomology of cobar resolution of the differential coalgebra $\DR{G}$ is isomorphic to the cohomology of $W(\g)^{bas}$, which is certainely true, since, as we have mentioned, both complexes are models for the cohomology of classifying space of the group $G$.
\end{Proof}\\ 

Now we can define a map  $c_\phi:F(W)\to\check{C}(\mathcal U,\,\Omega_{DR}(U))$, combining the Chern-Weil map on $W(\g)$ and the map, determined by the twisting cochain on the cobar resolution. Namely, put
\begin{equation}
\label{eqqqq}
c_\phi([x_1|\dots|x_n]\omega)=\phi(x_1)\cup\dots\cup\phi(x_n)\cup \{CW_\alpha(\omega)\}_\alpha,
\end{equation}
where $CW_\alpha:W(\g)\to\Omega^*_{DR}(U_\alpha)$ is the map, determined by the universal property of $W(\g)$ and the linear map
\begin{equation}
\label{transform5}
X^*\mapsto X^*(A_\alpha)\in\Omega_{DR}(U_\alpha),
\end{equation}
where $A_\alpha:U_\alpha\to\g$ are the local gauge potentials of a connection. So $\{CW_\alpha(\omega)\}_\alpha$ is a Chech $0$\/-cochain with values in $p$\/-forms on the open subsets, $p$ is the degree of $\omega$. It is not difficult to show that $c_\phi$ is a chain map. Indeed it is only necessary to check that $c_\phi|_{W(\g)}$ commutes with the differentials. But $CW_\alpha$ is by definition a homomorphism of differential graded algebras, so $d_{U_\alpha}CW_\alpha(\omega)=CW_\alpha(\partial\omega)$. On the other hand $\delta(\{CW_\alpha(\omega)\}_\alpha)=\{CW_\alpha(\omega)-CW_\beta(\omega)\}_{\alpha\prec\beta}$. But from the formulas \eqref{transform1}, \eqref{transform2}, the Bianchi identity \eqref{transform3} and the definition of $CW_\alpha$ (formula \eqref{transform5}) it follows that
$$
CW_\alpha(\omega^{(1)})\wedge g_{\alpha\beta}^*(\omega^{(2)})=CW_\beta(\omega),
$$
here $\tilde\Delta^*_W(\omega)=\omega^{(1)}\otimes\omega^{(2)}$. In effect, this formula is true for the generators $a_{X^*}$, one checks this by direct computation, and hence for all elements of $W(\g)$, since both sides of this identity are homomorphisms of d.g. algebras. So $\delta(\{CW_\alpha(\omega)\}_\alpha)=c_\phi(\tilde\Delta^*_W(\omega))$.

Thus $c_\phi:F(W)\to\check{C}(\mathcal U,\,\Omega_{DR}(U))$ is a chain map. It is easy to see that the composition of $c_\phi$ with \eqref{incl1} coincides with $\tilde\phi$ and the composition of $c_\phi$ with \eqref{incl2} coincides with the (localized) Cher-Weil homomorphism. Now, since the inclusions \eqref{incl1} and \eqref{incl2} are homotopic to each other, we conclude that for any basic element in Weil complex $x\in W(\g)^{bas}$, there exists a corresponding element $y\in F(\DR{G}$, which gives the same class in $F(W)$. Thus there must exist an element $z\in F(W)$, such that $dz=i_1(x)-i_2(y)$. Applying $c_\phi$ to both sides of this equality we conclude that for every closed form $\xi\in\Omega_{DR}^p(X)$ in the image of the Chern-Weil homomorphism, we can find a corresponding cocycle $\eta\in\sum_{i+j=p}\check{C}^i(\mathcal U,\,\Omega^j_{DR}(U))$ in the image of $\tilde\phi$ and a cochain $\zeta\in\sum_{i+j=p-1}\check{C}^i(\mathcal U,\,\Omega^j_{DR}(U))$ such that
$$
\xi-\eta=(d+\delta)\zeta.
$$
\begin{rem}
In fact, one can use proposition \ref{prop33} to find explicit formulas for the cocycles $\alpha\in F(\DR{G})$ cohomologuous to $Tr(F^{n})\in W(\g)^{bas}$, where $F$ is the ``universal curvature matrix'' in $W(\g)$. Hence, applying $\tilde\phi$ to these cocycless, one obtains formulas, expressing the Chern classes in the terms of the cocycle $g_{\alpha\beta}$. Namely, there exist an explicit way to write down a contracting homotopy for $W(\g)$. Applying this homotopy to $Tr(F^{n})$, we then use the second differential in $F(W)$. Subtract the part with the first filtration ($F'$, see the proof of \ref{prop33}) equal to $0$ from the result (these elements will belong to $F(\DR{G})$) and repeat the procedure. Each time you should add the elements with $0$ first filtration to the previously subtracted. After finitely many iterations we obtain $0$. Then the sum of all subtracted elements will correspond to the element of $F(\DR{G})$ we are looking for. This procedure can be rendered completely algorithmic with the help of the perturbation lemma, see equations \eqref{perturb}, \eqref{perteqn} and discussion that follows in section \ref{sectinfty} below. A reader, interested in details and concrete formulas can refer to \cite{mypapLF}. Besides this, there are other approaches that allow one find explicit formulas for the classes in cobar resolution and express the characteristic classes of a principal bundle in the terms of its transition cocycle $\{g_{\alpha\beta}\}$. For instance, there's a construction, due to Bott (see the last chapter of book \cite{BottTu} for example.) This Bott's construction is closely related to the one, presented here, which is also explained in paper \cite{mypapLF}.
\hfill\hfill$[]$\end{rem}

\subsection{Another geometric construction of twisting cochain associated with principal bundles}
\label{sectconnect}
It goes without saying that there exist more than only one way to write down a twisting cochain, corresponding to a princiapl bundle. In effect, using various gauge transformations, as it was explained in the previous section, one can multiply the given twisting map in many different ways. This subsection is devoted to a discussion of one particular way of deforming the twisting cochain we constructed in section \ref{sect31}. We shall first describe this deformation independently as an alternative variant of the twisting cochain associated with a principal bundle and then explain how it is related (in fact equivalent) to the previous one.

Let $\uc=\{U_i\}$ be the trivializing open cover of the base $X$ of the principal bundle $P$ with structure group $G$ (we used the same open cover throughout this and previous sections.) Let $\omega\in\DR{P}\otimes\g$ be a connection $1$\/-form on $P$. Choosing local sections $s_\alpha:U_\alpha\to\pi^{-1}(U_\alpha)$, we can pull $\omega$ down to a collection of local gauge potentials $A_{\alpha}=\sum_kA_{\alpha}^k\otimes X_k$, where $X_k$ is a linear base of \g\ and $A_{\alpha}^k$ are $1$\/-forms on $U_\alpha$. Let $F_\alpha=\sum_kF^k_\alpha\otimes X_k$ be the local curvature forms, corresponding to $A_\alpha$. Consider the following map $\xi:\DR{G}\to\check{C}(\uc,\,\DR{U})$ (in this formula we omit the summation sign with respect to the index $k$)
\begin{equation}
\label{nah-neh-nah}
\xi(x)=\{A^k_\alpha X_k(x)(e)-F^k_\alpha I_k(x)(e)\}_\alpha+\{1_{\alpha\beta}x(e)-x^0(g_{\alpha\beta})\}_{\alpha\prec\beta}.
\end{equation}
Here $e$ is the unit element of the group $G$, so that $\epsilon:x\mapsto x^0(e)=\epsilon(x)$ is the counit of the coalgebra $\DR{G}$, where $x^0$ denotes the projection of an element $x$ to the space of $0$\/-th degree forms (i.e. for a homogeneous element $x,\ \mathrm{deg}\,x=p$, $x^0=x,\ p=0$ and $x^0=0$, if $p>0$.) As before $I_k$ denotes the contraction (inner product) of an element $x\in\DR{G}$ with the left-invariant vector field, generated by $X_k$. Observe then that $\xi(x)=0$ for all $x,\ \mathrm{deg}\,x\ge1$.

Let us check, that $\xi(x)$ is in fact a twisting cochain. It is enough to check the equation $d\xi-\xi\ast\xi=0$ only for the arguments $x$ of degrees not greater than $2$ (indeed since $\xi(x)=0$ for $\mathrm{deg}\,x>1$, the same is true for $d\xi=0$ while $\xi\ast\xi$ should vanish for $\mathrm{deg}\,x\ge2$.) So, if $\mathrm{deg}\,x=0$, we get (in this and succeeding formulas we omit the restriction signs)
\begin{equation*}
\begin{split}
d(\xi)(x)&=(d+\delta)\xi(x)-\xi(dx)=\{dA_\alpha^kX_k(x)(e)-F_\alpha^kX_k(x)(e)\}_\alpha\\
             &\quad+\{A_\beta^kX_k(x)(e)-A_\alpha^kX_k(x)(e)-dx(g_{\alpha\beta})\}_{\alpha\prec\beta}\\
             &\quad+\{1_{\alpha\beta\gamma}x(e)-x(g_{\beta\gamma})+x(g_{\alpha\gamma})-x(g_{\alpha\beta})\}_{\alpha\prec\beta\prec\gamma}\\
            &=\{C_{ij}^kA_\alpha^i\wedge A_\alpha^j X_k(x)(e)\}_\alpha+\{A_\beta^kX_k(x)(e)-A_\alpha^kX_k(x)(e)-dx(g_{\alpha\beta})\}_{\alpha\prec\beta}\\
            &\quad+\{1_{\alpha\beta\gamma}x(e)-x(g_{\beta\gamma})+x(g_{\alpha\gamma})-x(g_{\alpha\beta})\}_{\alpha\prec\beta\prec\gamma}.
\end{split}
\end{equation*}
Here we used the identity $X_k=[I_k,\,d]$ and the Bianchi identity. Similarly, from the definition of comultiplication in the algebra of functions on a group, we have $x^{(1)}(g)x^{(2)}(h)=x(gh)$ for all $x\in\fc(G)$ and all elements $g,h\in G$ (here we use the standard notation for coproduct $x\mapsto x^{(1)}\otimes x^{(2)}$.) From this equality it follows that
$$
X(x^{(1)})(e)x^{(2)}(g)=\frac{d}{dt}_{|_{t=0}}x^{(1)}(\exp(tX))x^{(2)}(g)=\frac{d}{dt}_{|_{t=0}}x(\exp(tX)g)=X(x)(g),
$$
for all $X\in\g$. Similarly,
$$
x^{(1)}(g)X(x^{(2)})(e)=Ad_{g}X(x)(g)
$$
and
$$
X_i\ast X_j(x)(e)=X_i(x^{(1)})(e)X_j(x^{(2)})(e)=[X_i,\,X_j](x)(e).
$$
Thus we compute (in the second line we use the sign rule from \eqref{prod}, here we have to "`pull"' \v Cech $1$\/-cochain across $1$\/-form):
\begin{equation*}
\begin{split}
\xi\ast\xi(x)&=\{C_{ij}^kA_\alpha^i\wedge A_\alpha^j X_k(x)(e)\}_\alpha\\
                    &\quad+\{-A_\alpha^kX_k(x^{(1)})(e)(1_{\alpha\beta}x^{(2)}(e)-x^{(2)}(g_{\alpha\beta}))+(1_{\alpha\beta}x^{(1)}(e)-x^{(1)}(g_{\alpha\beta}))A_\beta^kX_k(x^{(2)})(e)\}_{\alpha\prec\beta}\\
                    &\quad+\{(1_{\alpha\beta}x^{(1)}(e)-x^{(1)}(g_{\alpha\beta}))(1_{\beta\gamma}x^{(2)}(e)-x^{(2)}(g_{\beta\gamma}))\}_{\alpha\prec\beta\prec\gamma}\\
                    &=\{C_{ij}^kA_\alpha^i\wedge A_\alpha^j X_k(x)(e)\}_\alpha\\
                    &\quad+\{-A_\alpha^kX_k(x)(e)+A_\alpha^kX_k(x)(g_{\alpha\beta})+A_\beta^kX_k(x)(e)-A_\beta^kAd_{g_{\alpha\beta}}X_k(x)(g_{\alpha\beta})\}_{\alpha\prec\beta}\\
                    &\quad+\{1_{\alpha\beta\gamma}x(e)-x(g_{\beta\gamma})+x(g_{\alpha\gamma})-x(g_{\alpha\beta})\}_{\alpha\prec\beta\prec\gamma}.
\end{split}
\end{equation*}
The difference of these two expressions is equal to
$$
d(\xi)(x)-\xi\ast\xi(x)=\{A_\beta^kAd_{g_{\alpha\beta}}X_k(x)(g_{\alpha\beta})-A_\alpha^kX_k(x)(g_{\alpha\beta})-dx(g_{\alpha\beta})\}_{\alpha\prec\beta}.
$$
But since $g_{\alpha\beta}A_\beta=A_\alpha g_{\alpha\beta}+dg_{\alpha\beta}$ (see \eqref{transform3},) this expression is equal to $0$.

Similarly, if $\mathrm{deg}\,x=1$, we have (we use the second Bianchi identity):
\begin{equation*}
\begin{split}
d(\xi)(x)&=\{-dF_\alpha^kI_k(x)(e)\}_\alpha+\{F_\beta^kI_k(x)(e)-F_\alpha^kI_k(x)(e)\}_{\alpha\prec\beta}\\
             &=\{-C_{ij}^kF^i_\alpha\wedge A^j_\alpha I_k(x)(e)\}_\alpha+\{F_\beta^kI_k(x)(e)-F_\alpha^kI_k(x)(e)\}_{\alpha\prec\beta}.
\end{split}
\end{equation*}
But from the definition of coproduct of de Rham forms on a group, it follows that
\begin{equation*}
\begin{split}
X_i\ast I_j(x)(e)-I_j\ast X_i(x)(e)&=X_i(x^{(1)})(e)I_j(x^{(2)})(e)-I_j(x^{(1)})(e)X_i(x^{(2)})(e)\\
                                                        &=[X_i,\,I_j](x)(e)=C_{ij}^kI_k(x)(e),
\end{split}
\end{equation*}
so using the identities, similar to the considered above, we get the following formula (in the second line we once again use the sign rule, this time we push a \v Cech $1$\/-simplex over a \emph{map} of degree $1$):
\begin{equation*}
\begin{split}
\xi\ast\xi(x)&=\{A^i_\alpha X_i(x^{(1)})(e)\wedge F^j_\alpha I_j(x^{(2)})(e)-F^j_\alpha I_j(x^{(1)})(e)\wedge A^i_\alpha X_i(x^{(2)})(e)\}_\alpha\\
                   &\quad+\{-F_\alpha^kI_k(x^{(1)})(e)(1_{\alpha\beta}x^{(2)}(e)-x^{(2)}(g_{\alpha\beta}))+(1_{\alpha\beta}x^{(1)}(e)-x^{(1)}(g_{\alpha\beta}))F_\beta^kI_k(x^{(2)})(e)\}_{\alpha\prec\beta}\\
                   &=\{C_{ij}^kA^i_\alpha\wedge F^j_\alpha I_k(x)(e)\}_\alpha+\{-F_\alpha^kI_k(x)(e)+F_\beta^kI_k(x)(e)\\
                   &\qquad\qquad\qquad F_\alpha^kI_k(x^{(1)})(e)x^{(2)}(g_{\alpha\beta})-x^{(1)}(g_{\alpha\beta})F_\beta^kI_k(x^{(2)})(e)\}_{\alpha\prec\beta}\\
                   &=\{-C_{ij}^kF^i_\alpha\wedge A^j_\alpha I_k(x)(e)\}_\alpha+\{F_\beta^kI_k(x)(e)-F_\alpha^kI_k(x)(e)\}_{\alpha\prec\beta}.
\end{split}
\end{equation*}
In the last line we used the equation \eqref{transform4}: $g_{\alpha\beta}F_\beta=F_{\alpha}g_{\alpha\beta}$. So, subtracting the last formula from the previous one, we get $0$. Finally, if $\mathrm{deg}\,x=2$, then $d(\xi)(x)=0$, and
$$
\xi\ast\xi(x)=\{F_\alpha^iI_i(x^{(1)})(e)\wedge F_\alpha^jI_j(x^{(2)})(e)\}_\alpha.
$$
But since
$$
I_i(x^{(1)})(e)I_j(x^{(2)})(e)+I_j(x^{(1)})(e)I_i(x^{(2)})(e)=[I_i,\,I_j](x)(e)=0,
$$
this expression is equal to $0$. This means that $\xi$ is a twisting cochain. One can show that the twisted tensor product, constructed for this particular cochain is homotopic to the cochain complex of the principal bundle. Then from the general results on twisting cochains it follows that $\xi$ should be equivalent to the cochain $\phi_P$, described in section \ref{sect31}. We shall prove this equivalence relations directly by constructing the corresponding gauge transformation. To this end, let us consider the following formula: $x\mapsto\{\exp{(\sum_kA_\alpha^kI_k)}(x)(e)\}_\alpha$. From the example \ref{counterex} it follows that this map induces a gauge transformation of the complex $\check C(\uc,\,\DR{U})$, such that
$$
\bigl(d(c)(x)+c\smile\phi_P(x)\bigr)^0=\bigl(d(c)(x)\bigr)^0=\{A^i_\alpha X_i(x)(e)-F^i_\alpha I_i(x)(e)\}_\alpha=\bigl(\xi\smile c(x)\bigr)^0.
$$
Here $(\,)^k$ denotes the $k$\/-\v Cech degree part of a cochain. Indeed since both sides of this equality are derivations in $x$, commuting with $d$ and are represented by the sums of finite number of locally-defined terms (i.e. of terms that lie in $\DR{U_\alpha}$), it is enough to check this equality only for $\mathrm{deg}\,x=0$, and only on local level, which is explained in example \ref{counterex}. Further, since $(\phi_P(x))^k=(\xi(x))^k=0,\ k\ge2$ and $(c(x))^k=0, k\ge1$, we should compare only the $1$\/-\v Cech degree part of the corresponding equality. Once again, it is enough to do it only for $\mathrm{deg}\,x=0$, which is simple.

Observe that the established relation between $\xi$ and $\phi_P$ also shows that the equivalence class of $\xi$ (modulo the gauge transformations) does not depend on the particular choice of connection and its local presentations $A_\alpha$. In addition, one can use this relation to find the homotopies between the expressions representing the characteristic classes of a principal bundle in terms of a connection and its curvature and the expressions, involving the cocycle $g_{\alpha\beta}$: just apply the universal construction of the chain homotopy between the characteristic maps $\tilde \phi_P$ and $\tilde\xi$ induced by a gauge transformation to the particular example described in this section and keep in mind that only the terms with $F_\alpha$ in $\tilde\xi$ will contribute to the image of $x\in\DR{G}$, when $\mathrm{deg}\,x>0$.
\begin{rem}\rm
To get a more conceptual picture of the two twisting cochains, described in this paper, one can use the following observation, similar to the note, concerned with example \ref{counterex}. Let $\widetilde{\DR{G}}$ be the vector bundle over the base $X$ with fibre $\DR{G}$, associated with the principal bundle $P$ (we let $G$ act on $\DR{G}$ by left translations). Let $\Gamma(\widetilde{\DR{G}})$ denote the space of differential forms on $X$ with values in this bundle. Consider the map $\tilde c:\Gamma(\widetilde{\DR{G}})\to\DR{P}$, given on the local level by the formula:
$$
f_\alpha\omega^{i_1}\wedge\dots\wedge\omega^{i_p}\otimes h_\alpha(x)\mapsto f_\alpha\exp{(\sum_i I_i\otimes A_\alpha^i)}(\omega^{i_1}\wedge\dots\wedge\omega^{i_p})\otimes h_\alpha(x).
$$
Here $x\in X$ is a generic point and $h_\alpha(x)$ is a differential form on $U_\alpha$. This map is well-defined, i.e. it doesn't depend on the choice of trivialization: when we pass from $U_\alpha$ to $U_\beta$, we get the same result. Moreover, it is an isomorphism of vector spaces (it is enough to prove this on local level, which is evident -- its inverse is given by $\exp{(-\sum_i I_i\otimes A_\alpha^i)}$.) The map $\tilde c$ intertwines the de Rham differential on $P$ with a covariant derivative on the differential (we use de Rham differential in the fibre) space of sections $\Gamma(\widetilde{\DR{G}})$. Now the twisting cochain $\xi$ is just the result of the construction of $\phi_P$ applied to $\Gamma(\widetilde{\DR{G}})$, which we equip with the differential, pulled back from $\DR{P}$.
\hfill$[]$\end{rem}
\begin{rem}\rm
The twisting cochain $\xi$ has one important advantage: one can restrict its domain to the subspace of \emph{right-invariant} differential forms on $G$, which can be thought of as the Chevallay-Eilenberg complex $C^*(\g)$ of the Lie algebra of $G$. Observe, that the Lie algebra \g\ acts on $C^*(\g)$ by Lie derivatives $L_X$, and that one can extend this action by the convolutions (inner multiplications) with elements of \g\ to obtain the collection of graded derivatives $L_X,\ I_X$, verifying the usual Cartan identities. Besides this, since for all $p,\ C^p(\g)$ is a finite-dimensional vectorspace, the map
$$
G\times C^p(\g)\to C^p(\g),\ (g,\, c)\mapsto c^g,
$$
where $c^g(X_1,\dots,X_p)=c(Ad_g(X_1),\dots,\,Ad_g(X_p))$, can be regarded dually as a map
$$
C^p(\g)\to C^p(\g)\otimes C^\infty(G).
$$
Thus, the following formula is indeed a formula of the \emph{twisting map} in the sense of section \ref{sectnext}, defined on the Chevvallay complex $C^*(\g)$, associated with the principal bundle $P$ and taking values in $\check C^*(\uc,\DR{U})$:
\begin{equation}
\label{Chevtwist}
c\mapsto\{\sum_i\{A_i\otimes L_{X^i}c-F_i\otimes I_{X^i}c\}\}_\alpha+\{1_{\alpha\beta}c(0)-c^{g_{\alpha\beta}}\}_{\alpha\prec\beta}.
\end{equation}
\hfill$[]$\end{rem}


\section{Globalization of twisting cochains}
\label{seccc}
This section is devoted to the description of two different approaches to the following natural question: ``Given a twisting cochain with values in \v Cech-de Rham algebra of the base, is it possible to find a global twisting cochain, e.g. a cochain with values in the algebra $\DR{X}$ of de Rham forms of the base or another \textit{commutative} algebra, equivalent to it?'' We give two constructions of this sort. One of them gives a twisting cochain with values in $\DR{X}$, but it uses such an inexplicit object as partition of unity. The other one gives a cochain with values in the \emph{Dupont's} algebra -- a commutative algebra, which appears as a result of the realization of the simplicial algebra of de Rham forms on the open subsets see sections \ref{sect41} and \ref{sect42} (some further details of this construction may be found in papers \cite{Dupont, mypapLF}.)

\subsection{From \v Cech complex to de Rham complex}
\label{sectinfty}
In this subsection we shall describe a construction, turning twisting cochains (and twisting maps) with values in \v{C}ech complexes into twisting cochains with values in the algebra of de Rham forms on the base. This will be based on the general methods and ideas of \cite{SmSOp}, which we shall suitably modify, or rather specialize to the complexes, used in our paper. The formulas, which we shall obtain here, will be used in the next section.

First of all, observe that there is an evident chain map $i':\DR{X}\to\check C(\uc,\DR{U})$, given by
$$
h\stackrel{i'}{\mapsto}\{h_{|_{U_\alpha}}\}_\alpha.
$$
In previous section we called this map "`localization"' of a differential form. It is a homomorphism of algebras, in fact, although the algebra $\check C(\uc,\DR{U})$ is not graded-commutative, its $0$\ \v Cech-degree part is, thus there is no contradiction here. Moreover, this map establishes an isomorphism on the level of cohomology. This can be seen very easily, from the spectral sequence, induced by the natural filtration on the bicomplex $\check C(\uc,\DR{U})$: since the sheaf $\DR{U}$ is feeble, its cohomology is concentrated in the $0$\/-th degree, and is isomorphic to the $\DR{X}$. It is also easy to see that the differential $d_1$ on $E_1$ term of this spectral sequence coincides with the de Rham complex.

Thus we conclude that the map $i'$ is a quasi-isomorphism of differential graded algebras and from the general theory of such algebras it follows that there exists a homotopy inverse $A_\infty$\/-map $\mathcal P:\check C(\uc,\DR{U})\to\DR{X}$, (explanation of what we mean by $A_\infty$\/-map is given below in remark \ref{inftymaps}, if you need more detailed explanations, you can refer for example to the book \cite{SmSOp}.) Our purpose now is to explain how one can use this map to define a $\DR{X}$\/-valued twisting cochain on $\DR{G}$. First of all, let us recall what is an $A_\infty$\/-homomorphism of two algebras and explain how it can be used to modify twisting cochains. We prefer to formulate this in the form of a remark rather than definition.
\begin{rem}\rm
\label{inftymaps}
We shall say that an $A_\infty$ map between two algebras $A$ and $A'$ is given, if there is a graded coalgebra homomorphism between their bar resolutions $B(A)\to B(A')$ commuting with differentials. Since $B(A')$ is a free coalgebra, this amounts to a series of maps $\pc(n):A^{\otimes n}\to A'$, verifying certain relations, given by commutator with bar-complex differentials (in our case $A=\check C(\uc,\DR{U})^{\otimes n},\ A'=\DR{X}$.) Now given an $A_\infty$\/-map $\mathcal P=\{\mathcal P(n)\}$, where $\mathcal P(n):A^{\otimes n}\to A'$, and a $A$\,-valued twisting cochain $\phi:K\to A$, we can define an $A'$\/-valued twisting cochain $\pc\circ\phi_P:K\to A'$ as follows:
$$
\pc\circ\phi(x)=\sum_{n\ge1}\pc(n)(\phi(x^{(1)})\otimes\dots\otimes\phi(x^{(n)})),
$$
where $x\to x^{(1)}\otimes\dots\otimes x^{(n)}$ is the abbreviation for the $n-1$\/-fold coproduct of an element $x\in\DR{G}$. Observe, that since $\mathrm{deg}\,(x^{(1)}\otimes\dots\otimes x^{(n)})=\mathrm{deg}\,x$ and the map $\pc(n)$ decreases the degree by $n-1$, there are only finitely many nonzero terms in this sum. It is then easy to derive identity \eqref{flatc} for $\pc\circ\phi$ from the corresponding identity for $\phi$ and the the properties of $A_\infty$ map \pc.

In effect, similar construction works in case when we change the coalgebra $K$ for a one homotopy-equivalent to it. Recall that there exist notions of $A_\infty$\/-coalgebras and their homomorphisms. One can regard them as dual to the notion of $A_\infty$\/-algebra and $A_\infty$\/-algebra homomorphism, described above. For instance, if $K$ and $K'$ are two coalgebras, then one can define an $A_\infty$\/-coalgebra homomorphism $c\mathcal P:K'\to K$ as a sequence of linear maps $c\pc(n):K'\to K^{\otimes n}$, which induce a homomorphism on the level of cobar-resolutions of $K'$ and $K$. If now $\phi:K\to\Omega$ is a twisting cochain, then one can show that the map $\phi\circ c\pc$
$$
\phi\circ c\pc(k)=\sum_{n\ge0}m_n(\underbrace{\phi\otimes\dots\otimes\phi}_n)c\pc(n)(k)
$$
here $m_n$ is the $n-1$\/-fold multiplication in $\Omega$ (observe that we have to assume the finiteness of the sum on the right hand side, or to introduce some sort of topology on the right hand side, to provide the convergence of this sum.)

Also observe that all $A_\infty$ algebra and coalgebra homomorphisms induce the maps of chain complexes:
$$
1\otimes\pc:K\otimes_\phi\Omega\to K\otimes_{\pc\circ\phi}\Omega'\ \mbox{and}\ 1\hat\otimes\pc:K\hat\otimes_\phi\Omega\to K\hat\otimes_{\pc\circ\phi}\Omega'
$$
and similarly for $c\pc$. In fact these maps are given by the formulas
\begin{align*}
1&\otimes\pc(k\otimes a)=\sum_{n\ge0}k^{(1)}\otimes\pc(n)(\phi(k^{(2)})\otimes\dots\otimes\phi(k^{(n)})\otimes a),\\
\begin{split}
1&\hat\otimes\pc(k\otimes a)\\
  &=\sum_{n\ge1}\sum_{m=1}^n k^{(m)}\otimes\pc(n)(\phi(k^{(m+1)})\otimes\dots\otimes\phi(k^{(n)})\otimes a\otimes\phi(k^{(1)})\otimes\dots\otimes\phi(k^{(m-1)})).
\end{split}
\end{align*}
\hfill$[]$\end{rem}
Thus, we need to find the higher maps that constitute the $A_\infty$\/-map, homotopy inverse to $i'$ and use them to obtain proper twisting cochain. In general, this would demand solving inductively all the equations, that define the maps $\pc(n)$. However, in the case we discuss here, it is possible to write down the map $\pc\circ\phi_P$ explicitly without solving all these equations. In fact, we shall obtain a formula for the $A_\infty$\/-map \pc\ expressed in the terms of certain linear maps. To this end, let us first define a chain map $p:\check C(\uc,\,\DR{U})\to\DR{X}$, homotopy inverse to $i'$. Note that $p$ need not commute with multiplication at all, we only want it commute with differentials. To define such $p$ we regard $\check C(\uc,\,\DR{U})$ as a bicomplex with two differentials, $\delta$ and $d$ given by the \v Cech and the de Rham differentials respectively. Let us forget about the differential $d$ (both in $\check C(\uc,\,\DR{U})$ and $\DR{X}$) for a while, tat is, consider $\DR{X}$ as a complex with zero differential. Then the map $i_0=i'$ still induces an isomorphism in cohomology, but this time it is possible to give an explicit inverse map. Namely, let $\{\varphi_\alpha\}$ be a partition of unity on $X$ associated with the open cover \uc. Then we have the following map
\begin{align*}
p_0:&\check C(\uc,\,\DR{U})\to\DR{X},\\
p_0(&\{h_{\alpha_1\dots\alpha_n}\}_{\alpha_1\prec\dots\prec\alpha_n})=
\begin{cases}
&\!\!\!0,\ n\ge2,\\
&\!\!\!\sum_\alpha h_\alpha\varphi_\alpha,\ n=1.
\end{cases}
\end{align*}
It is easy to check that $p_0\circ i_0=Id$ and hence $p_0$ is a homotopy inverse map for $i_0$. Moreover, it is possible to write down an explicit formula for the contracting homotopy. To this end we first extend a \v Cech cochain $\{h_{\alpha_1\dots\alpha_n}\}_{\alpha_1\prec\dots\prec\alpha_n}$ to an arbitrary (i.e. not necessarilly ordered) collection of $n$ indices $\alpha_1,\dots\alpha_n$ so that $h_{\alpha_{\sigma(1)}\dots\alpha_{\sigma(n)}}=(-1)^\sigma h_{\alpha_1\dots\alpha_n}$. We shall denote the anti-symmetrization procedure by $A$. Clearly, there is a one-to-one correspondence between the anti-symmetrized \v Cech cochains and the usual (ordered) \v Cech cochains. We shall define the inverse procedure, associating to every antisymmetric \v Cech cochain $A(h)$ its regular part $h$ by $S$. Then the \v Cech differential $\delta$ can be extended to the space of anti-symmetrized cochains by the same formula:
$$
\delta(\{h_{\alpha_1\dots\alpha_n}\})_{\alpha_1\dots\alpha_n\alpha_{n+1}}=\sum_i(-1)^ih_{\alpha_1\dots\widehat{\alpha_i}\dots\alpha_{n+1}}.
$$
It is easy to see that $\delta$ commutes with the anti-symmetrization, i.e. that $\delta A(h)= A(\delta h)$ (in fact, one can define the differential of an anti-symmetrized \v Cech cochain by this rule.) Moreover one can define the $\cup$\/-product of two antisymmetric \v Cech cochains so that $A(h'\cup h'')=A(h')\cup A(h'')$ (altough we shall not need this here.)

Having this in mind, let us define the following map:
\begin{align*}
H_0:&\check C^*(\uc,\,\DR{U})\to\check C^{*-1}(\uc,\,\DR{U}),\\
H_0(&\{h_{\alpha_1\dots\alpha_n}\})_{\alpha_1\prec\dots\prec\alpha_n\prec\alpha_{n-1}}=S\bigl(\sum_{\alpha_n}A(\{h_{\alpha_1\dots\alpha_n}\})\varphi_{\alpha_n}\bigr).
\end{align*}
In a more explicit way, this formula can be written as follows
$$
H_0(\{h_{\alpha_1\dots\alpha_n}\})_{\alpha_1\prec\dots\prec\alpha_n\prec\alpha_{n-1}}=\sum_{i=1}^{n-1}(-1)^{i-1}\Bigl\{\sum_{\alpha_{i-1}\prec\alpha\prec\alpha_i}h_{\alpha_1\dots\alpha_{i-1}\alpha\alpha_i\dots\alpha_n}\varphi_\alpha\Bigr\}+(-1)^{n-1}\sum_{\alpha_{n-1}\prec\alpha}h_{\alpha_1\dots\alpha_{n-1}}\varphi_\alpha.
$$
Then the following formulas can be easily checked:
\begin{align}
\label{perturb}
H_0\, H_0&=0, & H_0\, i_0&=0, & p_0\, H_0&=0, & i_0\, p_0-Id&=\delta\, H_0+H_0\,\delta.
\end{align}
Taking into consideration these equations and the equality $p_0\circ i_0=Id$, mentioned above, we observe that we are now in the situation, described by the so-called ``perturbation lemma'', see \cite{SmSOp}. This lemma ensures us that for an additional differential $d$ on $\check C(\uc,\,\DR{U})$, one can find a differential $d'$ on $\DR{X}$, and maps $i,\ p$ and $H$ with the same domains and ranges as those of $i_0,\ p_0$ and $H_0$ respectively, which will verify the same equations as before but with respect to the differential $d+\delta$ on $\check C(\uc,\,\DR{U})$ and the differential $d'=d'+0$ on $\DR{X}$ (in fact, if there were a nonzero differential $\delta'$ on $\DR{X}$ before perturbation, we would have to use $d'+\delta'$ here.) These maps are given by the formulas
\begin{equation}
\label{perteqn}
\begin{split}
d'&=p_0\,d\,i_0+p_0\,d\,H_0\,d\,i_0+\dots, \qquad i=i_0+H_0\,d\,i_0+H_0\,d\,H_0\,d\,i_0+\dots,\\
p&=p_0+p_0\,d\,H_0+p_0\,d\,H_0\,d\,H_0+\dots, \qquad H=H_0+H_0\,d\,H_0+H_0\,d\,H_0\,d\,H_0+\dots.
\end{split}
\end{equation}
In our case, since $i_0=i'$ commutes with the de Rham differential $d$ on $\DR{X}$ and $\check C(\uc,\,\DR{U})$, we see that $d'=d$ and $i=i_0=i'$. Thus $p$ is the homotopy inverse map for $i=i'$ we were looking for.

Further, recall that both cochain complexes $\check C(\uc,\,\DR{U})$ and $\DR{X}$ are in fact differential algebras. From a general theory (see \cite{SmSOp}, page 95), it follows that given maps $i,\ p$ and $H$, satisfying the equations \eqref{perturb} and the equality $p\circ i=Id$, one can construct a structure of $A_\infty$\/-algebra $\pi=\{\pi(n)\}$ on $\DR{X}$, two $A_\infty$ maps $\mathcal I=\{\ic(n)\}:\DR{X}\to \check C(\uc,\,\DR{U})$ and $\pc=\{\pc(n)\}:\check C(\uc,\,\DR{U})\to\DR{X}$ (where we regard $\check C(\uc,\,\DR{U})$ as an $A_\infty$\/-algebra with the $A_\infty$\/-structure, induced from the associative multiplication of the \v Cech cochains) and an $A_\infty$\/-homotopy $\hc=\{\hc(n)\}$. This collection verifies the usual set of relations for homotopy-inverse $A_\infty$\/-maps and the contracting homotopy that is, the equations similar to \eqref{perturb} and the equation $\pc\circ\ic=Id$, but all the differentials and compositions should be taken in the sense of $A_\infty$\/-algebras. These maps are given by the following formulas (after Smirnov we shall denote the multiplication on $\check C(\uc,\,\DR{U})$ by $\tilde\pi$, but we shall regard all the maps here as maps of the chain complexes, and not as maps of their suspensions, as it is done in the cited book):
\begin{equation}
\label{ainfexpl}
\begin{split}
\pi(n)&=p\tilde\pi(H\,\tilde\pi\otimes1+1\otimes H\,\tilde\pi)\dots(H\,\tilde\pi\otimes1\otimes\dots\otimes1+\dots+1\otimes\dots\otimes1\otimes H\,\tilde\pi)(i\otimes\dots\otimes i),\\
\ic(n)&=H\,\tilde\pi(H\,\tilde\pi\otimes1+1\otimes H\,\tilde\pi)\dots(H\,\tilde\pi\otimes1\otimes\dots\otimes1+\dots+1\otimes\dots\otimes1\otimes H\,\tilde\pi)(i\otimes\dots\otimes i),\\
\pc(n)&=p\,\tilde\pi(H\,\tilde\pi\otimes1+1\otimes H\,\tilde\pi)\dots\\
           &\quad\dots(H\,\tilde\pi\otimes\dots\otimes1+\dots+1\otimes\dots\otimes H\,\tilde\pi)(H\otimes i\,p\otimes\dots\otimes i\,p+\dots+1\otimes\dots\otimes1\otimes H),\\
\hc(n)&=H\,\tilde\pi(H\,\tilde\pi\otimes1+1\otimes H\,\tilde\pi)\dots\\
           &\quad\dots(H\,\tilde\pi\otimes\dots\otimes1+\dots+1\otimes\dots\otimes H\,\tilde\pi)(H\otimes i\,p\otimes\dots\otimes i\,p+\dots+1\otimes\dots\otimes1\otimes H).
\end{split}
\end{equation}
In the first two of these formulas the number of tensors $i$ in the last brackets is equal to $n$, and the last two of these formulas consist of the brackets of the form $H\,\tilde\pi\otimes1\otimes\dots\otimes1+\dots+1\otimes\dots\otimes1\otimes H\,\tilde\pi$, except for the last one, where we should put the sum $H\otimes i\,p\otimes\dots\otimes i\,p+1\otimes H\otimes i\,p\otimes\dots\otimes i\,p+\dots+1\otimes\dots\otimes1\otimes H$, i.e. the tensor $H$ moves from left to the right so that on the left hand of it there stand $1$-\/s, and on the right hand there stand maps $i\,p$. The total number of tensor legs in every term of this sum is $n$.

Now, since the map $i=i'$ is in fact a homomorphism of differential graded algebras, one readily sees that $\pi(n)=0$ for $n\ge3$ and $\pi(2)$ coincides with the usual product in $\DR{X}$. Similarly, $\ic(n)=0$ if $n\ge2$ and $\ic(1)=i$. Thus we conclude that the formulas \eqref{ainfexpl} give an explicit way to write down the homotopy inverse $\pc$ to the map $i'$. In principle, one can use these formulas to get few first terms of a twisting map $\pc\circ\phi_P$ explicitly.

\subsection{Simplicial spaces, algebras and their realizations}
\label{sect41}
In this and the following section we shall describe the globalization procedure in terms of the so-called Dupont's realization of a cosimplicial algebra (see \cite{Dupont} and \cite{mypapLF} for references.) We shall apply this procedure to the algebra of functions on the open subsets of a manifold. The main advantage of this approach is that this realization can be regarded as a locally euclidean topological manifold with explicitly given partition of unity (see discussion in the end of the next section.)

Let us first of all recal few definitions and results from the theory of simplicial sets and algebras. The standard reference for this subject is J.~P.~May's book \cite{MaySO}. Let $\Delta$ denote the simplicial ctegory, i.e. the category of (nonempty) finite ordered sets and their non-decreasing maps. Recall that for a category $\mathcal C$, a \emph{simplicial object} in $\mathcal C$ is a functor $X_*:\Delta^{op}\to\mathcal C$. If $X_*,\ Y_*$ are two simplicial objects in \cc, then a simplicial map $X_*\to Y_*$ is just the natural transformation of functors from $X_*$ to $Y_*$. Similarly, \emph{cosimplicial object} in \cc\ is a covariant functor $X^*:\Delta\to\cc$, and a cosimplicial map $X^*\to Y^*$ is again a natural transformation.

In more explicit terms, let $\mathbf n$ denote the ordered set $\{0,\,1,\,2,\dots,n\}\in\Delta$. Clearly, one can think that the objects of $\Delta$ coincide with the collection of all the oredered sets $\mathbf n,\ n\in\mathbb N_0$. Here $\mathbb N_0=\mathbb N\cup\{0\}$. Moreover, one can describe all the morphisms in this category in a pretty combinatoric way. Namely, let $\delta_i(n)=\delta_i,\ i=0,\dots,n$ be the unique non-decreasing injection $\mathbf{n-1}\to\mathbf{n}$, which omits the element $i\in\mathbf{n}$ and $\sigma_j(n)=\sigma_j,\ j=0,\dots,n$ the unique non-decreasing surjection $\mathbf{n+1}\to\mathbf{n}$, which hits the element $j\in\mathbf{n}$ twice. It turns out that all the morphisms in $\Delta$ are generated by $\delta_i(n)$ and $\sigma_j(n)$ for all $n,\ i$ and $j$.

Thus to define a simplicial object $X_*$ in \cc\ it is necessary and sufficient to determine a sequence of objects $\{X_n\}_{n\in\mathbb N_0}$ and two collections of \cc\/-morphisms between them: maps in the first collection will increase the dimension by $1$
$$
X_*(\delta_i)=\partial_i:X_n\to X_{n-1},\ i=0,\dots,n,
$$
they are called \emph{faces}, and those from the other collection will decrease the dimension,
$$
X_*(\sigma_j)=s_j:X_n\to X_{n+1},\ j=0,\dots,n,
$$
they are called \emph{degeneracies}. There are certain commutator relations, that must be verified by these maps. We shall not discuss them here. Everybody can obtain these relations by a direct inspection of definitions (see the first section of \cite{MaySO}.) In the terms of this approach, a simplicial map $f_*$ between $X_*$ and $Y_*$ is given by a collection of maps $f_n:X_n\to Y_n$ between the equal levels of $X_*$ and $Y_*$, commuting with the structure maps.

Similarly, cosimplicial object $X^*$ in \cc\ is determined by a collection of objects $\{X^n\}_{n\in\mathbb N_0}$ and maps
$$
X^*(\delta_i)=\delta^i:X^{n-1}\to X^n,\ \mbox{ and }\ X^*(\sigma_j)=\sigma^j:X^{n+1}\to X^n,\ i,j=0,\dots,n,
$$
called \emph{coface} and \emph{codegeneracy} respectively. They verify similar relations, as before, and as before we can regard cosimplicial maps as a collection of maps $f^n$, commuting with cofaces and codegeneracies.

Our primary interest is related with the categories \sca\tc\ and $co\sca\tc$ of simplicial and cosimplicial objects in the category \tc\ of topological spaces and their continuous maps, which are called \emph{simplicial spaces} (in particular, regarding an arbitrary set as a discrete topological space, we can embed the categories of simplicail and cosimplicial sets into \sca\tc\ or $co\sca\tc$.) For the objects in this category there exist a well-defined notion of \emph{geometric realization}. Under certain conditions on a simplcial space $X_*$ (see \cite{Segal74} for details,) its realization is defined by the formula
\begin{equation*}
|X_*|=\coprod_{i\ge0}X_i\times \bigtriangleup^i /\{(\theta_*(x),\,t)\sim(x,\,\theta^*(t))\}.
\end{equation*}
Here $\bigtriangleup^i$ is the standard geometric $i$\/-dimensional simplex, one can think of it as of the subset in $\mathbb R^{n+1}$, determined by conditions $t_0+\dots+t_n=1,\ t_i\ge0,\ i=0,\dots,n$. Observe that the collection of topological spaces $\{\bigtriangleup^i\}_{i\ge0}$ is in a natural way endowed with the structure of cosimplicial space, just identify the set of vertices of an $i$\/-dimensional simplex with $\mathbf i$ and notice that every map of vertices can be uniquely extended to a map of simpleces by linearity. If $\theta$ is a morphism $\mathbf i\to \mathbf j$ in $\Delta$, then $\theta_*:X_j\to X_i$ and $\theta^*:\bigtriangleup^i\to\bigtriangleup^j$ will denote the induced maps. Clearly, the geometric realization is a well-defined functor from the category of simplicial spaces to topological spaces. In particular, every simplicial map $f:X_*\to Y_*$ determines a continuous map $|f|:|X_*|\to|Y_*|$.

For example, let $\mathcal U=\{U_1,\dots,U_N\}$ be a finite (ordered) open cover of a topological space $B$. The order is given by the indices $1,\dots,N$ and wherever we have a string $\dots,U_\alpha,U_\beta,U_\gamma,\dots$ of elements of \uc, we shall assume that their order is not violated, i.e. $\ldots\preceq U_\alpha\preceq U_\beta\preceq U_\gamma\preceq\ldots$. Consider now the following simplicial topological space $\nc\uc_*$, closely related to the \emph{nerve} of \uc:
$$
\nc\uc_n=\coprod_{\{U_0,\dots, U_n\},\ U_i\in\uc}U_0\cap\dots\cap U_n.
$$
The face and degeneracy maps in $\nc\uc_*$ will be determined by the formulas
$$
\partial_p(x)=i_p(x),\qquad\qquad s_q(x)=j_q(x),
$$
$p,q=0,\dots,n$, where
$$
i_p:U_0\cap\dots\cap U_p\cap\dots\cap U_n\to U_0\cap\dots\cap\widehat U_p\cap\dots\cap U_n
$$
is the natural inclusion and
$$
j_p:U_0\cap\dots\cap U_p\cap\dots\cap U_n\to U_0\cap\dots\cap U_p\cap U_p\cap\dots\cap U_n=U_0\cap\dots\cap U_p\cap\dots\cap U_n
$$
is the identity map. The simplicial relations are readily verified.
\begin{df}
We shall call the topological space $\nc\uc_*$, described here the \emph{topological nerve} of the open cover \uc. The adjective ``topological'' is used to distinguish this construction from the usual nerve of an open cover of a space, which we shall call \emph{combinatoric}.
\end{df}

Observe that if \vc\ is a finer open cover, i.e. if for every $V_k\in\vc$ there exists an element $U_j\in\uc$, such that $V_k\subseteq U_j$, then one can define a simplicial map $\nc\vc_*\to\nc\uc_*$ by choosing $U_j$ with prescribed property for every $V_k\in\vc$ and extending the inclusion maps to all the levels of $\nc\vc_*$ in a natural way. In particular, if \bc\ is the open cover, consisting of the only one open set $U=B$, then $\nc\bc_*=B_*$ is the constant simplicial space, i.e. $b_n=B$ for all $n$ and all the structure maps are equal to identity, and from our discussion it follows that there exists a (unique) simplicial map $(Q_\nc)_*:\nc\uc_*\to\nc\bc_*=B_*$ for every open cover \uc.

On the other hand, since every simplex in $B_n$ is degenerate for $n\ge1$, so the realization of $B_*$ is canonically homeomorphic to $B$ itself. In this way we obtain a continuous map $Q_\nc:|\nc\uc_*|\to B$ for every open cover \uc\ of $B$. Our purpose is to prove the following statement
\begin{prop}
\label{hequivtop}
For every finite open cover\ \uc\ of $B$ the map $Q_\nc:|\nc\uc_*|\to B$ is a homotopy equivalence.
\end{prop}
\begin{Proof}
We shall prove this statement by induction on the numer $N$ of elements in \uc. The base of induction $N=1$ is clear. Let $N=2$, so $\uc=\{U_1,\,U_2\},\ B=U_1\cup U_2$. Let $B_1=U_1,\ B_2=U_2$ and $B_{12}=U_1\cap U_2$. The open cover \uc\ of $B$ generates the following open covers on $B_1,\ B_2$ and $B_{12}$: $\uc_1=\{U_1\},\ \uc_2=\{U_2\}$ and $\uc_{12}=\{U_1\cap U_2\}$. Then we have the following diagramm
\begin{equation}
\begin{CD}
{\nc(\uc_{12})_*}@>>>{\nc(\uc_1)_*}\\
@VVV                   @VVV\\
{\nc(\uc_2)_*}@>>>{\nc\uc_*}.
\end{CD}
\end{equation}
Clearly, this square is cocartesian in the category of simplicial spaces. So applying the geometric realization functor to this diagramm yields a homotopy cocartesian square in topological category (it is a general fact about geometric realizations on simplicial spaces, that it sends cartesian and cocartesian squares to homotopy cartesian/cocartesian squares, one can prove it by replacing the simplicial spaces by the equivalent category of bisimplicial sets and using the fact that diagonal realization of bisimplicial sets is homotopy equivalent to the consequitive realization in one and another simplicial directions
.) The maps $p_{\uc_{12}},\ p_{\uc_1},\ p_{\uc_2}$ and $Q_\nc$ commute with this diagramm and send it to the cocartesian square of spaces
\begin{equation}
\begin{CD}
{B_{12}}@>>>{B_1}\\
@VVV               @VVV\\
{B_2}@>>>{B}.
\end{CD}
\end{equation}
But as we know, the maps $p_{\uc_{12}},\ p_{\uc_1}$ and $p_{\uc_2}$ are homeomorphisms, in particular they are homotopy equivalences. Then the map $Q_\nc$ is also a homotopy equivalence.

The general case is treated in a similar way: if the cover \uc\ consists of $N$ elements, $\uc=\{U_1,\dots,U_N\}$, we put $B_1=U_1\cup\dots\cup U_{N-1}$, $B_2=U_N$ and $B_{12}=B_1\cap B_2$ and consider the cartesian square similar to what we had before.
\end{Proof}\\ 

Assume now that we have a contravariant functor $\Omega$ from topological spaces to commutative differential graded algebras. Applying it degree-wise to the levels of a simplicial space $X_*$, we obtain a \emph{co}simplicial object in the category of commutative differential graded algebras, or a cosimplicial CDGA for short. More concretely, let us suppose that the space $B$ in the previous construction was a smooth manifold. We can choose the functor $\Omega$ to be the functor of de Rham differential forms on manifolds. It is our purpose to define \emph{the realisation} of $\Omega(X_*)$ in the category of commutative differential algebras. In short, $|\Omega(X_*)|$ should be a commutative differential graded algebra, which will play the role of de Rham forms on the realisation of $X_*$; in particular, it should be an algebra homotopy equivalent to the de Rham algebra of the manifold $B$.

It is our purpose now to define this sort of realization. To this end we consider first the algebraic building blocks, similar to the geometric simplices $\bigtriangleup^i$ of geometric realisation. This time they should be differential algebras, modelling the de Rham forms on geometric simplices. This can be done directly by putting $\Omega_i=\DR{\bigtriangleup^I}$, but we shall prefer a more delicate approach, coming from Sullivan, Thom and others, and use the polynomial differential forms on simplices. So we put
$$
\Omega_n=\mathbb R[t_0,\dots,t_n;\,dt_0,\dots,dt_n]/(t_0+t_1+\dots+t_n=1),
$$
where $\mathbb R[t_0,\dots,t_n;\,dt_0,\dots,dt_n]$ is the free graded-commutative algebra, generated by the elements $t_i,\,dt_i,\ deg\, t_i=0,\ deg\,dt_i=1,\ i=0,\dots,n$ with differential $d$, generated by the relation $d(t_i)=dt_i$, which we factorise by the differential ideal, generated by the relation $t_0+\dots+t_n=1$ (e.g. it follows from this relation that $dt_0+\dots+dt_n=0$, etc.). Thus we obtain the the collection of differential graded algebras $\{\Omega_n(\bigtriangleup)\}_{n\ge0}$. We introduce the structure of simplicial differential graded algebra on it by the formulas
\begin{align*}
\partial_k(t_i)&=\begin{cases}
                            t_i,\ &i<k,\\
                            0,\ &i=k,\\
                            t_{i-1},\ &i> k;
                            \end{cases}\\
s_j(t_i)&=\begin{cases}
                  t_i,\ &i<j,\\
                  t_i+t_{i+1},\ &i=j,\\
                  t_{i+1},\ &i>j.
                  \end{cases}
\end{align*}
The usual simplicial relations are easily checked. As we have mentioned above, one can regard this simplicial differential algebra as the algebra of polynomial differential forms on geometric simplices.

Now in order to define the realization of a cosimplicial commutative differential graded algebra, we consider the example in which it is given by a contravariant functor applied to a simplicial space (e.g. by de Rham algebras on simplicial manifold,) consider the definition of the geometric realisation of a simplicial space. It is equal to the factor-space of the disjoint union of cartesian products $X_i\times\bigtriangleup^i$. Reasoning dually now we conclude that the realization we seek, should be equal to a subalgebra of the direct product of algebras $\DR{X_i}\otimes\Omega_i$. In order to find the conditions, determining the necessary subalgebra inside $\prod_{i=0}^\infty\DR{X_i}\otimes\Omega_i(\bigtriangleup)$, we look at the relations, used to define the geometric realization: an element of our realization should give a well-defined function on $|X_*|$. Thus we obtain the following definition:
\begin{df}
Let $\{\tilde\Omega^i\}_{i\ge0}$ be a cosimplicial commutative differential graded algebra (in particular, $\tilde\Omega^i=\DR{X_i}$ for a simplicial manifold $X_*$). Then we shall call a realization of $\{\tilde\Omega^i\}$ the DG algebra $|\tilde\Omega^*|$, defined as follows
$$
|\tilde\Omega^*|=\Bigl\{\{\omega_i\otimes\varphi_i\}\in\prod_{i=0}^\infty\tilde\Omega^i\otimes\Omega_i(\bigtriangleup)|\theta^*(\omega_i)\otimes\varphi_i=\omega_j\otimes\theta_*(\varphi_j)\Bigr\}.
$$
\end{df}

\vspace{1mm}
Here, as before $\theta:\mathbf i\to\mathbf j$ is a morphism in $\Delta$, and $\theta_*,\ \theta^*$ -- the corresponding simplicial and cosimplicial structure maps in $\Omega_*(\bigtriangleup)$ and $\tilde\Omega^*$. Clearly, since both maps $\theta_*$ and $\theta^*$ are homomorphisms of algebras, $|\tilde\Omega^*|$ is a commutative subalgebra of the direct product algebra (since the multiplication is defined degree-wise, there is no problems with convergeance).

Suppose now that the cosimplicial algebra $\{\tilde\Omega^*\}$ is equal to $\{\DR{\nc\uc_i}\}$ where \uc\ is an open cover of a smooth manifold $B$. Then there exists an evident map $Q_\nc^*:\DR{B}\to|\DR{\nc\uc_*}|$, given by
$$
\omega\mapsto\{(\bigoplus_i\omega_{|_{U_i}})\otimes 1,(\bigoplus_{i\le j}\omega_{|_{U_i\cap U_j}})\otimes 1,\dots\}\in\prod_{i=0}^\infty\DR{\nc\uc_i}\otimes\Omega_i(\bigtriangleup).
$$
Here $\uc=\{U_1,\dots,U_N\}$ and $1\in\Omega_0(\bigtriangleup)$ is the unit element, and the projection of $Q_\nc^*(\omega)$ on $\DR{\nc\uc_k}\otimes\Omega_k(\bigtriangleup)$ is equal to the direct sum of its restrictions on all the possible intersections $U_{i_0}\cap\ldots\cap U_{i_k}$ for all sequences of indices ${i_0\le\dots\le i_k}$, tensored by the identity in $\Omega_k(\bigtriangleup)$.

Clearly, the image of $Q_\nc^*$ belongs to $|\DR{\nc\uc_*}|$. Indeed, it is enough to show that $\delta_p(Q_\nc^*(\omega)_k)=Q_\nc^*(\omega)_{k+1}$ and $\sigma_q(Q_\nc^*(\omega)_k)=Q_\nc^*(\omega)_{k-1}$  for all the indices $k\ge 0$, $p=0,\dots, k+1$ and $q=0,\dots,k-1$. Here $Q_\nc^*(\omega)_k$ denotes the projection of $Q_\nc^*(\omega)$ onto $\DR{\nc\uc_k}\otimes\Omega_k(\bigtriangleup)$. But $\delta_p=(\partial_p)^*=i_p^*$ and 
$$
i_p^*(\omega_{|_{U_0\cap\dots\cap\widehat U_p\cap\dots\cap U_{k+1}}})=\omega_{|_{U_0\cap\dots\cap U_{k+1}}}.
$$
In the same way $\sigma_q=(s_q)^*=j_q^*$ and we have a similar identitty. Observe, that since the restriction of differential forms is a homomorphism of algebras, $Q_\nc^*$ is a homomorphism of algebras.

The following proposition is an algebraic analogue of the proposition \ref{hequivtop}.
\begin{prop}
For any finite open cover \uc\ of a manifold $B$, the map $Q_\nc^*$ is a quasi-isomorphism of differential graded algebras, $Q_\nc^*:\DR{B}\to|\DR{\nc\uc_*}|$.
\end{prop}
\begin{Proof}
Recall that a homomorphism of differential graded algebras is called quasi-isomorphism, if it induces an isomorphism in cohomology. So we should show that the map in cohomology, induced by $Q_\nc^*$ is an isomorphism. Once again we can use induction on the number $N$ of elements in \uc. If $N=1$, then from the very definition of $|\DR{\nc\uc_*}|$ it follows that $Q_\nc^*$ is an isomorphism of algebras. In the general case, we use the Meier-Vietoris exact sequence and the 5-lemma (and of course we use the fact that the corresponding diagramm is cartesian in the category of differential graded algebras.)
\end{Proof}\\

\subsection{Relation with previous constructions}
\label{sect42}
Let $\pi:P\to B$ be a principal bundle over a manifold $B$ with structure group $G$, and let \uc\ be an open cover of a manifold $B$, trivializing the bundle $P$. In this section we shall discuss a construction of the twisting cochain on $G$ with values in the realization $|\DR{\nc\uc_*}|$ of the cosimplicial algebra $\DR{\nc\uc_*}$. Then it will not be difficult to show that (under certain identification) this new cochain will be equivalent to the cochains, described in the previous sections. In fact, the way we shall obtain the formulas in this section is by transition of the previously found cochains into the algebra $|\DR{\nc\uc_*}|$ with the help of a suitable $A_\infty$\/-map. So first of all we shall need to compare the \v Cech complex of the previous paragraph and the algebra $|\DR{\nc\uc_*}|$. They cannot be isomorphic, since the former one is a non-commutative differential algebra while the latter one is a graded commutative algebra. However, we have the following result
\begin{prop}
The algebra $|\DR{\nc\uc_*}|$ is homotopy equivalent to $\check C(\uc,\,\DR{U})$.
\end{prop}
\begin{Proof}
In effect this statement is already proven above. Just recall, that two differential algebras $\ac=\check C(\uc,\,\DR{U})$ and $\bc=|\DR{\nc\uc_*}|$ are homotopy equivalent, if there exists a third differential algebra \cc, and two quasi-isomorphisms $\ac\leftarrow\cc$ and $\cc\to\bc$. In our case these ingredients can be chosen as follows: $\cc=\DR{B}$, $\cc\to\bc$ is the homomorphism $Q_\uc^*$, and $\ac\leftarrow\cc$ is the inclusion $\DR{B}\to\check C(\uc,\,\DR{U})$, given by $\varphi\mapsto\{\varphi_{|_{U_\alpha}}\}_\alpha$.
\end{Proof}\\ \\ \vspace{1mm}
For our purposes we shall need a more detailed description of the homotopy equivalence $\check C(\uc,\,\DR{U})\sim|\DR{\nc\uc_*}|$. To this end let us consider the map $w:\check C(\uc,\,\DR{U})\to|\DR{\nc\uc_*}|$, given by the following rule. First we define it on the \v Cech cochains of degree $0$ (here we omit the $\otimes$ sign between elements of $\DR{\nc\uc_*}$ and $\Omega_*(\bigtriangleup)$; we also omit $1\in\Omega_*(\bigtriangleup)$ from our formulas):
\begin{align*}
w(\{h_\alpha\}_\alpha)&=\{\bigoplus_i h_i,\,\bigoplus_{i\le j}\Bigl((h_i)_{|_{U_{ij}}}t_0+(h_j)_{|_{U_{ij}}}t_1\Bigr),\dots\}\\
\intertext{so that on the $n$\/-th simplicial level we have}
w(\{h_\alpha\}_\alpha)_n&=\bigoplus_{i_0\le\dots\le i_n}\Bigl((h_{i_0})_{|_{U_{i_0\dots i_n}}}t_0+\dots+(h_{i_n})_{|_{U_{i_0\dots i_n}}}t_n\Bigr).
\end{align*}
Let us show that this element belongs to $|\DR{\nc\uc_*}|$; we do it here for purely pedagogical purposes: it would be quite difficult to prove all the formulas of this papragraph in the same way, so below we shall give a more conceptual proof. As before, it is enough to check only the relations, involving the (co)face and (co)degeneracy operations. So if we apply the coface operation in $\DR{\nc\uc_*}$ direction, $\delta_p\otimes 1,\ p=0,\dots,n+1$, to $w(\{h_\alpha\}_\alpha)_n$, we obtain
\begin{align*}
i^*_p\otimes 1(\bigoplus_{i_0\le\dots\le i_n}&\{(h_{i_0})_{|_{U_{i_0\dots i_n}}}t_0+\dots+(h_{i_n})_{|_{U_{i_0\dots i_n}}}t_n\}\\
&=\bigoplus_{i_0\le\dots\le i_{p-1}\le j\le i_p\dots\le i_n}\{(h_{i_0})_{|_{U_{i_0\dots j\dots i_n}}}t_0+\dots+(h_{i_n})_{|_{U_{i_0\dots j\dots i_n}}}t_n\},
\end{align*}
where $j$ on the right hand side stands in the $p$\/-th place and the term with $h_j$ is missing. On the other hand, if we apply $1\otimes\partial_p$ to
\begin{align*}
w(\{h_\alpha\}_\alpha)_n&=\bigoplus_{i_0\le\dots\le i_{n+1}}\{(h_{i_0})_{|_{U_{i_0\dots i_n}}}t_0+\dots+(h_{i_n})_{|_{U_{i_0\dots i_n}}}t_n\}\\&=\bigoplus_{i_0\le\dots\le i_p\le\dots\le i_{n+1}}\{(h_{i_0})_{|_{U_{i_0\dots i_p\dots i_{n+1}}}}t_0+\dots+(h_{i_{n+1}})_{|_{U_{i_0\dots i_p\dots i_{n+1}}}}t_{n+1}\},
\end{align*}
we shall obtain the same result, this follows from the definition of $\partial_p(t_j)$ given above. Similarly, since for all $q=0,\dots, n-1$ we have
\begin{align*}
j_q^*(&\{(h_{i_0})_{|_{U_{i_0\dots i_n}}}t_0+\dots+(h_{i_n})_{|_{U_{i_0\dots i_n}}}t_n\})\\
          &=\begin{cases}
                (h_{i_0})_{|_{U_{i_0\dots i_n}}}t_0+\dots +(h_{i_q})_{|_{U_{i_0\dots i_n}}}(t_q+t_{q+1})+\dots+(h_{i_n})_{|_{U_{i_0\dots i_n}}}t_n\},\ i_q=i_{q+1},\\
                 0,\ \mbox{otherwise,}
              \end{cases}
\end{align*}
we see that applying $\sigma_q\otimes 1$ to $w(\{h_\alpha\}_\alpha)_n$, we obtain the same result as applying $1\otimes s_q$ to $w(\{h_\alpha\}_\alpha)_{n-1}$.

Further on for $x=\{h_{\alpha\beta}\}_{\alpha\prec\beta}$, we put
\begin{align*}
w(x)&=\{0,\,w(x)_1,\,w(x)_2,\,\dots\},\\
\intertext{where}
w(x)_1&=\bigoplus_{i\le j}\Bigl(t_0d(t_1h_{ij}-t_1d(t_0h_{ij}))\Bigr)=\bigoplus_{i\le j}-dt_0h_{ij},\\
\begin{split}
w(x)_2&=\bigoplus_{i\le j\le k}\Bigl(t_0d(t_1(h_{ij})_{|_{U_{ijk}}}+t_2(h_{ik})_{|_{U_{ijk}}})\\
            &\qquad\quad+t_1d(-t_0(h_{ij})_{|_{U_{ijk}}}+t_2(h_{jk})_{|_{U_{ijk}}})\\
            &\qquad\quad+t_2d(-t_0(h_{ik})_{|_{U_{ijk}}}-t_1(h_{jk})_{|_{U_{ijk}}})\Bigr).
\end{split}\\
\intertext{In general, in the $n$\/-th simplicial degree there will stand}
w(x)_n&=\bigoplus_{i_0\le\dots\le i_n}\Bigl(\sum_{k=0}^nt_kd(\sum_{j=0}^nt_j(\tilde h_{i_ki_j})_{|_{U_{i_0\dots i_n}}})\Bigr).
\end{align*}
Here we put $\tilde h_{kl}=h_{kl},\ k< l$ and $\tilde h_{kl}=-h_{kl}$ otherwise, in particular, $\tilde h_{ii}=0$. In general, for $x=\{h_{\alpha_0\dots\alpha_n}\}_{\alpha_0\prec\dots\prec\alpha_n}\in\check C(\uc,\,\DR{U})$, the element $w(x)$ will consist of the following elements $w(x)_p$ belonging to the $p$\/-th simplicial degree (for the sake of brevity we omit the signs of restricions of differential forms, defined on some open subsets, to the smaller subsets):
\begin{equation}
\label{hequi}
w(x)_p=\begin{cases}
                  &\!\!\!0,\qquad\ p<n,\\
                  &\!\!\!\bigoplus_{i_0\le\dots\le i_p}\Bigl\{\sum_{k_0=0}^pt_{k_0}d\Bigl\{\sum_{k_1=0}^pt_{k_1}\dots d\Bigl\{\sum_{k_n=0}^pt_{k_n}(-1)^{\sigma(k_0,\dots,k_n)}h_{k_0\dots k_n}\Bigr\}\dots\Bigr\}\Bigr\},\ p\ge n.
                 \end{cases}
\end{equation}
It is possible to check that these formulas give well-defined elements in $|\DR{\nc\uc_*}|$ and that the map in general does commute with the differentials. This way of proving is related with harsh computational difficulties. However there is a more conceptual way of looking on $|\DR{\nc\uc_*}|$, which makes the map $w$ a particular case of the maps, described in the  paragraph \ref{sectinfty}. In fact, this is precisely the way these formulas were obtained.

The idea is very simple: let us look on the geometric realization $|\nc\uc_*|$ as on topological space, glued from open subsets of $B$ \emph{and cylinders} $U_{\alpha\beta}\times\bigtriangleup^1,\ U_{\alpha\beta\gamma}\times\bigtriangleup^2,\dots$; we consider only nondegenerate simplices $(\alpha\beta\gamma)$ in the combinatoric nerve of the open cover \uc\ here. These subspaces are not open, unless we throw away the boundaries of the simplices, but we can use them to define an open cover of $|\nc\uc_*|$: put
$$
\widetilde{U_\alpha}=\bigcup_{\alpha\in\sigma}U_\sigma\times\bigtriangleup^{|\sigma|}_\alpha.
$$
Here the union is taken over all the simplices in combinatorical nerve of \uc, which contain $\alpha$, and $\bigtriangleup^{|\sigma|}_\alpha$, where $|\sigma|$ stands for the dimension of $\sigma$, denotes the union of all the open hyper-faces of $\bigtriangleup^{|\sigma|}$, containing $\alpha$ as one of their vertices. In other words, $\bigtriangleup^{|\sigma|}_\alpha=\bigtriangleup^{|\sigma|}-\partial_i\bigtriangleup^{|\sigma|}$, where $\sigma=(\alpha_0\prec\dots\prec\alpha_{i-1}\prec\alpha\prec\alpha_{i+1}\prec\dots\prec\alpha_{|\sigma|})$.

More generally, one can extend this definition to arbitrary combinatoric simplex $\tau=(\beta_0,\dots,\beta_k)$:
$$
\widetilde{U_\tau}=\bigcup_{\tau\subseteq\sigma}U_\sigma\times\bigtriangleup^{|\sigma|}_\tau,
$$
where $\bigtriangleup^{|\sigma|}_\tau$ denotes the union of all (hyper-)faces of $\bigtriangleup^{|\sigma|}$, which contain $\tau$ as their subface. Then one can embed $U_\tau$ into $\widetilde{U_\tau}$ regarding it as the ``bottom'' of all the cylinders in the union, i.e. putting all the coordinates $t_i$, which correspond to the vertices not in $\tau$, equal to $0$. Moreover the following equality holds:
$$
\widetilde{U_{\tau_1}}\bigcap\widetilde{U_{\tau_2}}=\widetilde{U_{\tau_1\cup\tau_2}}.
$$
Thus sending $U_\tau$ to $\widetilde{U_\tau}$, we establish a morphism $\nc^{comb}\uc_*\to\nc\widetilde{\uc}(|\nc\uc_*|)_*$, where we use the open cover $\widetilde{\uc}=\{\widetilde{U_\alpha}\}$ to define the nerve of the geometric realization $|\nc\uc_*|$. Here $\nc^{comb}$ denotes the combinatoric nerve of \uc.

Applying the functor of de Rham forms to both sides of this construction, we obtain a homomorphism of chain complexes (and even of differential graded algebras): $\tilde i:\check C(\uc,\,\DR{U})\to\check C(\widetilde{\uc},\,\DR{\widetilde{U}})$, which sends the differential form $x$ on $U_\tau$ to the differential form $\tilde x$ on $\widetilde{U_\tau}$, equal to $x\otimes 1$.
\begin{prop}
The map $\tilde i$ is a quasi-isomorphism.
\end{prop}
\begin{Proof}
It is enough to observe that its composition with $i':\DR{B}\to\check C(\uc,\,\DR{U})$, see section \ref{sectconnect}, is equal to $\widetilde{i'}\circ Q_\nc^*$. Here $\widetilde{i'}$ is the map analoguous to $i'$ defined on $|\DR{\nc\uc_*}|$. Since all the other maps here are quasi-equivalences, so is $\tilde i$.
\end{Proof}\\ \\ \vspace{1mm}
Now one can construct the map $w$ as the composition of $\tilde i$ and the inverse $\tilde p$ of $\widetilde{i'}$. In particular, in this way one obtains the formulas \eqref{hequi} and we can conclude that they are well-defined and commute with differentials.
\begin{rem}\rm
One can give an explicit formula for the homotopy inverse map of $w$. To this end consider the definite integral $\int_{\bigtriangleup^n}:\Omega_n\to\mathbb R$. It is easy to check, using the Stokes formula
$$
\int_{\bigtriangleup^n}d\omega=\int_{\partial\bigtriangleup^n}\omega,
$$
that the map
$$
v=\prod_n n!(Id\otimes\int_{\bigtriangleup^n}):|\DR{\nc\uc_*}|\to\check C(\uc,\,\DR{U})
$$
commutes with the differentials. On the other hand, its composition with $w$ is equal to the identity. In fact $\int_{\bigtriangleup^n}t_0dt_1\dots dt_n =\frac{1}{n!}$, while $\mathrm{deg}_t\,w(x)_n\le\mathrm{deg}\,x$ and the degrees are equal, only when $\mathrm{deg}\,x=n$, so $v$ is actually homotopy inverse to $w$.
\hfill$[]$\end{rem}
\begin{rem}\rm
In effect, one can use the similar construction to define a map, homotopy inverse to $Q_\nc^*$ (see previous section.) It is enough just to replace $t_i$ by $\varphi_i$, where $\{\varphi_i\}$ is a partition of unity, associated with \uc.
\hfill$[]$\end{rem}


\section{Comparison with the Getzler-Jones classes}
\label{sect40}
In this section we shall explain in what sense the cyclic Chern class of Getzler, Jones and Petrack, Bismut and others (see \cite{GJP},
and references therein) are variants of the twisting cochain constructions, described in the previous sections. To this end we shall need to make a couple of intermediate steps, which relate the twisting cochains and the corresponding chain maps of the previous sections to the Getzler-Jones-Petrack's results.

\subsection{The Getzler-Jones-Petrack map}
\label{sect43}
Let $E$ be a rank $n$ complex vector bundle over a compact manifold $X$. One can associate to it the principal $GL(n)$\/-bundle, frame bundle $P_E$ of $E$, so that $E=P_E\times_{GL(n)} \mathbb C^n$. One can apply to $P_E$ the methods, described in previous sections and get formulas for the twisting cochain, for the characteristic classes, etc.. In particular, one can obtain the map
$$
\tilde{\hat\phi}:K\hat\otimes_\phi\Omega\to NH^*(\Omega),
$$
where $K$ is a differential coalgebra model of the group $GL$ (for instance the coalgebra of polynomial differential forms) and $\Omega$ -- a differential algebra model of the base $X$, see proposition \ref{propleftright} and discussion that follows it. For example, one can take $K=\Omega_{DR}(G)$ and $\Omega=\check C(\uc,\DR{U})$. Using the gluing construction, described in the paragraph \ref{sectinfty}, one can substitute $\DR{X}$ instead of $\check C(\uc,\DR{U})$. Besides this as we explained in section \ref{seccc} (first of all paragraph \ref{sectinfty},) the domain of the map can be replaced with the quasi-equivalent complex $\DR{\tilde P_E}$, where $\tilde P_E$ is the associated gauge bundle of $P_E$.

On the other hand, the fact, that $P_E$ is associated with a vector bundle can be used to obtain another map with the same domain and range, given in pretty explicit terms. Here we shall describe this construction. Let $E\leftrightarrows X\times\mathbb C^N$ be an inclusion/projection of $E$ to/from a  trivial bundle. Let $p\in Mat_N(C^\infty(X))$ be the corresponding projector. One can use it to construct the morphism we need as follows.

First, let $X\times GL(N)$ be the gauge bundle of the trivial vector bundle $\underline N=X\times\mathbb C^N$. Since $E$ is inside $\underline N$, we can regard $\tilde P_E$ as a subbundle of $X\times GL(N)$. Namely chose the complementary projection $1-p=q\in Mat_N(C^\infty(X))$, such that $pq=qp=0$ and $p\oplus q=1$. Then one can identify $\tilde P_E$ with the following subbundle of $X\times GL(N)$:
$$
\check P_E=\{(x,g)\in X\times GL(N)|g(Im\, p(x))=Im\,p(x),\, g(Im\,q(x))=1_{Im\,q(x)}\}.
$$
Equivalently, the latter condition can be written in the matrix form as follows: $gq(x)=q(x)g=q(x)$. Let us denote the complementary subbundle by $\bar E$

In the terms of functions and differential forms on the bundles, this inclusion of gauge bundles induces a restriction epimorphism $\DR{X\times GL(N)}\to \DR{\tilde P_E}$, the fact that it is epimorphic is clear from the local considerations. To render these and following considerations rigoruos, one should consider the algebras of vertically-polynomial differential forms $\Omega^*_{vpoly}(\tilde P_E),\ \Omega_{vpoly}^*(X\times GL(N))$ etc., see section \ref{sect17}, i.e. all the forms should belong to $\Omega_{poly}(G)$ on each fibre. This is what we shall always assume from now on, although we shall not encumber our text with this redundant notation below.

One can now regard $\DR{\tilde P_E}$ as a factor algebra of $\DR{X\times GL(N)}$ modulo the kernel of the restriction map. On the other hand one can describe the generating set of relations of the kernel, thus we obtain a representation of $\DR{\tilde P_E}$ as factor-algebra of $\DR{X\times GL(N)}$ by the ideal generated by the following set of functions 
\begin{equation}
\label{eqideal}
\sum_j u_{ij}p_{jk}-u_{ik},\ \sum_j p_{ij}u_{jk}-p_{ik},\ i,k=1,\dots,N.
\end{equation}
In effect, there is an alternative point of view on this construction, which probably makes it easier for understanding: consider the evident diagramm of group bundles:
\begin{equation}
\label{eqdia1}
\xymatrix{ & {\quad \tilde P_{E,\bar E}}\ar[dl]_-{p_1}\ar[dr]^-{i_{E,\bar E}} &\\
{\tilde P_E} & & {X\times GL(N)=\tilde P_{\underline N}.}
}
\end{equation}
Here $\tilde P_{E,\bar E}=\tilde P_E\times_X\tilde P_{\bar E}$, $p_1$ is the projection on the first factor, and $i_{E,\bar E}$ is the natural embedding of $\tilde P_{E,\bar E}$ into $\tilde P_{\underline N}$, induced by the embeddings of $E$ and $\bar E$ into $\underline N$. The left leg of this diagramm is a fibration, and the right one -- an inclusion. Then one can regard the functions on $\tilde P_E$ as a function on $\tilde P_{\underline N}$ whose restriction to $\tilde P_{E,\bar E}$ doesn't depend on coordinates in the $\bar E$ directions.

Recall that choosing an embedding $E\to\underline N$, enables one to write down global expressions, representing connection on $E$ and its curvature (see \cite{GJP}): 
$$
A=pdp-(1-p)d(1-p),\ F=pdpdp.
$$
These expressions should be regarded as $Mat_N(\mathbb C)$\/-valued differential forms on $X$, such that $A_{x,\xi}(Im\, p(x))\subseteq Im\, p(x)$ and $F_{x,\xi,\eta}(Im\, p(x))\subseteq Im\, p(x)$ for all points $x\in X$ and all vectors $\xi,\,\eta\in T_xX$. On the other hand, regarding $Mat_N(\mathbb C)=\mathfrak{gl}_N$ as the space of left-invariand vector fields on $GL_N$, one can think of $A$ and $F$ as of $\DR{X}$\/-valued derivatives on $\DR{X\times GL(N)}$ and apply to them the construction of twisting map, desctribed in the example \ref{example1}. The following proposition is easy but important for the understanding of the following results.
\begin{prop}
One can restrict the maps $A,\ F$ to the factor-algebra $\DR{\tilde P_E}=\DR{X\times GL(N)}/I$, where $I$ is the ideal, generated by relations \eqref{eqideal}, so that they become graded derivatives on this algebra with values in $\DR{X}$.
\end{prop}
\begin{Proof}
It is enough to check, that the kernel of the factorisation map $\DR{X\times GL(N)}\to\DR{\tilde P_E}$ is killed by the derivatives $A$ and $F$. But this is clear from the fact that this kernel is nothing but the ideal, generated by the matrix elements of the difference $1_N-p(x)$.
\end{Proof}\\ \\ \vspace{1mm}
Thus one can use $A$ and $F$ (or rather their restrictions to $\DR{\tilde P_E}$) to define a map 
$$
\tilde\phi_E:\DR{\tilde P_E}\to NH(\DR{X}),
$$
given by the formula, similar to the one we used in section \ref{sectbegin}:
$$
\tilde\phi_E(a)=\sum_{n\ge0}\phi_E^{\otimes n}(a)|_{1_{GL(N)}}.
$$
Here $a$ is an arbitrary differential form on $\tilde P_E$, $\phi_E$ is the twisting map, determined by $A$ and $F$, and $1_{GL(N)}$ is the section of $\tilde P_E$, corresponding to the identic morphism of the bundle $E$. From now on we shall call the map $\tilde\phi_E$ \emph{the Getzler-Jones-Petrack map}. As a matter of fact, in order to obtain the genuine construction of \cite{GJP}, one needs to introduce additional equivariant structure. We shall discuss this refinement in section \ref{sect44}. But before we do it, let us discuss the algebraic properties of the map $\tilde\phi_E$.
\begin{prop}
\label{prt_GJP}
 The map $\tilde\phi_E$ commutes with differentials. If one considers the usual multiplication in $\DR{\tilde P_E}$ and the cyclic shuffle product in $NH(\DR{X})$, then this map is a homomorphism of graded algebras. If one changes the projections $p$ and $q$ that determine $\tilde\phi_E$, then this map will be changed by an inner conjugation of its domain by a section of the gauge bundle.
\end{prop}
\begin{Proof}
 All the properties listed here follow from the results of sections \ref{sect2} and \ref{sect3}. In effect, one readily sees that the map $\tilde\phi_E$ is just disguised the gauge bundle monodromy map $\widetilde T_A$, see section \ref{sect17}, thus all the properties of the latter applies to the former.
\end{Proof}\\ \\ \vspace{1mm}
\begin{rem}
One can introduce the map
\begin{equation}
\label{eqcomulte}
\DR{\tilde P_E}\to\DR{\tilde P_E}\otimes_{\DR{X}}\DR{\tilde P_E},
\end{equation}
induced by the comultiplication in the algebra $\DR{GL(N)}$: one needs to check that the ideal $J$ generated by the elements of the projector $p$ is maped to $J\otimes\DR{GL(N)}+\DR{GL(N)}\otimes J$ by this comultiplication, which indeed follows from the equation $p^2=p$. Then by a slight modification of the reasonings of section \ref{sectbegin}, one shows that the map $\tilde\phi_E$ intertwines the comultiplication \eqref{eqcomulte} and the map $NH(\DR{X})\to NH(\DR{X})\otimes_{\DR{X}}NH(\DR{X})$, given by \eqref{coprodom2}.\hfill$[]$
\end{rem}

\subsection{The comparison theorem}
\label{sectcompth}
The main purpose of the present section is to prove the theorem \ref{omagn}.
\begin{theorem}
\label{omagn}
Under the identifications made in previous sections, the Jones-Getzler map $\tilde\phi_E$ is homotopic to the map $\tilde{\hat\phi}$.
\end{theorem}
\begin{Proof}
Before we begin proving this statement, let us make clear, which identifications we imply at. First of all, if $\phi$ is a twisting cochain on $K$ with values in $\Omega$, then $\tilde{\hat\phi}$ will denote the map $K\hat\otimes_\phi\Omega\to NH(\Omega)$ (see \eqref{eqloopmap},) determined by this data. In our case we take $\Omega$ to be $\DR{X}$ or $|\DR{\nc\uc_*}|$ (we chose these complexes for the sake of commutativity of the corresponding differential algebras and also because the map $\tilde\phi_E$ takes values in $NH(\DR{X})$) and take $\phi$ to be the twisting cochain on $\Omega_G$ with values in $\Omega$. To this end we need to start with the twisting cochain $\phi_P$ or $\xi$, described above, with values in \v Cech complex and use the higher homotopy maps, described in the previous section to pass from the \v Cech cochains to differential forms. Observe, that in fact the statement of this theorem doesn't depend on the choice of $\Omega$ and $K$, taking into consideration the homotopy equivalence of all the constructions. That is the map $\tilde{\hat\phi}$ is homotopic to $\tilde\phi_E$ for all $\Omega$ and $K$.

Second, we need to use the homotopy equivalence of the proposition \ref{propleftright} that relates the complex $K\hat\otimes_\phi\Omega$ with $\DR{\tilde P_E}$. In effect, we have assumed that $P$ is a principal $GL(n)$\/-bundle, so we take $K$ to be a model of this group (e.g. $K=\Omega_{poly}^*{GL(n)}$; we shall often abbreviate $\Omega_{poly}^*{GL(n)}$ to $\Omega_G$ below.) Besides this we have chosen $E$ to be the vector bundle associated with $P$ by the canonical $n$\/-dimensional complex representation of $GL(n)$. Then $P\times_{Ad} GL(n)=\tilde P_E$ as bundles of groups. On the other hand, one can obtain $P$ as the reper bundle of $E$, which we shall denote by $P_E$. In effect, we shall show that the equivalence of \ref{propleftright} intertwines the coproduct structures \eqref{eqcomulte} and \eqref{coprodom1}, (up to homotopy) so that $\tilde{\hat\phi}$ and $\tilde\phi_E$ will be equivalent as homomorphisms of coalgebras.

First of all, let us prove the homotopy equivalence of $\DR{\tilde P_E}$ and $\Omega_G\hat\otimes_\phi\DR{X}$ as coalgebras over $\DR{X}$. To this end we consider two intermediate objects. First object is the groupoid space $\mathscr GP=\hat P\times_G P$ (see the definition of $\hat P$ in the proof of proposition \ref{propleftright}; here and everywhere below $G$ and $P$ will mean the group $GL(n)$ and the principal bundle $P_E$ respectively.) Projections $s,r:\mathscr GP\to X$ are given by the projections on the first and the second factor respectively, and the composition is given by the following rule: 
$$
(u,\,v)\ast(x,\,y)=(ug,\,y),
$$
for all $u,x\in\hat P,\ v,y\in P$ such that $\pi(v)=\hat\pi(x)$ and $g=g_1g_2$, where $\varphi_\alpha(v)=(g_1,\,\pi(v)),\ \hat\varphi_\alpha(x)=(\hat\pi(x),\,g_2)$ ($\varphi_\alpha,\ \hat\varphi_\alpha$ are the local trivializations in $P$ and $\hat P$ respectively.) Then $P\times_{Ad}G=\tilde P_E$ can be embedded as a subbundle of $\mathscr GP$ lying above the diagonal of $X\times X$.

The second object we shall consider here is an algebraic analog of $\mathscr GP$, the complex $GP=GP(\Omega,\,K)=\Omega\,{}_\phi\!\otimes K\otimes_\phi\Omega$. It is equal to the tensor product $\Omega\otimes K\otimes\Omega$ as a vector space, while the differential is given by the following formula (we omit the tensor signs for the sake of brevity)
$$
{}_\phi\, d_\phi(a\, k\, b)=da\,k\,b+(-1)^{|a|}(a\,dk\,b+a\phi(k^{(1)})\,k^{(2)}\,b)+(-1)^{|a|+|k|}(a\,k\,db+(-1)^{|k^{(2)}|}a\,k^{(1)}\,\phi(k^{(2)})b).
$$
This complex is a bimodule over $\Omega$ when we let this algebra act on the left and right tensor legs by the left and the right multiplication respectively. Then there is a map $GP\to GP\otimes_\Omega GP$, given by
$$
a\,k\,b\mapsto (a\,k^{(1)}\, 1)\otimes_\Omega(1\,k^{(2)}\,b).
$$

Our first claim is that $GP$ is in effect an algebraic model for $\mathscr GP$, i.e. we claim that there is a quasi-equivalence $\mathscr I$ of $GP$ and the de Rham algebra of $\mathscr GP$ regarded as grupoid coalgebras over $\Omega$, i.e. this quasi-isomorphism should intertwine the ``coproducts'' described above. The first statement, i.e. that the complexes are quasi-equivalent, follows from the isomorphism
$$
GP=(\Omega\,{}_\phi\!\otimes K)\otimes^K(K\otimes_\phi\Omega),
$$
where $\otimes^K$ is the tensor product over coalgebra, in our case over the Hopf algebra $K$. Recall that $K\otimes_\phi\Omega$ is a model of $P$ as $G$\/-space. Besides this, since the action of $G$ and coaction of $K$ are free, the tensor product over $K$ and homotopy tensor product over $K$ are equivalent and we can apply the reasoning we used in the proof of proposition \ref{propleftright}. Thus, we have established the quasi-isomophism of the complexes. It is this quasi-isomorphism that we shall denote by $\mathscr I$. Observe, that the quasi-isomorphism $\mathscr J$ of proposition \ref{propleftright} is a restriction of $\mathscr I$ onto the diagonal.

Next we prove that $\mathscr I$ commutes with the action. Similarly to what we did above, one shows that $GP\otimes_\Omega GP$ is a model of $\mathscr GP\times_X\mathscr GP$. Then, since the maps
$$
k\otimes a\mapsto k^{(1)}\otimes(k^{(2)}\otimes a)\ \mbox{and}\ a\otimes k\mapsto(a\otimes k^{(1)})\otimes k^{(2)}
$$
correspond under the identifications of theorem \ref{propleftright}, which we use here, to the action of $G$ on $P$ and $\hat P$ respectively, and it is these two actions that actually determine the groupoid structure on $\mathscr GP$, we conclude that $\mathscr I$ is a map of coalgebras over $\Omega$.

Observe now that $\Omega_G\hat\otimes_\phi\DR{X}$ embeds into the corresponding $GP$ diagonally as a coalgebra over $\DR{X}$. Of course, one should take $GP$ with $\Omega=\DR{X}$. Namely, we just send $a\,k\,b$ to $(-1)^{|a||k|}k\otimes ab$. We shall denote this map by $\tilde m^*$. This embedding, on one hand, intertwines the coproduct structures on both sides, and on the other hand, it makes the following diagramm commute:
\begin{equation}
\label{diagcomm}
\begin{CD}
{GP(\DR{X},\,\Omega_G)}                  @>{\mathscr I}>> {\DR{\mathscr GP}}\\
@V{\tilde I}VV                                                                             @V{\Delta^*}VV\\
{\Omega_G\hat\otimes_\phi\DR{X}}@>{\mathscr J}>>{\DR{\tilde P_E.}}
\end{CD}
\end{equation}
Since the left vertical, the top horizontal and the right arrows commute with the coproducts, and the vertical maps are epimorphic, we conclude that the bottom arrow also should commute with coproduct.

Now we can prove that $\tilde\phi_E$ is homotopic to $\tilde{\hat\phi}$. First of all, observe, that the homotpoy equivalence of proposition \ref{propleftright} can be constructed as follows: take the space 
$$
\mathscr P(X;\,\hat P\times_G P)=\{(\gamma,\,[\hat p,\,p])|\gamma:[0;\,1]\to X,\ \gamma(0)=\pi_1([\hat p,\,p]),\ \gamma(1)=\pi_2([\hat p,\,p])\},
$$
then, on one hand this space is homotopy equivalent to $\hat P_E$, the homotopy being given by pulling back to $\mathscr P(X;\,\hat P\times_G P)$ the homotopy that contracts to $X$ the free path space of a topological space $X$. On the other hand, its cohomology can computed by the complex
$$
\DR{\hat P\times_G P}\otimes_{\DR{X}^{ev}} B(\DR{X}),\ \mbox{where}\ B(\DR{X})=\bigoplus_{n=0}^\infty \DR{X}\otimes\overline{\DR{X}}{\,}^{\otimes n}\otimes\DR{X}.
$$
Here  $\overline{\DR{X}}$ denotes the factor-space $\DR{X}/1$ and for an algebra $A$, $A^{ev}=A\otimes A^o$; $A^o$ denotes $A$ with inverted multiplication: $a\circ b=(-1)^{|a||b|}ba$. There's an evident generalization of the Chen's iterated integral map, see the construction described in the section \ref{sectbegin}, and \cite{Chen}, which gives the map $\DR{\hat P\times_G P}\otimes_{\DR{X}^{ev}} B(\DR{X}),\to\DR{\mathscr P(X;\,\hat P\times_G P)}$. In fact, to construct this map one should regard $\mathscr P(X;\,\hat P\times_G P)$ as a pullback:
\begin{equation*}
\begin{CD}
{\mathscr P(X;\,\hat P\times_G P)} @>>> {\hat P\times_G P}\\
@VVV                                                               @VVV\\
{\mathscr P(X)}                                   @>>>{X\times X}.
\end{CD}
\end{equation*}
Also observe, that instead of $\DR{\hat P\times_G P}$ one can use $GP(\DR{X},\,\DR{GL(n)})$, which models $\DR{\hat P\times_G P}$ as $\DR{X}$\/-bimodule: this substitution preserves the $\DR{X}$\/-bimodule structure.
\begin{rem}
In effect, one can also use the following complex to compute the cohomology of $\mathscr P(X;\,\hat P\times_G P)$:
$$
B(\DR{P},\,\DR{X},\,\DR{\hat P})^{co \DR{GL(n)}}=\bigoplus_{n=0}^\infty (\DR{P}\otimes\overline{\DR{X}}{\,}^{\otimes n}\otimes\DR{\hat P})^{co \DR{GL(n)}}.
$$
Here for a Hopf algebra $K$ and a $K$\/-bicomodule $M$, $M^{co K}=Ker(\triangledown:M\to K\otimes M)$, $\triangledown(m)=m^{(-1)}\otimes m^{(0)}-S(m^{(1)})\otimes m^{(0)}$. We abbreviate the left (resp. right) $K$\/-coaction on $M$ to $m\mapsto m^{(-1)}\otimes m^{(0)}$ (resp. $m\mapsto m^{(0)}\otimes m^{(1)}$.)\hfill$[]$
\end{rem}

\vspace{1mm}
In order to construct the homotopy equivalence $\mathscr P(X;\,\hat P\times_G P)\leftrightarrows \tilde P_E$, we interptrete the groupoid $\mathscr GP=\hat P\times_G P$ as follows: for any pair of points $(x,x')\in X\times X$, $\mathscr GP_{(x,x')}=Iso(E_x,\,E_{x'})$, the space of linear isomorphisms between the fibre of $E$ over $x$ and its fibre over $x'$. Let $A$ be a linear connection in $E$, for example, one can take $A$ to be a Grassmanian connection, induced by an embedding of $E$ into a trivial bundle. From now on we shall denote this grouppoide by $\mathscr G_E$ to underline its close relation with the bundle $E$. Then the homotopy equivalence is induced by the following maps:
\begin{align*}
\Phi:&\tilde P_E\to \mathscr P(X;\,\hat P\times_G P), \qquad i(a)=(\ast_{\pi(x)},\,[\hat p,p])\\
\intertext{where $\pi:\tilde P_E\to X$ is the projection, $\ast_x$ for $\pi(a)=x\in X$ is a constant path $\ast_x(t)=x,\ t\in[0;1]$, and $\hat p\in\hat P,\ p\in P$ are defined as the elements in $\hat P_{\pi(a)}$ and $P_{\pi(a)}$ respectively, which in the corresponding local charts $\hat\varphi_\alpha,\ \varphi_{\alpha}$ of $\hat P$ and $P$ respectively are represented by such elements $(\pi(a),\hat h)$ and $(h,\pi(a))$, that $\hat hh=g$ where $\tilde\varphi_\alpha(a)=(\pi(a),g)$. Its homotopy inverse is given by}
\Psi:&\mathscr P(X;\,\hat P\times_G P)\to\tilde P_E, \qquad p(\gamma,[\hat p,p])(e)=(\int_\gamma A)^{-1}[\hat p,p](e)\\
\intertext{where $e$ is an arbitrary element in the fibre of $E$, which lies above $s([\hat p,p])$ (we interprete $\hat P\times_G P$ as grouppoide space $\mathscr GP$,) $\int_\gamma A$ is the parallel transport in $E$ along the path $\gamma$ with respect to the connection $A$, and $[\hat p,p](e)$ is the action of $\mathscr GP$ on $E$, where we interprete this grouppoide as we described it above. Clearly, $\Psi\Phi=Id_{\tilde P_E}$. On the other hand, $\Phi\Psi$ is homotopic to the identity via the homotopy $H_s,\ s\in[0;1]$, given by}
H_s:&\mathscr P(X;\,\hat P\times_G P)\to\mathscr P(X;\,\hat P\times_G P) \qquad H_s(\gamma,[\hat p,p])=(\gamma'_s,(\int_{\gamma''_s}A)^{-1}[\hat p,p]).
\end{align*}
Here $\gamma'_s(t)=\gamma(st)$ and $\gamma''_s(t)=\gamma(s+(1-s)t)$ for $s,t\in[0;1]$ and we once again interprete $\mathscr GP$ as the set of linear isomorphisms between fibres of $E$.

The purpose of these considerations is in the following remark: one can write down an explicit formula for the inverse image map between the de Rham complexe of $\tilde P_E$ and a suitable model for $\mathscr P(X;\,\hat P\times_G P)$, induced by $\Psi$. To this end we first of all observe that there's a map $M:\mathscr G_E\times_{X^{op}}\mathscr C_E\to\tilde P_E$. Here
$$
\mathscr G_E\times_{X^{op}}\mathscr C_E=\{(a,b)\in\mathscr G_E\times\mathscr C_E|r(a)=s(b),\,r(b)=s(a)\},
$$
and the map is just composition of $a$ and $b$.

Next we consider the embedding of $E$ into a trivial bundle $\underline N=X\times\mathbb C^N$ we used above. Let $p$ be the projection $X\times\mathbb C^N\to E$, and $q=1-p$ be the complementary projection, defining a vector bundle $\bar E$, such that $E\oplus \bar E=\underline N$. Then there exists an evident diagramm of grouppoides $\mathscr G$, associated with the bundles $E,\ F$ and $\underline N$ (compare it with the diagramm \eqref{eqdia1}):
\begin{equation}
\label{diaggroup}
\xymatrix{ & {\quad\mathscr G_E\times_{X\times X}\mathscr G_{\bar E}}\ar[dl]_-{p_1}\ar[dr]^-{i_E\times i_{\bar E}} &\\
{\mathscr G_E} & & {X\times GL(N)\times X.}
}
\end{equation}
Here $\mathscr G_E\times_{X\times X}\mathscr G_{\bar E}$ is the fibred product of $\mathscr G_E$ and $\mathscr G_{\bar E}$ over $X\times X$:
$$
\mathscr G_E\times_{X\times X}\mathscr G_{\bar E}=\{(e,f)\in\mathscr G_E\times\mathscr G_{\bar E}|s(e)=s(f),\,r(e)=r(f)\},
$$
and $X\times GL(N)\times X$ is just $\mathscr G_{\underline N}$. The map $i_E\times i_{\bar E}$ is an embedding and $p_1$ is the evident projection onto the first factor. Below we shall use the same letters for the corresponding embeddings and projections of $\tilde P_E,\ \tilde P_{\bar E}$ and $\tilde P_{\underline N}$. Observe that when we identify $\mathscr G_{\underline N}$ with $X\times GL(N)\times X$ and $\tilde P_{\underline N}$ with $X\times GL(N)$, then $\mathscr G_{\underline N}\times_{X\times X}\mathscr G_{\underline N}$ becomes naturally isomorphic to $X\times GL(N)\times GL(N)\times X$ and the corresponding map $M$ will be given by the product of matrices:
$$
M(x,g,h,y)=(x,gh).
$$

Now let us choose a connection in $\underline N$, which preserves the subbundles $E$ and $\bar E$. For example, one can take $A=pdp-qdq$, the Grassmanian connection, generated by $p$ (see \cite{GJP}.) Then $\Psi^*$ will be equal to the following composition:
\begin{equation}
\label{constrphie}
\begin{split}
\DR{\tilde P_E}\stackrel{p_1^*}{\longrightarrow}\DR{\tilde P_E\times_X \tilde P_{bar E}}\stackrel{M^*}{\longrightarrow}&\DR{\mathscr G_E\times_{X\times X}\mathscr G_{\bar E}}\otimes_{\DR{X}^{ev}}\DR{\mathscr G_E\times_{X\times X}\mathscr G_{\bar E}}\longrightarrow
\\ &\qquad\quad\stackrel{\int_\gamma A}{\longrightarrow}\DR{\mathscr G_E\times_{X\times X}\mathscr G_{\bar E}}\otimes_{\DR{X}^{ev}}\DR{\mathscr P(X)}
\end{split}
\end{equation}
Here we used an implicit quasi-isomorphism between $\DR{(\mathscr G_E\times_{X\times X}\mathscr G_{\bar E})\times_{X^{op}}(\mathscr G_E\times_{X\times X}\mathscr G_{\bar E})}$ and $\DR{\mathscr G_E\times_{X\times X}\mathscr G_{\bar E}}\otimes_{\DR{X}^{ev}}\DR{\mathscr G_E\times_{X\times X}\mathscr G_{\bar E}}$, and $A$ is a connection in $E\oplus \bar E$, described above. Observe that this step can be realized as follows: for a form $\omega$ in $\DR{\mathscr G_E\times_{X\times X}\mathscr G_{\bar E}}$ we first consider an element $\tilde\omega$ in the differential forms on $\mathscr G_{\underline N}$, which maps into $\omega$ under $(i_1\times i_2)^*$ (see the diagramm \eqref{diaggroup},) apply $\int_\gamma A$ to $\tilde\omega$ and obsevre that the result doesn't depend on the choice of $\tilde\omega$. Finally to complete the construction we should observe in the end of this sequence of homomorphisms, that the differential forms that we obtain in $\DR{\mathscr G_E\times_{X\times X}\mathscr G_{\bar E}}$ are basic with respect to the projection $p_1$, hence they belong to $\DR{\mathscr G_E}$.

Now, if we use the Chen's iterated integral map, we shall obtain the following description of $\Psi^*$ (more accurately, of its cmposition with the inverse of the Chen's map $\sigma$, which gives quasi-isomorphism of $\DR{\mathscr P(X;\,\hat P\times_G P)}$ and $\DR{\mathscr G_E}\otimes_{\DR{X}^{ev}}B(\DR{X})$): for any polynomial in the fibre direction differential form $\omega$ on $\tilde P_E$, we choose a representative $\omega'\otimes\omega''\in\DR{X}\otimes\Omega^*_{poly}(GL(N))=\Omega_{vpoly}^*(\tilde P_{\underline N})$. Then $M^*(\omega'\otimes\omega'')=\omega'\otimes({\omega''}^{(1)}\otimes{\omega''}^{(2)})\otimes 1$. Next we apply a map, similar to the monodromy map from \ref{sectbegin} to ${\omega''}^{(2)}$, thus obtaining an element in $\DR{\mathscr G_{\underline N}}\otimes B(\DR{X})$. Finally, we apply $(i_1\times i_2)^*$ to the first tensor factor and observe that the result doesn't actually depend on coordinates in the $\bar E$ direction.

An important consequence of this construction is the following
\begin{lemma}
Let the connection $A$ in $\underline N$ be given by the grassmanian construction, associated with the projector $p$, defining $E$, then the map $\tilde\phi_E$ is equal to the following composition:
\begin{multline*}
\DR{\tilde P_E}\stackrel{\sigma^{-1}\Psi^*}{\longrightarrow}\DR{\mathscr G_E}\otimes_{\DR{X}^{ev}}B(\DR{X})\longrightarrow\\ \stackrel{I^*\otimes m}{\longrightarrow}\DR{\tilde P_E}\otimes_{\DR{X}} NH(\DR{X})\stackrel{1^*\otimes Id}{\longrightarrow}NH(\DR{X}),
\end{multline*}
where $I:\tilde P_E\to\mathscr G_E$ is a natural inclusion, $m:B(\DR{X})\to NH(\DR{X})$ is the map, given by multiplication of the the leftmost and the rightmost tensors and normalization, and $1:X\to\tilde P_E$ is the unit section.
\end{lemma}
\begin{Proof}
The proof of this fact follows directly from the construction of map $\Psi^*$ described in \eqref{constrphie}, where on the right end of the diagramm we use the construction of section \ref{sectbegin}, compared with the definition of $\tilde\phi_E$.
\end{Proof}\\ \\ \vspace{1mm}

Let us now obtain a similar description of the map $\tilde{\hat\phi}$, determined by a twisting cochain $\phi$. To this end we first of all observe that the map $I^*$ is homotopy equivalent to $\tilde I:GP(\DR{X},\,\Omega_G)\to \Omega_G\hat\otimes_\phi\DR{X}$ from the diagramm \eqref{diagcomm}, and the map $1^*$ is equal to $\tilde\epsilon:\Omega_G\hat\otimes_\phi\DR{X}\to\DR{X}$, for $\epsilon$ -- the counit in $\Omega_G$ (this map extends to $:\Omega_G\hat\otimes_\phi\DR{X}$ because $\DR{X}$ is commutative.) It is now our purpose to find a map 
$$
\tilde\Psi^*:\Omega_G\hat\otimes_\phi\DR{X}\to GP(\DR{X},\,\Omega_G)\otimes_{\DR{X}^{ev}} \tilde B(\DR{X}),
$$
such that $\tilde{\hat\phi}=(\tilde\epsilon\otimes Id)(\tilde I\otimes m)\tilde\Psi^*$ (here $\tilde B(\DR{X})$ is the double-sided bar-resolution: $\tilde B(\DR{X})=\oplus_{n\ge0}\DR{X}\otimes\overline{\DR{X}}{}^{\otimes n}\otimes\DR{X}$.)

In order to find such a map, observe, that the complex $GP(\DR{X},\,\Omega_G)\otimes_{\DR{X}^{ev}} B(\DR{X})$ is quasi-isomorphic to $\Omega_G\hat\otimes_\phi\DR{X}$, the quasi-isomorphism being given by the map
\begin{align*}
\begin{split}
\tilde\Phi^*=Id\otimes_{\DR{X}^{ev}} Proj_0:GP(\DR{X},\,\Omega_G)&\otimes_{\DR{X}^{ev}}\tilde B(\DR{X})\to\\
                                                                              &\qquad\qquad \to\Omega_G\hat\otimes_\phi\DR{X},
\end{split}\\ 
\intertext{where}
Proj_0(a[a_1|\dots|a_n]b)=\begin{cases} 0,\ &n\ge1\\
                                                                    ab,\ &n=0.
                                                                    \end{cases}&
\end{align*}
In order to prove that $\tilde\Phi^*$ is quasi-isomorphism, one can employ the standard spectral sequence technique. Here we shall describe an explicit homotopy inverse map for $\tilde\Phi^*$; it is this map, that shall play the r\^ole of $\tilde\Psi^*$. So, we put:
\begin{equation}
\label{psitild}
\tilde\Psi^*(k\otimes a)=k\otimes (1[\,]a)+\sum_{n=1}^\infty k^{(1)}\otimes (1[\phi(k^{(2)})|\dots|\phi(k^{(n+1)})]a),
\end{equation}
where $k\otimes a$ is an arbitrary element in $\Omega_G\hat\otimes_\phi\DR{X},\ k\in\Omega_G,\ a\in\DR{X}$, and the tensor product on the right hand side is taken over $\DR{X}^{ev}$ (and we omit the left and right $\DR{X}$ tensors in $GP(\DR{X},\,\Omega_G)$.) It is easy to show that the map $\tilde\Psi^*$, given by this formula, commutes with the differential. Indeed, it is only necessary to chek that the map $\tilde\Psi^*$ commutes with the ``left'' part of the twisting, i.e.:
\begin{equation}
\label{eqcheckkk}
\tilde\Psi^*(dk\otimes a+(-1)^{|k|}k\otimes da+(-1)^{|k^{(1)}|}k^{(1)}\otimes \phi(k^{(2)})a)=(d\otimes Id+Id\otimes b+1\cap\phi)\tilde\Psi^*(k\otimes a),
\end{equation}
where $b$ is the differential in the double-sided bar-resolution $\tilde B(\DR{X})$ (we assume that the usual sign convention holds for the tensor product of differentials) and
$$
1\cap\phi(k\otimes (a1[a_2|\dots|a_n]a_{n+1}))=(-1)^{|k^{(1)}|}k^{(1)}\otimes(\phi(k^{(2)})a_1[a_2|\dots|a_n]a_{n+1}).
$$
So, we compute, step by step. On the left hand side, we have (we write down only the first three terms from the formula \eqref{psitild} and denote the argument of the left hand side of \eqref{eqcheckkk} by $L(k\otimes a)$):
\begin{equation}
\label{calcul}
\begin{split}
\tilde\Psi^*(L(k\otimes a))&=dk\otimes (1[\,]a)+(-1)^{|k|}k\otimes (1[\,]da)+(-1)^{|k^{(1)}|}k^{(1)}\otimes(1[\,]\phi(k^{(2)})a)\\
                                              &\quad+dk^{(1)}\otimes(1[\phi(k^{(2)})]a)+(-1)^{|k^{(1)}|}k^{(1)}\otimes(1[\phi(dk^{(2)})]a)\\
                                              &\qquad\quad+(-1)^{|k|}k^{(1)}\otimes (1[\phi(k^{(2)})]da)+(-1)^{|k^{(1)}|}k^{(1)}\otimes(1[\phi(k^{(2)})]\phi(k^{(3)})a)\\
                                              &\quad+dk^{(1)}\otimes(1[\phi(k^{(2)})|\phi(k^{(3)})]a)+(-1)^{|k^{(1)}|}k^{(1)}\otimes(1[\phi(dk^{(2)})|\phi(k^{(3)})]a)\\
                                              &\qquad\quad+(-1)^{|k^{(1)}|+|k^{(2)}|}k^{(1)}\otimes(1[\phi(k^{(2)})|\phi(dk^{(3)})]a)\\
                                              &\qquad\qquad\quad+(-1)^{|k|}k^{(1)}\otimes (1[\phi(k^{(2)})|\phi(k^{(3)})]da)\\
                                              &\qquad\qquad\qquad\quad+(-1)^{|k^{(1)}|}k^{(1)}\otimes(1[\phi(k^{(2)})|\phi(k^{(3)})]\phi(k^{(4)})a)+...
\end{split}
\end{equation}
On the other hand,
$$
\tilde\Psi^*(k\otimes a)=k\otimes(1[\,]a)+k^{(1)}\otimes(1[\phi(k^{(2)})]a)+k^{(1)}\otimes(1[\phi(k^{(2)})|\phi(k^{(3)})]a)+....
$$
If we apply $(d\otimes Id+Id\otimes b+1\cap\phi)$ to the first two terms, we obtain: 
\begin{equation*}
\begin{split}
&dk\otimes(1[\,]a)+(-1)^{|k|}k\otimes (1[\,]da)+(-1)^{|k^{(1)}|}k^{(1)}\otimes (\phi(k^{(2)})[\,]a)\\
&\qquad+dk^{(1)}\otimes(1[\phi(k^{(2)})]a)+(-1)^{|k^{(1)}|}k^{(1)}\otimes (1[d\phi(k^{(2)})]a)+(-1)^{|k|}k^{(1)}\otimes (1[\phi(k^{(2)})]da)\\
&\qquad+(-1)^{|k^{(1)}|+1}k^{(1)}\otimes (\phi(k^{(2)})[\,]a)+(-1)^{|k^{(1)}|}k^{(1)}\otimes (1[\,]\phi(k^{(2)})a)\\
&\qquad+(-1)^{|k^{(1)}|}k^{(1)}\otimes (\phi(k^{(2)})[\phi(k^{(3)})]a)
\end{split}
\end{equation*}
The first, the second, the fourth, the sixth and the seventh terms in this formula coincide with the first, the second, the fourth, the sixth and the third terms of the formula \eqref{calcul} respectively, the third and the sixth terms of this formula cancel each other, and the fifth term of this formula is equal to $(-1)^{|k^{(1)}|}(k^{(1)}\otimes (1[\phi(dk^{(2)})]a)+(-1)^{|k^{(1)}|+1}(k^{(1)}\otimes (1[\phi(k^{(2)})\phi(k^{(3)})]a)$, since $d\phi+\phi\cup\phi=0$. The first term here cancells with the fifth term in \eqref{calcul}, and the second one will cancell with a term, coming from $(Id\otimes b)(k^{(1)}\otimes(1[\phi(k^{(2)})|\phi(k^{(3)})]a))$. Further calculations will give the same results (in effect, this reasoning coincides with the proof that the map $\tilde\phi$ commutes with differentials.)

Thus, $\tilde\Psi^*$ is a chain map. One can find a chain homotopy $\tilde\Upsilon$, such that $\tilde\Psi^*\tilde\Phi^*-Id=d\tilde\Upsilon-\tilde\Upsilon d$, but this is not quite necessary, since we know, that $\tilde\Phi^*$ induces an isomorphism in homology, and on the other hand, it is evident that $\tilde\Phi^*\tilde\Psi^*=Id$.

The reason why we didn't pay much attention to the proof of the fact that $\tilde\Psi^*$ commutes with differentials is that in effect, one can obtain the map $\tilde\Psi^*$ and the corresponding chain homotopy with the help of the Perturbation Lemma, mentioned in previous sections. To this end, consider a ``smaller'' differential on $GP(\DR{X},\,\Omega_G)\otimes_{\DR{X}^{ev}} \tilde B(\DR{X})$ given by the sum of the usual differential in the tensor product and the additional part, given by the twisting cochain. Now, there's an evident homotopy equivalence between $\Omega_G\otimes\DR{X}$ and $GP(\DR{X},\,\Omega_G)\otimes_{\DR{X}^{ev}} \tilde B(\DR{X})$ if we throw away the twisting cochain from the differential on the left hand side too. Now, the additional terms, containig $\phi$ can be regarded as a small perturbation of the differential and the formulas, similar to \eqref{perteqn} now give the desired homotopy inverse $\tilde\Psi^*$ and $\tilde\Upsilon$.

It is evident, that the map $\tilde\Psi^*$ verifies the equation $\tilde{\hat\phi}=(\tilde\epsilon\otimes Id)(\tilde I\otimes m)\tilde\Psi^*$. Now the theorem follows from the next homotopy-commutative diagramm:
\begin{equation*}
\xymatrix{
{\Omega_G\hat\otimes_\phi\DR{X}} \ar[r]^{\mathscr J} \ar@/^1pc/[dd]^{\tilde\Psi^*} & {\DR{\tilde P_E}} \ar@/_1pc/[d]_{\Psi^*} \\
& {\DR{\mathscr P(X;\,\hat P\times_G P)}} \ar[d]^{\sigma^{-1}} \ar@/_1pc/[u]_{\Phi^*}\\
{GP(\Omega_G,\,\DR{X})\otimes_{\DR{X}^{ev}}\tilde B(\DR{X})}\ar@/^1pc/[uu]^{\tilde\Phi^*} \ar[r]^(.55){\mathscr I\otimes Id} \ar[d]_{\tilde I\otimes m} & {\DR{\mathscr G_E}\otimes_{\DR{X}^{ev}}\tilde B(\DR{X})} \ar[d]^{I^*\otimes m}\\
{(\Omega_G\hat\otimes_\phi\DR{X})\otimes_{\DR{X}}NH(\DR{X})} \ar[d]_{\tilde\epsilon\otimes Id} \ar[r]^(.54){\mathscr J\otimes Id} & {\DR{\tilde P_E}\otimes_{\DR{X}}NH(\DR{X})} \ar[d]^{1^*\otimes Id}\\
 {NH(\DR{X})}\ar@{=}[r] & {NH(\DR{X}).}
}
\end{equation*}
Here $\mathscr I,\ \mathscr J$ are the homotopy equivalences from the diagramm \eqref{diagcomm}, and $\mathscr J$ is besides this the homotopy equivalence from the proposition \ref{propleftright}. We have the following sequence of equalities up to homotopy:
$$
\tilde\phi_E\mathscr J=(1^*\otimes Id)(I^*\otimes m)\sigma^{-1}\Psi^*\mathscr J=(1^*\otimes Id)(I^*\otimes m)(\mathscr I\otimes Id)\tilde\Psi^*=(\tilde\epsilon\otimes Id)(\tilde I\otimes m)\tilde\Psi^*=\tilde{\hat\phi}.
$$
Observe that we have in fact showed, that $\tilde\phi_E$ and $\tilde{\hat\phi}$ are homotopic as homomorphisms of coalgebras over $\DR{X}$. Indeed, all the maps, involved here commute (up to a homotopy) with the coproducts.
\end{Proof}\\

\subsection{Equivariantization and Bismut classes}
\label{sect44}
In this last (but not the least) section we endeavor to explain in what way one can reproduce the cyclic Chern character of Bismut $ch(\nabla)$ (or, rather its Getzler-Jones-Petrack version $ch(p,\,\mathscr A)$) in the framework of our theory. Recall that this class belongs to the $\mathbb C[[u^{-1},\,u]$\/-linear cyclic cohomology of the equivariant complex $\DR{X\times \mathbb T}^{\mathbb T}[[u^{-1},\,u]$. The latter algebra is the $\mathbb T$\/-equivariant complex of the cartesian product $X\times\mathbb T$, where $\mathbb T$ denotes the unit circle. In effect, usual $\mathbb T$\/-equivariant complex should be an algebra over the polynomial ring $\mathbb C[u]$, but in our case we consider the tensor product over $\mathbb C[u]$ of this complex with the formal Laurent series ring. The condition of $\mathbb C[[u^{-1},\,u]$\/-linearity of the cyclic complex means that in making it one should consider all the tensor products over this ring. 

In this section we shall give a description of the cyclic Chern character class in terms of the Getzler-Jones map (see previos section) and discuss an analogous construction that involves twisting cochains. For any twisting cochain $\phi$ with domain $K$, first of all for $K=\Omega_{poly}{GL(n)}$, these constructions will give a series of cocycles in normalized cyclic complex of the algebra $A_{\mathbb T}$, where $A$ is the range of $\phi$. Here $A_{\mathbb T}$ is the equivariantized version of $A$, see description below. Similarly to the discussion of the paper \cite{GJP}, one can prove that the homology of the normalized cyclic complex is equal to the equivariant cohomology of the free loop space. However, these classes depend on the choice of $\phi$, and in particular on $K$ and $A$. As for its dependence on $\phi$ one can hope to get rid of it in the manner similar to Zamboni's, \cite{Zamboni}, proof of the independence of Getzler-Jones-Petrack's class of the choice of projector $p$, see its description below. We shall discuss this matter in a separate paper. However it is necessary to observe that in the case we consider here, we are faced with another sort of non-invariance, that is of the dependence on the coalgebra $K$. It would be nice to have a more invariant description of these classes, i.e. a description, independent of the choice of model coalgebra. But this should be done in the more general context of homotopy theory of maps between algebras and coalgebras and their invariants, which we postpone for future.

So let us first of all recall the definition of Getzler-Jones-Petrack's class $ch(p,\,\mathscr A)$. By definition it is the $0$\/-degree class in cyclic homology of complex $\DR{X\times \mathbb T}^{\mathbb T}[[u^{-1},\,u]$ (see its description below) determined by the following closed chain:
$$
ch(p,\,\mathscr A)=\sum_{n} Tr(p,\,\mathscr A,\dots,\mathscr A),
$$
where for an $n+1$\/-tuple of $N\times N$ matrices $A^0,A^1,\dots,A^1$ (entries of these matrices can have arbitrary nature) we put
$$
Tr(A^0,A^1,\dots,A^n)=\sum_{i_0,\dots,i_n=1}^N A^0_{i_0i_1}\otimes A^1_{i_1i_2}\otimes\dots\otimes A^n_{i_ni_0}.
$$
In our case we take $p$ to be just the projector-valued function on $X$, which determines $E$, and $\mathscr A$ a flat $\mathbb T$\/-equivariant connection on $X\times\mathbb T$, given by the formula:
$$
\mathscr A=A+\frac{Fdt}{u}.
$$
Here $A=pdp+(1-p)d(1-p)$ is the connection on $E\subseteq\underline N$, determined by $p$, and $F=pdpdp$ is its curvature. Both $A$ and $F$ should be regarded as matrix-valued differential forms on $X$. One can easily check, that the following equations hold (see \cite{GJP}, \S6)
\begin{align*}
d_up+[\mathscr A,\,p]&=0,\\
d_u\mathscr A+\frac12[\mathscr A,\,\mathscr A]&=0.
\end{align*}
The differential $d_u$ is the $\mathbb T^1$\/-equivariant differential on $\DR{X\times\mathbb T}^{\mathbb T}\otimes_{\mathbb C[u]}\mathbb C[u,\,u^{-1}]]\cong \DR{X}[1,\,dt]\otimes\mathbb C[u,\,u^{-1}]]$ (this algebra is denoted by $\Omega_{\mathbb T}(X\times\mathbb T)$ in the cited paper; for the sake of brevity we shall stick to the same notation below):
$$
d_u(a\, p(u,\,u^{-1}))=da\,p(u,\,u^{-1}),\qquad d_u(bdt\,q(u,\,u^{-1}))=dbdt\,\,q(u,\,u^{-1})+(-1)^{|b|}b\,uq(u,\,u^{-1}).
$$
It follows from these equations that $ch(p,\,\mathscr A)$ is a cocycle in the normalized cyclic complex $\tilde N(\Omega_{DR}^{X})$, where we once again use the notation of \S6 of the article \cite{GJP}:
\begin{align*}
\tilde N(\Omega_{DR}^{X})&=\tilde C(\DR{X})/\tilde D(\DR{X})\\
\intertext{where}
\tilde C(\DR{X})&=C_{\Omega_{\mathbb T}}(\Omega_{\mathbb T}(X\times\mathbb T))\widehat\otimes_{\Omega_{\mathbb T}}\mathbb C[[u^{-1},\,u]
\end{align*}
is the cyclic complex over the algebra $\Omega_{\mathbb T}=\mathbb C[u]$. Note that as a matter of fact it doesn't matter, whether we take tensor products of $\Omega_{\mathbb T}(X\times\mathbb T)$ over $\mathbb C[u]$ or of $\Omega_{\mathbb T}(X\times\mathbb T)\otimes_{\mathbb C[u]}\mathbb C[[u^{-1},\,u]$ over $\mathbb C[[u^{-1},\,u]$. Besides this $\tilde D(\DR{X})$ is the subcomplex generated by cyclic operations applied to degenerate cohchains in $\tilde C(\DR{X})$, i.e. to cochains, containing an element of $C^\infty(X)$ as one of their tensor legs. Our purpose now is to obtain $ch(p,\,\mathscr A)$ as an image of some element under the map, similar to $\tilde\phi_E$. To this end we need to generalize first $\tilde\phi_E$ so that its domain and range become $\mathbb T$\/-equivariant complexes and then pass to the cyclic complex on the right. Take notice that in our case we shall write the non-normalized tensor leg in the definition of $\tilde N(\Omega_{DR}^X)$ in the right-most leg of the tensor product, rather than in the left-most, as it is done in the paper \cite{GJP}. Of course, this doesn't affect the final result.

So let $A=pdp+(1-p)d(1-p),\ F=pdpdp$ be the grassmanian connection on $E$ and its curvature. We regard them as matrix-valued differential forms on $X$. Let $E_{ij}$ denote the matrix units, a basis in $Mat_N(\mathbb C)$ (so that $A=\sum_{i,j}A^{ij}\otimes E_{ij},\ F=\sum_{i,j} F^{ij}\otimes E_{ij}$,) which we regard as the Lie algebra $\mathfrak{gl}_N$, and let $e_{ij},\ \iota_{ij}$ be the Lie derivatives and contractions by the left-invariant vector fields on $GL(N)$, corresponding to $E_{ij}$. Now recall, that $\tilde\phi_E$ is generated by the derivative $\phi_E=\sum_{i,j}(A^{ij}\otimes e_{ij}+F^{ij}\otimes\iota_{ij})$ (it acts on the elements of $\DR{X}\otimes\DR{GL(N)}$ and preserves the ideal, which determines $\tilde P_E$.)

Observe that one can extend $\phi_E$ to a derivative on the equivariant complex $\Omega_{\mathbb T}(X\times\mathbb T)$: we just put $\phi_E(u)=\phi_E(dt)=0$, i.e. we extand $\phi_E$ to a $\mathbb C[dt,\,u,\,u^{-1}]]$\/-linear derivative $\Omega_{\mathbb T}(X\times\mathbb T)$. One can easily check that this extension preserves the usual equations, verified by the twisting maps (observe that $\phi_E$ is a twisting map in the sense of the section \ref{sectnext}, when we regard it as a derivative on $\DR{GL(N)}$, and hence the same equations are preserved by its restriction to $\DR{\tilde P_E}$,) i.e.
$$
[\phi_E^u,\,d_u]+(1\otimes\phi_E^u)\phi_E^u=0.
$$
This means that we can extend the usual formula for $\tilde\phi_E$ to obtain a map $\tilde\phi_E^u:\Omega_{\mathbb T}(\tilde P_E\times\mathbb T)\to NH_{\Omega_{\mathbb T}} (\Omega_{\mathbb T}(X\times\mathbb T))$ (the index $\Omega_{\mathbb T}$ on the right means that we consider the $\mathbb C[u]$\/-linear version of the Hochschild complex.) We put:
$$
\tilde\phi_E^u=\sum_n(\phi_E^u)^{\otimes n}|_{1}
$$
(here $1$ is the unit section of $\tilde P_E$.) On the other hand, one can consider another twisting map on the same complex $\DR{GL(N)}$ with values in the equivariant complex of $X\times\mathbb T$ and use it to induce a map on $\Omega_{\mathbb T}(X\times\mathbb T)$. Namely we put
$$
\psi_E^u=\sum_{i,j}(A^{ij}\otimes e_{ij}-u^{-1}F^{ij}dt\otimes e_{ij}).
$$
Since the equivariant curvature of the equivariant connection $\mathscr A$ vanishes, we conclude (see section \ref{sectnext}) that $\psi_E^u$ is a twisting map. Now we can restrict it to the algebra $\DR{\tilde P_E}\subseteq\DR{X}\otimes\DR{GL(N)}$ and then pass to the equivariant complex $\Omega_{\mathbb T}(\tilde P_E\times\mathbb T)$. Finally, we put:
$$
\tilde\psi_E^u=\sum_n(\psi_E^u)^{\otimes n}|_1.
$$
The following proposition is an important step towards the fuller understanding of Getzler-Jones-Petrack class.
\begin{prop}\label{proplbnl}$\ $\\
\vspace{-.7cm}
\begin{description}
\item{{\rm({\it i\/})}} The twisting maps $\phi_E^u$ and $\psi_E^u$ are gauge equivalent in the sence of section \ref{sectnnext};
\item{{\rm({\it ii\/})}} The image of class $ch(p,\,\mathscr A)$ in the Hochschild cohomology (under the natural inclusion of cyclic complex into the Hochschild complex) is equal to $\tilde\psi_E^u(tr_u)$, where $tr_u$ is certain equivariantly-closed element in $\Omega_{\mathbb T}(\tilde P_E\times\mathbb T)$.
\end{description}
\end{prop}
\begin{Proof}
First of all we shall give an explicit formula for the gauge transformation, relating  $\phi_E^u$ and  $\psi_E^u$. Put
$$
C=\sum_{i,j}u^{-1}F^{ij}dt\otimes\iota_{ij}.
$$
This is a degree $0$ differentiation on $\DR{GL(N)}$ with values in $\Omega_{\mathbb T}(X\times\mathbb T)$, and we can extend it to a differentiation on $\Omega_{\mathbb T}(\tilde P_E\times\mathbb T)$ with values in the same algebra, as described above. An easy computation shows that the following equation holds for the differentiations $\phi_E^u,\ \psi_E^u$ and $C$:
$$
\phi_E^u-[\phi_E^u,\,C]-dC=\psi_E^u,
$$
and hence if we put $c=\exp(C)=1+C$, then $c^{-1}=1-C$ and $\psi_E^u=c^{-1}\,\phi_E^u\, c+c^{-1}dc$.

Next, consider the function $Tr_n$ on $GL(n)$, equal to the sum of the diagonal elements of a matrix: $Tr_n=\sum u_{ii}$. This function is $Ad$\/-invariant, so we can extend it from $GL(n)$, which we regard as a fibre of $\tilde P_E$, to the whole $\tilde P_E$. Indeed, this is possible since one passes from one chart of the bundle $\tilde P_E$ to another by means of the adjoint action of the group on itself. We shall denote this extension by $tr_E$. Clearly this function is equal to the restriction of the function $Tr_N$, or rather of $1\otimes Tr$, on $X\times GL(N)$ on $\tilde P_E\subseteq X\times GL(N)$, the latter inclusion being induced by the inclusion of $E$ into trivial bundle $X\times\mathbb C^N$. In other words, $tr_E\equiv1\otimes Tr_N(\mathrm{mod(I)}$, where $I$ is the ideal in $C^\infty(X)\otimes\ac(GL(N))$, generated by $1_N-p(x)$. Similarly, one can extend the differential of $Tr_n$, $dTr$, to a degree $1$ differential form $dtr_E$ on $\tilde P_E$, so that $d(tr_E)=dtr_E$, where on the left we have the usual de Rham differential on $\tilde P_E$. Moreover, this form is equal to the restriction of $dTr_N$ onto $\tilde P_E$. Then we put $tr_u=tr_E+u^{-1}dtr_E\,dt$. Clearly $tr_u$ is an equivariantly-closed form, $d_u(tr_u)=0$. In effect, $tr_u$ is equivariantly-exact: $tr_u=d_u(u^{-1}tr_E\,dt)$.

Now it is a matter of simple calculations to see that $\tilde\psi_E^u(tr_u)=ch(p,\,\mathscr A)$. The only thing one should use here is that the Lie algebra of $GL(N)$ is a linear space spanned by the left-invariant vector fields $T_{ij}(g)=\sum_k u_{ki}(g)\dfrac{\partial}{\partial u_{kj}}$ on $GL_n(k)$, and hence
$$
\sum_{i,j}A_{ij}\otimes T_{ij}(u_{mn})=\sum_{i,j,k}A_{ij}\otimes u_{ki}\dfrac{\partial}{\partial
u_{kj}}(u_{mn})=\sum_{i,j,k}A_{ij}\otimes u_{ki}\delta_{km}\delta_{nj}=\sum_iA_{in}\otimes u_{mi}.
$$
Besides this, one should take into consideration, that Lie derivative commutes with the de Rham differential and that the restriction of $dtr_E$ on the unit section is equal to $0$ (because this form has nontrivial entries only in the vertical direction.)
\end{Proof}\\ \\ \vspace{1mm}
In fact, from the observations we made, it follows that the element $ch(p,\,\mathscr A)$ is exact as the element in the Hochschild complex. Namely, $ch(p,\,\mathscr A)=b(\tilde\psi_E^u(u^{-1}tr_E\,dt))$. However, it is not exact in cyclic complex: an easy computation shows that $B(\tilde\psi_E^u(u^{-1}tr_E\,dt))$ is not a degenerate element, and hence $(b+uB)(\tilde\psi_E^u(u^{-1}tr_E\,dt))\ne ch(p,\,\mathscr A)$ in $\tilde N(\Omega_{DR}^{X})$. 

It would be tempting to assume that $\phi_E^u(tr_u)$ gives the same cyclic class. Then we could perform the same equivariantization procedure with the map $\tilde{\hat\phi}$ associated with a twisting cochain $\phi$ and obtain the same class, applying $\tilde{\hat\phi}_u$ to the element $tr_n\otimes1 +dtr_n\otimes u^{-1}dt\in \DR{GL(n)}\hat\otimes_\phi\Omega_{\mathbb T}(X\times\mathbb T)$, which is evidently equal to $tr_u$ under the equivalence $\mathscr I$. Unfortuantely, the element $\tilde\phi_E^u(tr_u)$ is in fact not even closed with respect to the cyclic differential $B$\ -- its first tensor leg is proportional to $u^{-1}dt$. So, we cannot regard $\tilde\phi_E^u(tr_u)$ as the image of some cyclic cocycle with respect to the natural embedding of the cyclic complex into the Hochschild one, and one should change the construction of $\tilde{\hat\phi}_u$ a little bit. In fact, one can amend this situation, using the following idea.

Let $\phi:K\to A$ be a twisting cochain, where $K$ is a differential graded coalgebra and $A$\ -- a differential algebra. Consider the formal equivariantization of $A$: $A_{\mathbb T}=A\otimes(\mathbb C[1,dt]\otimes\mathbb C[u,u^{-1}]])$ with a $\Omega_T$\-linear differential $d_u$:
$$
d_u(a'+a''dt)=da'+da''dt+(-1)^{|a''|}ua''.
$$
We can now consider the twisted tensor product of $K$ and $A_{\mathbb T}$, where $A_{\mathbb T}$ is considered as the $A$\-module via the natural embedding $A\to A_{\mathbb T}$. Let $K_o$ denote the coalgebra $K$ regarded as coalgebra with trivial differential. Then the following proposition holds:
\begin{prop}
The twisted tensor product $K\otimes_\phi A_{\mathbb T}$ is isomorphic to the following chain complex: $K_o\otimes_{\phi_o} A_{\mathbb T}$, where $\phi_o$ is the twisting cochain $\phi_o:K_o\to A_{\mathbb T}$, given by the formula:
$$
\phi_o(k)=\phi(k)+u^{-1}\phi(dk)dt.
$$
\end{prop}
\begin{Proof}
Consider the map $O:K\otimes A_{\mathbb T}\to K\otimes A_{\mathbb T}$:
$$
O(k\otimes a)=k\otimes a-dk\otimes u^{-1}adt.
$$
This map is invertible, its inverse is given by
$$
O^{-1}(k\otimes a)=k\otimes a+dk\otimes u^{-1}adt.
$$
Using this map, we can deform the differential $d_\phi$, replace it by $O^{-1}\circ d_\phi\circ O$. An easy calculation now shows that the latter differential is equal to $(-1)^{|k|}(k\otimes da+k^{(1)}\otimes \phi(k^{(2)})a+k^{(1)}\otimes u^{-1}\phi(k^{(2)})dt\,a)=d_{\phi_o}(k\otimes a)$ on arbitrary element $k\otimes a$.
\end{Proof}\\ \\ \vspace{1mm}
Similarly, one can deform the differential in $K\hat\otimes A$ so that it become isomorphic to $K_o\hat\otimes_{\phi_o}A$. On the other hand, since the differential in $K_o$ vanishes, we conclude that an element $k\otimes 1\in K_o\hat\otimes_{\phi_o}A$ is closed if and only if it verifies the equation
$$
k^{(1)}\otimes \phi_o(k^{(2)})=k^{(2)}\otimes \phi_o(k^{(1)}).
$$
In particular, this equation holds, if $k$ is in the cocenter of $K$ (i.e. if $k^{(1)}\otimes k^{(2)}=k^{(2)}\otimes k^{(1)}$,) thus $k\otimes 1$ is closed in $K_o\hat\otimes_{\phi_o} A$ for all $k$ in cocenter. For instance, if $K=\DR{GL(n)}$ one can take $k=\sum u_{ii}=tr_n$.

Using the standard construction, discussed in this paper, one obtains (for arbitrary twisting cochain $\phi$ on $\DR{GL(n)}$) a closed element $\tilde{\hat\phi}_o(tr_n\otimes 1)$ in $\tilde N(A_{\mathbb T})$, where the normalized cyclic complex $NH(A_{\mathbb T})$ is given by the same construction as $\tilde N(\Omega_{DR}^{X})$ before. Repeating the reasonong of the \S5 \cite{GJP} one concludes that the cohomology of $\tilde N(A_{\mathbb T})$ is equal to the equivariant cohomology of the free loop space of $X$ as soon as $A$ is quasi-equivalent to $\DR{X}$. It is now obvious, that this elemnt is in effect closed in reduced complex with respect to the cyclic operator $B$ (in effect, its last tensor leg is equal to $1\in A_{\mathbb T}$.) We shall denote the element, determined by this construction by $ch(\phi,\,k)$ (here $\phi$ is a twisting cochain and $k\in K$ is a cocentral element.) In particular, we obtain the class $ch(\phi_P,\,tr_n)$, where $\phi_P$ is the cochain, defined in previous chapter, given by the formula
\begin{equation}
\begin{split}
ch(\phi_P,\,tr_n)&=\sum_{n=0}^\infty \sum_{i_1,\dots,i_n}\{\delta^{i_1i_2}-g_{\alpha\beta}^{i_1i_2}+u^{-1}dg_{\alpha\beta}^{i_1i_2}dt\}\otimes\\
                             &\qquad\qquad\quad\{\delta^{i_2i_3}-g_{\alpha\beta}^{i_2i_3}+u^{-1}dg_{\alpha\beta}^{i_2i_3}dt\}\otimes\dots\otimes \{\delta^{i_ni_1}-g_{\alpha\beta}^{i_ni_1}+u^{-1}dg_{\alpha\beta}^{i_ni_1}dt\}\otimes 1.
\end{split}
\end{equation}
Here $A=\check C_\uc^\ast(X,\,\DR{U})$, $g_{\alpha\beta}:U_{\alpha\beta}\to GL(n)$ are the transition functions of the principal bundle $P$ and $\delta^{ij}$ is the Kronecker symbol. This formula defines a class in the cyclic cohomology of $A_{\mathbb T}$, which is isomorphic to the cohomology of free loop space of $X$ (we assume that $X$ is 1-connected.)

Observe that in general the class of an element $ch(\phi,\,k)$ can depend on the choice of the twisting cochain $\phi$. Indeed, if $\phi'$ is anothe cochain, equivalent to $\phi$ via a transformation map $c:K\to A$, then one can extend $c$ to an equivalence $c_o$ between the deformed cochains $\phi_o$ and $\phi'_o$. Just put $c_o(k)=c(k)+u^{-1}c(dk)dt$. Then the formula \eqref{hombad} gives a connecting homotopy that allows one compare the elements $\tilde{\hat\phi}_o(k\otimes1)$ and $\tilde{\hat\phi}'_o(k\otimes 1)$. As one can see, this formula does not necessarily respect the cyclic operator $B$, i.e. $B(\hat H_{c_o}(k\otimes 1))$ can be different from $0$, since it can contain such elements elements as $c(dk)dt$ at the last tensor leg. Thus the classes we obtain depend on $\phi$ pretty much as the Chern-Symons clases depend on the choice of flat connections. One can try to prove their independence on the choices made by a more delicate reasoning, similar to that of Zamboni \cite{Zamboni}. On the other hand, if we chose the twisitng cochain $\xi$ from the section \ref{sectconnect}, then we evidently obtain a class, that will be equal to Bismut's class after globalization (since they coincide on every local plot in $\mathcal LX$ in the sense of Chen.)


\begin{thebibliography}{99}
\bibitem{GJP}{E.~Getzler, J.D.S.~Jones, S.~Petrack. \emph{Differential forms on loop space and the cyclic bar complex}, Topology {\bf 30} (1991), 339-371}
\bibitem{Br59}{E.~Brown. \emph{Twisted tensor products}, Ann. of Math, {\bf 69} (1959), 223-246}
\bibitem{Smirn76}{В.А.~Смирнов. \emph{Функтор $D$ для скрещенных тензорных произведений}, Мат.~Заметки {\bf 20} вып.~4, (1976), 465-472, in Russian}
\bibitem{Smirn77}{В.А.~Смирнов. \emph{Функтор $D$ и сильная гомотопия}, Мат.~Заметки {\bf 21} вып.~4, (1977), 557-564, in Russain}
\bibitem{SmSOp}{V.A.~Smirnov. \emph{Simplicial and Operad Methods in Algebraic Topology}, AMS, Transl. Math. Monogr. 198, 2001}
\bibitem{mypapLF}{Г.И.~Шарыгин. \emph{Локальные формулы для характеристических классов главного $GL_n$\/-расслоения}, Мат.~Сборник, {\bf 190} вып.~10 (2008), 127-158, in Russian}
\bibitem{BottTu}{R.~Bott, L.W.~Tu. \emph{Differential Forms in Algebraic Topology}, Springer-Verlag, Berlin-Heiderberg-New York, 1982}
\bibitem{Kan}{D.~Kan. \emph{A combinatorial definition of homotopy groups}, Ann. of Math. {\bf 67} (1958), 282-312}
\bibitem{Chen}{K.T.~Chen. \emph{Iterated integrals of differential forms and loop space homology}, Ann. Math. {\bf 97}(1973), 217-246}
\bibitem{Chen2}{K.T.~Chen. \emph{Extension of $C^\infty$\/-function algebra by integrals and Malcev completion of $\pi_1$}, Adv.
Math. {\bf 23} (1977), 181-210}
\bibitem{Jones}{J.D.S.~Jones. \emph{Cyclic homology and equivariant homology}, Invent. Math. {\bf 87} (1987), 403-423}
\bibitem{Shih}{W.T.~Shih. \emph{Homology des espaces fibr\'es}, Publications Math. de l'IHES {\bf 13} (1962), 93-176}
\bibitem{McL}{S.~McLane. \emph{Homology}, Springer-Verlag, Berlin-Heidelberg-New York, 1963}
\bibitem{BottSeg}{R.~Bott, G.~Segal. \emph{The Cohomology of the Vector Fields on a Manifold}, Topology, {\bf 16} (1977), 285-298}
\bibitem{Anders}{D.W.~Anderson. \emph{A generalization of the Eilenberg-Moore spectral sequence}, Bulletin of the Amer. Math. Soc. {\bf 78} Number 5 (1972), 784-786}
\bibitem{Onisch}{А.Л.~Онищик. \emph{Топология транзитивных групп преобразований}, Москва, Физ. мат. литература, 1995, in Russian}
\bibitem{Sweedler}{M.E.~Sweedler. \emph{Hopf algebras}, W. A. Benjamin, New-York, 1969}
\bibitem{Dupont}{Dupont~J.-L. Simplicial de Rham cohomology and characteristic classes of flat bundles. Topology {\bf 15} (1976), 233-245}
\bibitem{MaySO}{J.P.~May. \emph{Simplicial methods in algebraic topology}, U.Chicago, 1967}
\bibitem{Segal74}{G.~Segal. \emph{Categories and cohomology theories}, Topology {\bf 13} (1974), 293-312}
\bibitem{Zamboni}{L.Q.~Zamboni. \emph{A Chern character in cyclic homology}, Trans. of AMS {\bf 331} (1992), 157-163}
\end{thebibliography}
\end{document}